\documentclass[11pt,a4paper]{amsart}
\usepackage{amsfonts,amsmath,amssymb}
\usepackage{hyperref}
\usepackage[capitalise]{cleveref}
\usepackage{latexsym,fullpage,amsfonts,amssymb,amsmath,amscd,graphics}
\usepackage[all]{xy}
\usepackage{amssymb,amsthm,amsxtra}
\usepackage[usenames]{color}
\usepackage{amscd}
\usepackage{amsthm}
\usepackage{amsfonts}
\usepackage{amssymb}
\usepackage{mathrsfs,comment}
\usepackage{url}
\usepackage{bbm}
\usepackage{wasysym}
\usepackage{enumitem}
\usepackage{framed}
\usepackage{mathtools}
\usepackage{booktabs,array,longtable}

\usepackage[style=alphabetic, backend=biber, url=false, doi=false]{biblatex}
\usepackage{epigraph}
\usepackage{tikz}
\usetikzlibrary{decorations.markings}
\usetikzlibrary{cd}
\usetikzlibrary{patterns}
\usetikzlibrary{shapes}

\usetikzlibrary{positioning,arrows.meta,calc}

\tikzset{anchorbase/.style={baseline={([yshift=-0.5ex]current bounding 
box.center)}}}
\tikzset{wipe/.style={white,line width=4pt}}

\tikzset{->-/.style={decoration={
  markings,
  mark=at position #1 with {\arrow{>}}},postaction={decorate}}}
\tikzset{-<-/.style={decoration={
  markings,
  mark=at position #1 with {\arrow{<}}},postaction={decorate}}}

\tikzset{darkg/.style={green!70!black}}

\tikzset{->-/.style={decoration={
  markings,
  mark=at position #1 with {\arrow{>}}},postaction={decorate}}}
\tikzset{-<-/.style={decoration={
  markings,
  mark=at position #1 with {\arrow{<}}},postaction={decorate}}}

\theoremstyle{plain}% default
\newtheorem*{theorem*}{Theorem}
\newtheorem*{remark*}{Remark}
\newtheorem*{example*}{Example}
\newtheorem{lemma}{Lemma}[subsection]
\newtheorem{proposition}[lemma]{Proposition}
\newtheorem{corollary}[lemma]{Corollary}
\newtheorem{theorem}[lemma]{Theorem}

\newtheorem*{conjecture*}{Conjecture}

\newtheorem{introtheorem}{Theorem}
\newtheorem{introcorollary}[introtheorem]{Corollary}

\theoremstyle{definition}
\newtheorem{definition}[lemma]{Definition}

\newtheorem{example}[lemma]{Example}

\theoremstyle{remark}
\newtheorem{remark}[lemma]{Remark}

\newtheorem{notation}[lemma]{Notation}

\newtheorem{question}[lemma]{Question}

\newcommand{\Tr}{\operatorname{Tr}}
\newcommand{\Hom}{\operatorname{Hom}}

\newcommand{\eps}{\varepsilon}
\newcommand{\triv}{{\mathbbm 1}}

\newcommand{\Ind}{\operatorname{Ind}}
\newcommand{\Res}{\operatorname{Res}}

\newcommand{\id}{\operatorname{id}}

\newcommand{\Ker}{\operatorname{Ker}}

\newcommand{\Span}{\mathrm{span}}
\newcommand{\Aut}{{\operatorname{Aut}}}

\newcommand{\rk}{{\operatorname{rk}}}
\newcommand{\End}{\operatorname{End}}

\newcommand{\Z}{{\mathbb Z}}

\newcommand{\bR}{{\mathbb R}}

\newcommand{\V}{\mathbf{V}}

\newcommand{\uRep}{\underline{{\rm Rep}}}
\newcommand{\Perm}{\underline{{\rm Perm}}}
\newcommand{\Rep}{\operatorname{Rep}}

\newcommand{\F}{{\mathbb F}}

\newcommand{\kk}{\mathbbm{k}}

\newcommand{\Vect}{\mathrm{Vect}_{\kk}}
\newcommand{\one}{{\mathbf{1}}}

\newcommand{\T}{\mathcal{T}}

\newcommand{\TID}{\mathcal{TID}}
\newcommand{\ID}{\mathcal{ID}}
\newcommand{\Kar}{\mathrm{Kar}}

\newcommand{\Down}{\mathcal{D}^{\downarrow}}

\newcommand{\Emb}{\operatorname{Emb}}
\newcommand{\rad}{\operatorname{rad}}

\newcommand{\Stab}{\operatorname{Stab}}
\newcommand{\I}{\mathcal{I}}

\newcommand{\N}{\mathcal N}
\newcommand{\ord}{\operatorname{ord}}
\newcommand{\boxt}{\mathbin{\boxtimes}}

\newcommand{\nd}{\mathrm{nd}}
\newcommand{\rmrad}{\mathrm{rad}}
\newcommand{\cC}{\mathcal{C}}

\newcommand{\Rel}{\mathrm{Rel}}
\newcommand{\Syl}{\operatorname{Syl}}
\newcommand{\Sub}{\mathbf{Sub}}
\newcommand{\cE}{\mathcal{E}}
\newcommand{\cP}{\mathcal{P}}
\newcommand{\mfrak}{\mathfrak{m}}
\newcommand{\Set}[1]{\mathbf{S}(#1)}
\newcommand{\longtwoheadrightarrow}{\relbar\joinrel\twoheadrightarrow}

\newcommand{\oG}{{\overline{G}}}
\newcommand{\oX}{{\overline{X}}}
\newcommand{\oU}{{\overline{U}}}
\newcommand{\For}{\operatorname{For}}

\def\quotient#1#2{%
    \raise1ex\hbox{$#1$}\Big/\lower1ex\hbox{$#2$}%
}

\begin{document}
\title{Tensor ideals of oligomorphic permutation categories}
\date{\today}

  \author{Inna Entova-Aizenbud, Thorsten Heidersdorf}
% \address{}

 \address{I. E.: Department of Mathematics, Ben Gurion University of the Negev, Beer-Sheva, Israel}
 \email{entova@bgu.ac.il}
 \address{T. H.: School of Mathematics, Statistics and Physics, Newcastle University, UK}
 \email{heidersdorf.thorsten@gmail.com}

\begin{abstract} We classify tensor ideals for oligomorphic permutation categories $\Perm(G,\mu)$. To this end we define an equivalence relation on the class of transitive smooth $G$-sets. We show that there is a bijection between the tensor ideals of $\Perm(G,\mu)$ and the downward-closed subsets of the poset of equivalence classes. As a corollary we obtain that there is a bijection between tensor ideals and thick tensor ideals. We determine the poset in many cases (mostly in characteristic $0$): the infinite symmetric group, wreath products, products of oligomorphic groups, equivariantizations, finite groups in positive characteristic and the classical groups $ GL, O, Sp, U$ over a finite field.
\end{abstract}

\maketitle

\epigraph{"A new Power is rising. Against it the old allies and policies will not avail us.... [...] This then is one choice before you, before us. We may join with that Power. [...] Its victory is at hand; and there will be rich reward for those that aided it. As the Power grows, its proved friends will also grow; and the Wise, such as you and I, may with patience come at last to direct its courses, to control it. [...] There need not be, there would not be, any real change in our designs, only in our means."}{\textit{The Lord of the Rings \cite{Tolkien}}}

\setcounter{tocdepth}{1}
\tableofcontents

\section{Introduction}

\subsection{Tensor ideals}
Let $\kk$ be a field. An important problem in any $\kk$-linear monoidal category $\mathcal{T}$ is to classify the tensor ideals (on morphisms) and the thick tensor ideals (on objects). These arise as kernels of monoidal functors and give important structural information about the given category. Let us denote the poset of tensor ideals by $\TID(\mathcal{T})$ and the poset of thick tensor ideals by $\ID(\mathcal{T})$. These are related by the decategorification map \[ Ob: \TID(\mathcal{T}) \to \ID(\mathcal{T}), \ \ \I \mapsto \I = \{ X \in \mathcal{T} \ | \ \id_X \in \I(X,X) \}.\]

Among tensor ideals in a $\kk$-linear rigid symmetric monoidal category $\mathcal{T}$ with $\End(\one)=\kk$, one has the trivial tensor ideals  (the zero ideal and the entire category) and the ideal $\mathcal{N}$ of negligible morphisms, the analogue of the radical of an algebra. The ideal $\mathcal{N}$ is the maximal proper tensor ideal, and it consists of all morphisms $f \in \Hom(X, Y)$ such that $\mathrm{Tr}(g \circ f) = 0$ for all $g \in \Hom(Y, X)$. The quotient of the Karoubi envelope (idempotent completion) $\mathcal{T}^{\Kar}$ by $\mathcal{N}$ is called the \emph{semisimplification of $\mathcal{T}$}, and it is a semisimple $\kk$-linear rigid symmetric monoidal category if it satisfies the additional condition that any nilpotent endomorphism has zero trace.

A famous example of the classification of tensor ideals is Ostrik's classification \cite{Ostrik} of the thick tensor ideals of $\mathrm{Tilt}(U_q(\mathfrak{g}))$, the category of tilting modules of {Lusztig}'s quantum group at a root of unity $q$, in terms of two-sided cells of the affine Weyl group of $\mathfrak{g}$. As shown in \cite{Coulembier-ideals}, the decategorification map is in general not a bijection in this case.

Even for diagrammatic monoidal categories, the problem of classifying tensor or thick tensor ideals is complicated and only a few cases are known. It was shown in \cite{CO-abelian} that there is only one proper non-trivial tensor ideal for $\uRep(S_t)$ if $t\in \mathbb{Z}_{\geq 0}$ in characteristic $0$. The proof uses the classification of blocks of $\uRep(S_t)$ in \cite{CO-blocks}. Other known cases include the classification of tensor ideals and thick tensor ideals for $\uRep(GL_t)$ and $\uRep(O_t)$ ($t \in \mathbb{Z}$) in \cite{Comes}\cite{Comes-Heidersdorf,Coulembier-ideals}. In this case the decategorification map is a bijection, and the proper non-trivial tensor ideals form a strictly descending infinite chain \[ \mathcal{N} = \mathcal{I}_1 \supset \I_2 \supset {\I_3} \ldots \] This classification is of fundamental importance in the construction of abelian envelopes for these categories \cite{EHS, Coulembier-abelian} and in applications to the representation theory of supergroups \cite{heidersdorf2017mixed}.

\subsection{Oligomorphic permutation categories} 
Oligomorphic groups come from permutation group theory and model theory, where they were introduced by P. Cameron \cite{Cameron-orbits}. A topological group $G$ is called (pro-)\emph{oligomorphic} if it satisfies three conditions: (i) the open subgroups form a base of {neighborhoods} of the identity (that is, $G$ is non-Archimedean), (ii) the topology is Hausdorff, (iii) for any open subgroups $U, V\leq G$, the set $U\backslash G/V$ is finite.

A typical source of oligomorphic groups is the {following} setting: let $G$ act faithfully on a set $\Omega$ so that $G$ has finitely many orbits on $\Omega^n$ for every $n$; the topology on $G$ is then determined by the base of {open} subgroups $\Stab(A)$, the pointwise stabilizers of finite subsets $A\subset \Omega$.

The natural class of interesting $G$-sets in this setting are the \emph{smooth finitary} $G$-sets: namely, $G$-sets $X$ which contain finitely many $G$-orbits, and such that $\Stab(x)$ is an open subgroup of $G$ for any $x\in X$. When considering the representation theory of an oligomorphic group $G$, it makes sense to first study permutation representations of $G$ on the spaces of functions of such smooth finitary $G$-sets.

In their seminal paper \cite{HS-oligomorphic}, Harman and Snowden proposed considering a measure $\mu$ on $G$ to construct the category $\Perm(G, \mu)$ of permutation representations of an oligomorphic group $G$. Such a measure is not a measure on $G$ in the analytic sense: instead, it is a collection of ``compatible'' locally constant measures on the finitary smooth $G$-sets (and more generally, on the finitary smooth $U$-sets where $U\leq G$ is an open subgroup). The measure takes values in a base field $\kk$ of our choosing.

Oligomorphic permutation categories $\Perm(G, \mu)$, as defined in \cite{HS-oligomorphic}, are a genuinely new way to construct rigid symmetric monoidal categories, and in some cases, pre-Tannakian categories. They incorporate many known cases of diagrammatic monoidal categories such as the Deligne category $\uRep(S_t)$ for $t \in \kk$ (where $\kk$ is a field) or more generally all of Knop's tensor envelope categories \cite{Knop, Snowden-Knop}. The latter includes in particular the interpolation categories for wreath products $\Gamma \wr S_n$ \cite{Nyobe-Savage} where $\Gamma$ is a finite group and the interpolation category of the finitary linear group $GL_n(\F_q)$ \cite{EAHGL}. 

However, $\Perm(G,\mu)$ gives many new examples. It was shown in \cite{HS-oligomorphic} that there are four different measures on the group $G = \Aut(\mathbb{R}, <)$ of order preserving automorphisms of the real line, leading to the definition of Delannoy categories, as well as its circular variant. Delannoy categories were further studied in \cite{HSS-delannoy, CS-delannoy, CHS-delannoy}. Other examples of oligomorphic groups are the automorphism group of the Rado graph and the {oligomorphic} group associated with arboreal trees (see \cite{HNS}).

\subsection{Classification of tensor ideals}

The world of oligomorphic groups is very rich and encompasses many different types of groups, and it might seem hopeless to classify the tensor ideals of $\Perm(G,\mu)$, given the difficulty one encounters already in a special case such as $\mathrm{GL}_{\infty}(\F_q)$. Somewhat to our surprise, there is a general classification theorem of tensor ideals in $\Perm(G,\mu)$ for any $G$, any $\mu$ and any ground field $\kk$.

For fixed measure $\mu$ we define an equivalence relation $\sim_{\mu}$ on transitive smooth $G$-sets and denote the poset of equivalence classes by $\mathcal{P}_{\mu}^G$. We say a subset $D \subseteq \mathcal{P}^G_\mu$ is \emph{downward-closed} if $[Y] \in D$ and $[Z] \preceq_\mu [Y]$ implies $[Z] \in D$. 

\begin{introtheorem} [\cref{thm:classification}]
    The poset of tensor ideals of $\Perm(G, \mu)$ is in bijection with the poset of downward-closed subsets of $\mathcal{P}^G_\mu$. A tensor ideal $\mathcal{I}\subset\Perm(G, \mu)$ corresponds to the subset $D\subseteq \mathcal{P}^G_\mu$ if {and} only if $\mathcal{I}$ is generated by $\id_{\mathcal{C}(X)}$ where $[X]\in D$.
\end{introtheorem}

As any tensor ideal is therefore generated by the identities it contains, we get:

\begin{introcorollary} [\cref{cor:tensor_ideals_are_thick}] The decategorification map \[ Ob: \TID(\Perm(G, \mu)) \to \ID(\Perm(G, \mu))\] is a bijection.
\end{introcorollary}

The crucial tool of the proof of the classification theorem is an abstract description of the tensor ideals in a general braided right-rigid $\kk$-linear preadditive category $\mathcal{T}$ in \cite{Coulembier-ideals}. Restriction to morphisms out of the tensor unit induces an isomorphism of posets
\[
 \mathcal{TID}(\mathcal T)
 \xrightarrow{\sim}
 \operatorname{Sub}_{\mathcal T}\Hom_{\mathcal T}(\one,-),
 \qquad
 \I\longmapsto \I(\one,-)
\]
where the right-hand side consists of the submodules of $\Hom(\one,-)$. As the proof of our main classification theorem shows, these subfunctors admit for $\Perm(G,\mu)$ an explicit description in terms of the poset $\mathcal{P}_{\mu}^G$.

An important class of measures studied in \cite{HS-oligomorphic} are \emph{regular measures}; an easy consequence of the above classification is that these are precisely the measures for which the categories $\Perm(G, \mu)$ have only trivial tensor ideals. It was shown in \cite{HS-oligomorphic, HS-classical} that when $\mu$ is regular, $\Perm(G, \mu)^{\Kar}$ is a semisimple tensor category if and only if any nilpotent morphism in $\Perm(G, \mu)$ has zero trace (such property of $\mu$ is called the (Nil) property in \cite{HS-classical}).

A wider class of measures defined in \cite{HS-oligomorphic} are \emph{quasi-regular measures}: these are measures which restrict to a regular measure on some open subgroup of $G$. We show in \cref{prop:quasiregular} that there always exists a minimal non-zero tensor ideal when $\mu$ is quasi-regular.

\subsection{Examples}

It turns out that the poset $\mathcal{P}_{\mu}^G$ can be calculated very efficiently, leading to a large class of $\Perm(G,\mu)$-categories with explicit list of tensor ideals.
The following easy fact makes all these calculations feasible: in characteristic zero, finite covers $f:X\to Y$ do not change the equivalence class. As a consequence, it is enough to compute the equivalence classes for a distinguished class of $G$-sets such that every other transitive $G$-set has a finite cover from this class. Such distinguished classes can be easily found in many examples such as for $G = S_{\infty}$, wreath products $\Gamma \wr S_{\infty}$ or the finitary classical groups.

Let {us} describe briefly two types of examples of oligomorphic groups which we will use below. 

The first is given {by} the infinite symmetric group $S_{\infty}$, or, more generally, by wreath products $\Gamma \wr S_{\infty}$ of a finite group $\Gamma$ with $S_{\infty}$. For a field $\kk$ of characteristic zero, the $\kk$-valued measures on each of these oligomorphic groups are parametrized by $t\in \kk$, and they are regular if and only if $t\notin \Z_{\geq 0}$ (see \cite{HS-oligomorphic}, \cref{sec:wreath-products_intro}).

Another nice source of examples are the classical oligomorphic groups\footnote{There are two interesting topologies on the group $\mathrm{GL_{\infty}}$, giving two distinct oligomorphic groups: the \emph{parabolic} topology and the \emph{Levi} topology, which we denote here by $\mathrm{GL^{par}_{\infty},~GL^{Levi}_{\infty}}$ respectively. As previously shown in \cite{HS-classical}, for $t\neq 0$ we have $\Perm(\mathrm{GL^{par}_{\infty}},\mu_t)^{\Kar}\cong \Perm(\mathrm{GL^{Levi}_{\infty}},\mu_t)^{\Kar}$.} $\mathrm{GL^{par}_{\infty},~GL^{Levi}_{\infty}, ~Sp_{\infty},~ U_{\infty},~ O_{\infty}}$ over a finite field $\mathbb{F}_q$, described explicitly in \cref{ssec:setup_classcial_grps}. For a field $\kk$ of characteristic zero, the $\kk$-valued measures on each of the oligomorphic groups $\mathrm{GL^{par}_{\infty},~GL^{Levi}_{\infty}, ~Sp_{\infty},~ U_{\infty}}$ are parametrized by $t\in \kk$, while $\mathrm{O_{\infty}}$ admits two families of $\kk$-valued measures, $\mu^{even}_t, \mu^{odd}_t$, which coincide for $t=0$. These measures are known to be regular if $t\notin \{0, \pm q^n~;~ n\in \Z_{\geq 0}\}$ (see \cite{HS-classical}). The full list of values for which these measures are non-regular was established in \cite{HS-classical} (see \cref{ssec:singular_param}).

\begin{introtheorem}[See \cref{subsec:wreath-products-iedals}, \cref{sec:singular-non-zero}]\label{introthm:boring} 
Let $\mathrm{char}~\kk = 0 $. Assume one of the following settings:
\begin{itemize}
    \item Let $G = \Gamma \wr S_{\infty}$ for some finite group $\Gamma$ (e.g. $G=S_{\infty}$), or the classical group $\mathrm{GL^{par}_{\infty}}$ over a finite field. Let $\mu$ be a non-regular measure on $G$. 
    
    \item Let $G$ be one of the classical groups $\mathrm{GL^{Levi}_{\infty}, ~Sp_{\infty},~ U_{\infty},~ O_{\infty}}$ over a finite field. {Let} $\mu$ be a non-regular measure on $G$ with measure parameter $t\neq 0$.
\end{itemize}
 Then $\Perm(G,\mu)$ has only one non-trivial tensor ideal: the ideal of negligible morphisms $\mathcal{N}$.

\end{introtheorem}

We also point out that \cite{CS-delannoy} showed that the second Delannoy category also has just one proper non-trivial tensor ideal. At that point one might wonder whether the classification is maybe a bit boring. Fortunately, the poset $\mathcal{P}_{\mu}^G$ can be much more complicated. A quick way to get more complex examples is by looking at permutation categories of finite groups in positive characteristic.

\begin{introtheorem}[Theorem \ref{thm:E-p-subgroups}]
For a finite group $G$ over a field of characteristic $p$, one has
\[
\mathcal P_\mu^G\cong \Sub_p(G)/G,
\]
where $\Sub_p(G)/G$ is the poset of $G$-conjugacy classes of $p$-subgroups, ordered by subconjugacy.  
\end{introtheorem}

Another way to get {more} complex examples is by taking products of oligomorphic groups.  

\begin{introtheorem}[\cref{thm:products}]
Let $\mathrm{char}~\kk=0$. Let $\mu$ be a measure on $G_1\times G_2$, and let $\mu_1,\mu_2$ be the restricted measures on $G_1$ and $G_2$. Then there is a
monoidal equivalence
\[ \Perm(G_1,\mu_1)^{\Kar}\boxtimes \Perm(G_2,\mu_2)^{\Kar}
 \simeq \Perm(G_1\times G_2,\mu)^{\Kar}. \]
\end{introtheorem}

Given this it is easy to see that we have
\[
  \mathcal{P}_{\mu}^G\cong \mathcal{P}_{\mu_1}^{G_1}\times \mathcal{P}_{\mu_2}^{G_2}
\]
as posets. This can of course be iterated to $d$ factors. Already in an easy situation where the posets $\mathcal{P}_{\mu_i}^{G_i}$ consist of just two elements, say for $G_i = S_\infty$ and non-regular measure $\mu_i$, we obtain that the poset $\mathcal{P}_{\mu}^G$ of $G = S_\infty^{\times d}$ is the Boolean lattice, and the number of its downward-closed subsets is given by the $d$-th Dedekind number (these are $3, 6, 20, 168, 7581, \ldots$).

\subsection{Classical groups with singular measures}

The last part of this article is devoted to the study of the categories $\Perm(G,\mu)$ where $G$ is a classical group over a finite field $\F$ and the ground field $\kk$ for $\Perm(G,\mu)$ has characteristic $0$. When $t \neq 0$ and the measure is not regular, Theorem \ref{introthm:boring} together with the regular measure case describe all the ideals in $\Perm(G, \mu)$, so we now focus on the case when $t=0$. 

 For $t=0$, the measure $\mu_{t=0}$ on $\mathrm{GL^{par}_{\infty}}$ is regular (so $\Perm(\mathrm{GL^{par}_{\infty}}, \mu_{0})$ has only trivial ideals) and the Karoubi completion $$\mathscr{G}_{\F}=\Perm(\mathrm{GL^{par}_{\infty}},\mu_0)^{\Kar}$$ of $\Perm(\mathrm{GL^{par}_{\infty}},\mu_0)$ is semisimple. So we restrict our attention to the groups $G=\mathrm{GL^{Levi}_{\infty}, ~Sp_{\infty},~ U_{\infty},~ O_{\infty}}$.

At $t =0$, the measures $\mu_{t=0}$ on these four classical groups are no longer quasi-regular, and the categories $\Perm(G,\mu_0)$ are much more complicated.

\begin{introtheorem}[\cref{thm:zero-chain} and \cref{thm:orthogonal-zero}] 
\mbox{}

Let $\mathrm{char}~\kk = 0$. For $G = \mathrm{GL^{Levi}_{\infty}, ~Sp_{\infty},~ U_{\infty}}$ the non-trivial tensor ideals form a linear chain \[ \mathcal{N} = \I_1 \supset \I_2 \supset \ldots\] For $G = \mathrm{O_{\infty}}$ the non-trivial tensor ideals form an infinite chain of diamonds.
\end{introtheorem}

Apart from the mere classification, we can realize these tensor ideals as kernels of natural monoidal functors into abelian tensor categories. First of all, we compute explicitly the semisimplification of the idempotent completion $\Perm(G, \mu_0)^{\Kar}$:
\begin{introtheorem}[\cref{thm:rankzero}]
    The category $\Perm(G,\mu_0)^{\Kar}/\mathcal{N}$ is semisimple and admits an explicit description in terms of the parabolic general linear group: 
\[
\Perm(G,\mu_0)^{\Kar}/\mathcal{N}\cong
 \begin{cases}
 \mathscr{G}_{\mathbb{F}}\boxt\mathscr{G}_{\F},&\text{ if } G= \mathrm{GL^{Levi}_{\infty}},\\
 \mathscr{G}_{\F},&\text{ if }G=\mathrm{Sp_{\infty},~ O_{\infty}}\\
  \mathscr{G}_{\F'},&\text{ if }G=\mathrm{U_{\infty}} \text{ (here } \F' \text{ is a quadratic extension of }\F).\\
 \end{cases}
\]

\end{introtheorem}
Next, let $U\leq G$ be an open subgroup which is isomorphic to $G$ as an oligomorphic group, and consider the composition of functors
\[
 \Phi_{U}:\Perm(G,\mu_0) \xrightarrow{\Res_{U}}
 \Perm(U,\mu_0|_U)\cong \Perm(G,\mu_0)\xrightarrow{~\text{ semisimplification }}\mathscr \Perm(G,\mu_0)^{\Kar}/\mathcal{N}.
\]

\begin{introtheorem}[\cref{thm:kernels}] Let $G = \mathrm{GL^{Levi}_{\infty}, ~Sp_{\infty},~ U_{\infty}}$. For any fixed $r\geq 1$ there exists an open subgroup $U_r\leq G$ isomorphic to $G$, for which we have $\I_r = \ker(\Phi_{U_r})$.
\end{introtheorem}

This shows in particular that all $\I_r$ are tensor ideals of abelian type in the sense of \cite{CEO-abelian}.

There is a similar picture for $\mathrm{O_{\infty}}$ but here not all tensor ideals in this case are kernels of {functors} into abelian tensor categories (see \cref{thm:kernels}, \cref{cor:abelian-nonzero}). We also give a deformation theoretic interpretation of the tensor ideals in Section \ref{sec:higher-negligible}.

\subsection{Aftermath} We would like to point out at least three interesting directions for further research.
Firstly, in case $\mu$ is quasi-regular and satisfies the (Nil) property (the trace of any nilpotent endomorphism vanishes), $\Perm(G,\mu)$ admits an abelian envelope in the sense of \cite{CEAH,Coulembier-abelian}. In this setting the minimal thick tensor ideal of $\Perm(G,\mu)^{\Kar}$ constructed in Proposition \ref{prop:quasiregular} coincides with the ideal of splitting objects, and therefore with the projective objects in the abelian envelope $\uRep(G,\mu)$. In the general situation there is little hope to understand the tensor ideals of $\uRep(G,\mu)$ (consider for example the case of finite groups in positive characteristic - the Balmer spectrum provides an abstract description of the thick tensor ideals in the stable category, but little else is known).

Secondly, there are many more examples out there. While it would be possible to include more of these, we felt that the article was long enough and gave a reasonably rich account of the beautiful theory of tensor ideals in $\Perm(G,\mu)$. 

In an upcoming paper, joint with Harman and Snowden, we use the classification of ideals in the $t=0$ case to study the local abelian envelopes of the permutation category for the Levi version of the general linear group.

\subsection{Acknowledgments} We thank Nate Harman and Andrew Snowden for helpful discussions. The epigraph was added to the draft in December 2025 as an ironic bow to the enormous success of the Harman-Snowden construction and our subsequent move into oligomorphic country. It got a far more sinister meaning in August 2026 when we tried GPT Sol for the first time. While the theoretical part was largely finished by that time, several sample calculations were initially done with the help of GPT Sol. Claude Opus and GPT Astra were used to proofread the paper.

%%%%%%%%%%%%%%%%%%%%%%%%%%%%%%%%%

\part{Preliminaries}

\section{SM categories and SM functors} 
\subsection{SM categories}
Let $\kk$ be any field. We adopt the same notion of a $\kk$-linear symmetric monoidal (SM) category as \cite{Deligne-Milne}. In particular we have a binatural family of braiding morphisms $\gamma_{XY}:X\otimes Y\stackrel{\sim}{\to} Y\otimes X$ which satisfy the constraints of \cite{Deligne-Milne}. For an object $X\in \mathcal{T}$, let $X^*$ {denote} its dual (if it exists). If every object has a dual, we call $\mathcal{T}$ \emph{rigid}.

\medskip
A $\kk$-linear functor $F:\mathcal{T}\to \mathcal{T}'$ between two SM categories $\mathcal{T}$ and $\mathcal{T}'$ is symmetric monoidal (\emph{SM} for short) if it is equipped with natural isomorphisms $c_{XY}^F:F(X)\otimes F(Y)\stackrel{\sim}{\to}F(X\otimes Y)$ and $\triv\stackrel{\sim}{\to}F(\triv)$ satisfying the usual  compatibility conditions \cite[Definition 1.8]{Deligne-Milne}. If $F$ is a monoidal functor and $X^*$ the dual of $X$, then $F(X^\ast)$ is a dual of $F(X)$.

\subsection{Presheaves}

Let $\mathcal{T}$ be a small $\kk$-linear category. A \emph{left $\mathcal{T}$-module} is a $\kk$-linear functor $M:\mathcal{T}\to \mathrm{Mod}_\kk$, while a \emph{right $\mathcal{T}$-module}, or $\kk$-linear presheaf on $\mathcal{T}$, is a $\kk$-linear functor $\mathcal{T}^{op}\to \mathrm{Mod}_\kk$. We denote the corresponding functor categories by $\mathcal{T}\text{-}\mathrm{Mod}$ and $\mathrm{Mod}\text{-}\mathcal{T}$, respectively. Their morphisms are natural transformations, and kernels and cokernels are computed pointwise; in particular these categories are abelian.

For $X\in\mathcal{T}$, the covariant representable functor
\[
P_X=\Hom_{\mathcal{T}}(X,-)\in \mathcal{T}\text{-}\mathrm{Mod}
\]
is projective: the $\kk$-linear Yoneda lemma gives a natural isomorphism $\operatorname{Nat}(P_X,M)\cong M(X)$, and evaluation at $X$ is exact \cite[Chapters II, \S4 and III, \S2]{MacLane}, \cite[\S1.3.3--1.3.4]{Coulembier-ideals}.

%%%%%%%%%%%%%%%%

\subsection{Karoubi completion}

Let $\mathcal{T}$ be a $\kk$-linear monoidal category, and let $\mathcal{T}^{\oplus}$ be its additive envelope. We define
\[
\mathcal{T}^{\Kar}=\mathcal{T}^{\oplus, \Kar},
\]
the additive idempotent completion of $\mathcal{T}$ (called the \emph{Karoubi completion} or \emph{Karoubi envelope}).  Its objects are pairs $(X,e)$, where $X$ is a finite direct sum of objects of $\mathcal{T}$ and $e\in\End(X)$ satisfies $e^2=e$, and \[\Hom_{\mathcal{T}^{\Kar}}((X,e),(Y,f))=f\Hom_{\mathcal{T}^{\oplus}}(X,Y)e. \]
The identity of $(X,e)$ is $e$. The monoidal structure extends by
\[ (X,e)\otimes(Y,f)=(X\otimes Y,e\otimes f), \] so $\mathcal{T}^{\Kar}$ is again monoidal; braidings and duals, when present, extend as well \cite[\S2.6]{Coulembier-ideals} \cite[Lemme 1.3.10]{Andre-Kahn}.

%%%%%%%%%%%%%

\subsection{Tensor ideals}

Let $\mathcal{T}$ be a monoidal category with $\End_{\mathcal{T}}(\triv)=\kk$.

\begin{definition} A \emph{tensor ideal} $\I$ in $\mathcal{T}$ consists of a {$\kk$-subspace} $\I(X,Y) \subseteq \Hom(X,Y)$ for all $X,Y \in {\mathcal{T}}$ such that 
\begin{itemize} 
\item for all $X,Y,Z,W \in {\mathcal{T}}$ and $g \in \Hom(X,Y)$ and $h \in \Hom(Z,W)$ \[ f \in \I(Y,Z) \text{ implies } f \circ g \in \I(X,Z) \text{ and } h \circ f \in \I(Y,W); \] 
\item $f \in \I(X,Y)$ implies $\id_Z \otimes f \in \I(Z \otimes X,Z \otimes Y)$ and likewise from the right.
\end{itemize}

\end{definition}
The \emph{trivial} tensor ideals are the zero ideal ($\I(X,Y)=0$ for any $X, Y$) and the full category ($\I(X,Y)=\Hom(X, Y)$ for any $X, Y$).
\begin{definition} Let $\mathcal{T}$ be an additive braided monoidal category. A \emph{thick tensor ideal} of ${\mathcal{T}}$ is a replete full additive subcategory $\mathcal{I}\subseteq\mathcal{T}$ such that:
\begin{enumerate}
\item[(i)] $X\otimes Y\in\mathcal{I}$ whenever $X\in\mathcal{T}$ and $Y\in\mathcal{I}$;
\item[(ii)] $\mathcal{I}$ is closed under direct summands: if $X\in{\mathcal{T}}$, $Y\in\mathcal{I}$ and there exist $\alpha:X\to Y$, $\beta:Y\to X$ such that $\beta\circ\alpha=\id_X$, then $X\in\mathcal{I}$.
\end{enumerate}
\end{definition}

\begin{definition} For a given monoidal category $\mathcal{T}$ we denote the set of tensor ideals by $\TID(\mathcal{T})$ and the set of thick tensor ideals by $\ID(\mathcal{T})$.
\end{definition}

To any tensor ideal $\I$ we can associate the thick ideal $\I$ given by \[ \I = \{ X \in \mathcal{T} \ | \ \id_X \in \I(X,X) \}.\] This defines a surjective map, the decategorification map, \[ Ob: \TID(\mathcal{T}) \to \ID(\mathcal{T}) \] from the set of tensor ideals to the set of thick tensor ideals. It is not always injective \cite{Coulembier-ideals}.

%\medskip
%Let $K_0^+(\mathcal{T})$ be the split Grothendieck semiring of $\mathcal{T}$. Its elements are isomorphism classes of objects, with $[X] + [Y] = [X \oplus Y]$ and $[X][Y] = [X \otimes Y]$. The following is well-known.

%\begin{lemma}\label{lem:groth}
%Assume that $\mathcal{T}$ is additive and Karoubian. There is an inclusion-preserving bijection between thick tensor ideals $\mathcal{I}\subseteq\mathcal{T}$ and subsets $\mathfrak{j}\subseteq K_0^+(\mathcal{T})$ satisfying
%\begin{enumerate}
%\item $0\in\mathfrak{j}$ and $a,b\in\mathfrak{j}$ imply $a+b\in\mathfrak{j}$;
%\item $a\in\mathfrak{j}$ and $c\in K_0^+(\mathcal{T})$ imply $ca\in\mathfrak{j}$;
%\item $a+b\in\mathfrak{j}$ implies $a,b\in\mathfrak{j}$.
%\end{enumerate}
%Such a subset is called an order-ideal of $K_0^+(\mathcal{T})$.
%\end{lemma}

Assume $\mathcal{T}$ is rigid SM with $\End_{\mathcal{T}}(\triv)=\kk$.
\begin{definition}\label{def:negligible}\cite{Barrett-Westbury} A morphism $f \in \Hom(X, Y)$ is \emph{negligible} if $\mathrm{Tr}(g \circ f) = 0$ for all $g \in \Hom(Y, X)$. The collection of all such morphisms is denoted $\mathcal{N}$, the \emph{tensor ideal of negligible morphisms}.
\end{definition}

By \cite{Andre-Kahn} the tensor ideal {$\mathcal{N}$} is the largest proper tensor ideal of $\mathcal{T}$.

%%%%%%%%%%%%%%%%
\subsection{Tensor ideals and projective presheaves}

Assume $\mathcal{T}$ to be rigid SM, Krull-Schmidt, and $\End(\one) = \kk$ a field. 

Denote by $\TID(\mathcal{T})$ the set of tensor ideals. For $Z \in \mathcal{T}$ let $P_{Z} = \Hom_{\mathcal{T}}(Z,-) \in \mathcal{T}\text{-}\mathrm{Mod}$. For a tensor ideal $\mathcal{J}$ let $\mathcal{J}(X,Y)$ denote the usual submodule of $\Hom(X,Y)$.

We use the following theorem of Coulembier
\cite[Theorem~3.1.1 and Proposition~3.2.1]{Coulembier-ideals}.

\begin{theorem}(\cite{Coulembier-ideals})\label{thm:Coulembier}
Let $\mathcal T$ be a right-rigid braided monoidal preadditive category. Restriction
to morphisms out of the tensor unit induces an isomorphism of posets
\[
 \mathcal{TID}(\mathcal T)
 \xrightarrow{\sim}
 \operatorname{Sub}_{\mathcal T}\Hom_{\mathcal T}(\one,-),
 \qquad
 \I\longmapsto \I(\one,-).
\]
If $M\subset\Hom(\one,-)$ is a subfunctor, the inverse ideal is, under the
rigidity isomorphism
\[
 \iota_{A,B}:\Hom(\one,A^{\vee}\otimes B)\xrightarrow{\sim}\Hom(A,B),
\]
given by
\[
 \I_M(A,B)=\iota_{A,B}\bigl(M(A^{\vee}\otimes B)\bigr).
\]
\end{theorem}
\subsection{Tensor ideals in the Karoubi completion}

The following is well-known; see \cite[\S2.6.2]{Coulembier-ideals} and \cite[Lemme 1.3.10]{Andre-Kahn}.

\begin{lemma}\label{lem:ideals-kar}
There is a canonical, inclusion-preserving bijection between the tensor ideals of $\mathcal{T}$ and the tensor ideals of $\mathcal{T}^{\Kar}$.
\end{lemma}

\begin{proof}
Let $\mathcal{I}_0$ be a tensor ideal of $\mathcal{T}$ and denote by $\mathcal{I}_0^{\oplus}(X,Y)$ its extension to the additive envelope $\mathcal{T}^{\oplus}$. For objects $(X,e),(Y,f)$ of $\mathcal{T}^{\Kar}$ define
\[
\Phi(\mathcal{I}_0)((X,e),(Y,f))
   =f\,\mathcal{I}_0^{\oplus}(X,Y)\,e.
\]

Conversely, if $\mathcal{I}$ is a tensor ideal of $\mathcal{T}^{\Kar}$, let
\[
\Psi(\mathcal{I})(X,Y)
 =\mathcal{I}((X,\id_X),(Y,\id_Y))
\]
for $X,Y\in\mathcal{T}$. This is a tensor ideal of $\mathcal{T}$, and clearly $\Psi\Phi(\mathcal{I}_0)=\mathcal{I}_0$ and $\Phi\Psi(\mathcal{I})=\mathcal{I}$. 
\end{proof}

\subsection{Semisimplification} 

The \emph{semisimplification} of a rigid SM category $\mathcal{T}$ is the quotient $\bar{\mathcal{T}} = \mathcal{T}/\mathcal{N}$ \cite{Barrett-Westbury}. Its objects {consist} of the objects of $\mathcal{T}$ and its morphism spaces are given by \[ \Hom_{\mathcal{T}/\mathcal{N}}(X,Y) = \Hom(X,Y)/\mathcal{N}(X,Y).\] If $\mathcal{T}$ is abelian, $\bar{\mathcal{T}}$ is always semisimple. If the categorical trace of every nilpotent morphism in $\mathcal{T}^{\Kar}$ vanishes, the quotient $\mathcal{T}^{\Kar}/\mathcal{N}$ is semisimple.

\section{Oligomorphic groups}
\subsection{Definitions}\label{ssec:olig_defs}
We briefly assemble some background about oligomorphic groups and the associated monoidal categories. For more information on oligomorphic groups see \cite{Cameron-book} and \cite[Definitions 2.1 and 2.3 and Proposition 2.4]{HS-oligomorphic} for the monoidal categories.
\begin{definition}

A topological group $G$ is called \emph{pro-oligomorphic} (we will call them \emph{oligomorphic} for short) if the following $3$ conditions {hold}: 
\begin{enumerate}
    \item open subgroups of $G$ form a neighborhood base of the identity (that is, $G$ is non-Archimedean),
    \item the intersection of all open subgroups of $G$ is $\{1_G\}$ (that is, $G$ is Hausdorff),
    \item the set of double cosets $U\backslash G/V$ is finite for any two open subgroups $U,V\leq G$ (that is, $G$ is Roelcke precompact).
\end{enumerate}
\end{definition}
Whenever we write $U\leq G$, this will mean that $U$ is an open subgroup of $G$. We will denote by $\Sub(G)$ the set of all open subgroups of $G$.
\begin{example}
    Any finite group with discrete topology is an oligomorphic group.
\end{example}

Let us explain the main source of non-trivial examples of oligomorphic groups.
Consider a faithful action of a group $G$ on a set $\Omega$. This action is called \emph{oligomorphic}, and the pair $(G, \Omega)$ is called \emph{an oligomorphic pair}, if $G$ has only finitely many orbits on $\Omega^n$ for every $n\geq0$. The pointwise stabilizers of finite subsets of $\Omega$ define the permutation topology on $G$, forming a base of open subgroups in this topology. With this topology $G$ is pro-oligomorphic. 

A $G$-set is called \emph{smooth} if all point stabilizers are open and \emph{finitary} if it has finitely many $G$-orbits; a transitive smooth $G$-set is of the form $G/U$ with $U\in \Sub(G)$.

In this text, all the $G$-sets we consider will be smooth and finitary; the category of such $G$-sets will be {denoted} by $\mathbf{S}(G)$, while the full subcategory of transitive smooth $G$-sets will be {denoted} $\mathbf{S}(G)^{tran}$. We denote by $\bullet_G$ the trivial one-point $G$-set, which is the terminal object in $\mathbf{S}(G)$. An important feature of oligomorphic groups is that $\mathbf{S}(G)$ is closed under taking direct products of $G$-sets. For $X\in \mathbf{S}(G)$, we denote by $G\backslash X$ the set of $G$-orbits on $X$.

\begin{definition}[{\cite[Definition~3.1]{HS-oligomorphic}}]\label{def:measure}
Let $G$ be pro-oligomorphic and let $\kk$ be a commutative ring. A $\kk$-valued
\emph{measure} on $G$ assigns a scalar $\mu_U(X)\in \kk$ to every $X\in \Set{U}$, for every $U\in \Sub(G)$, subject to the following axioms:
\begin{enumerate}
\item isomorphic $U$-sets have the same measure;
\item $\mu_U({\bullet_U})=1$;
\item $\mu_U(X\sqcup Y)=\mu_U(X)+\mu_U(Y)$;
\item the measure is invariant under conjugation: if $g\in G$, then
      $\mu_{gUg^{-1}}({}^{g}X)=\mu_U(X)$;
\item if $f:X\to Y$ is a map in $\Set{U}^{tran}$ and $F$ is
      a fiber of $f$, regarded as a set for the stabilizer of a point of $Y$, then
      $\mu_U(X)=\mu(F)\mu_U(Y)$.
\end{enumerate}
When the acting subgroup is clear, we write simply $\mu(X)$.
\end{definition}

Additionally, we require measures to be compatible with restriction of the acting
subgroup:
\begin{equation}\label{eq:measure-restriction-compatibility}
  \mu_V(\Res^U_V X)=\mu_U(X)
  \qquad (V\leq U\leq G),
\end{equation}
for every finitary smooth $U$-set $X$ (this is implicit in the
$\widehat G$-set formulation of
\cite[\S2.5 and Definition~3.1]{HS-oligomorphic}).

\begin{notation}
Let $f:Z\to X$ be a $G$-map between $G$-sets. Assume that $X$ is transitive. Then the fibers of $f$ all have the same measure, which we will denote by $\mu(f)$ for short.
\end{notation}

A measure $\mu$ is \emph{regular} if
$\mu_G(G/U)\in\kk^\times$ for every $U\in\Sub(G)$, and
\emph{quasi-regular} if its restriction $\mu|_U$ is regular for some
$U\in\Sub(G)$. This is equivalent to
$\mu_U(U/V)\in\kk^\times$ for every open $V\leq U$
\cite[Definition~3.18]{HS-oligomorphic}.
Both regularity and quasi-regularity are preserved by restriction to open subgroups \cite[Proposition~3.20]{HS-oligomorphic}.

The \emph{universal measure ring} $\Theta(G)$ is the commutative unital
$\mathbb Z$-algebra generated by symbols $[X]_U$, indexed by isomorphism
classes of finitary smooth $U$-sets for all $U\in\Sub(G)$, subject to the
relations obtained from Definition~\ref{def:measure} and
\eqref{eq:measure-restriction-compatibility} by replacing
$\mu_U(X)$ with $[X]_U$
\cite[Definition~3.6]{HS-oligomorphic}.
Thus $\mu_{\mathrm{univ},U}(X)=[X]_U$ is a measure, and
\[
  \Hom_{\mathrm{CRing}}(\Theta(G),\kk)
  \cong \{\text{$\kk$-valued measures on $G$}\},
\]
via $\mu_U(X)=\varphi([X]_U)$ (we usually suppress $U$ in the notation). The multiplication satisfies $[X][Y]=[X\times Y]$ \cite[Remark~3.7]{HS-oligomorphic}.

%%%%%%%%%%%%%%

\subsection{Permutation representations for oligomorphic groups}
Let $G$ be pro-oligomorphic and let $\mu$ be a $\kk$-valued measure in the sense of \cite[Definition 3.1]{HS-oligomorphic}. Thus $\mu$ assigns a scalar to every finitary smooth $G$-set, is invariant under isomorphism and conjugation, is normalized by $\mu(\bullet_G)=1$, is additive under disjoint unions, and is multiplicative in transitive fibrations: if $X\to Y$ is a map of transitive smooth $U$-sets with fiber $F$, then $\mu(X)=\mu(F)\mu(Y)$.

For a smooth $G$-set $X$, a $\kk$-valued function is a \emph{Schwartz function} if it is invariant under some open subgroup and if it has finitary support; the space of such functions is denoted $\cC(X)$ \cite[\S4.2]{HS-oligomorphic}. The category $\Perm(G,\mu)$ has objects $\mathcal{C}(X)$, for $X \in \mathbf{S}(G)$ and morphisms given by
\[
\Hom_{\Perm(G,\mu)}(\mathcal{C}(X),\mathcal{C}(Y))=\cC(Y\times X)^G.
\]
If $A\in\cC(Y\times X)^G$ and $B\in\cC(Z\times Y)^G$, composition is convolution
\[
(B\circ A)(z,x)=\int_Y B(z,y)A(y,x)\,d\mu(y).
\]
Because $Y \times X$ is finitary, it decomposes into a finite number of $G$-orbits $O_1, \dots, O_n$. Then $\Hom(\mathcal{C}(X), \mathcal{C}(Y))$ has a basis given by the indicator functions $\{\one_{O_1}, \dots, \one_{O_n}\}$.

Disjoint union gives direct sum, Cartesian product gives tensor product, and the tensor unit is $\mathcal{C}(\bullet_G)$; this produces an additive rigid SM category \cite[Definition 8.1 and Propositions 8.2--8.3]{HS-oligomorphic}.

We shall repeatedly use the following special case of multiplicativity in
fibrations. If $K\leq H\leq G$ are open subgroups, then the canonical map
$G/K\to G/H$ has fibre $H/K$, and
\begin{equation}\label{eq:tower-formula}
  \mu_G(G/K)=\mu_H(H/K)\,\mu_G(G/H).
\end{equation}
We refer to \eqref{eq:tower-formula} as the \emph{Tower formula}. If
$[H:K]<\infty$, then $\mu_H(H/K)=[H:K]\cdot 1_\kk$.

\begin{lemma}\label{lem:fibre-splitting}
Let $p:X\to Y$ be a map of finitary $G$-sets, with $Y$ transitive. Let
$c = \mu(p)$ be the measure of a fibre. If $c\in\Bbbk^\times$, then
$\mathcal{C}(Y)$ is a retract of $\mathcal{C}(X)$ and therefore a direct summand of $\mathcal{C}(X)$ in the Karoubi envelope $\Perm(G, \mu)^{\Kar}$.
\end{lemma}

\begin{proof}
This follows directly from \cite[Proposition~7.21(c)]{HS-oligomorphic}.
\end{proof}

\subsection{Main examples of oligomorphic groups}\label{subsec:main-examples}

The following families are the fundamental examples of oligomorphic groups.
\begin{enumerate}[itemindent=0pt,leftmargin=*]
\item Every finite group $G$ is oligomorphic. The unique measure for a finite group $G$ is the counting measure $\mu(X)=|X|\cdot 1_{\kk}$ for $X\in \Set{G}$. The category $\Perm(G,\mu)$ in this case agrees with the category of permutation representations of $G$ over $\kk$. Their tensor ideals are studied in Section \ref{sec:finite-groups}.

\item Assume $\mathrm{char} ~ \kk = 0$. For $G=S_\infty= \Aut(\mathbb{N})$ acting on the countably infinite set $\Omega=\mathbb{N}$, the universal parameter  is $t=\mu(\Omega) \in \kk$, and the additive Karoubi completion of $\Perm(S_\infty,\mu_t)$ is Deligne's category $\uRep(S_t)$ for $t \in \kk$ \cite[\S14]{HS-oligomorphic}, \cite{Del07}.
      
\item The order-preserving automorphism group $\Aut(\bR, <)$ of the real line, acting on $\Omega=\bR$. On this group there are four measures \cite[Theorem 16.7]{HS-oligomorphic}, resulting in four permutation categories called the Delannoy categories. The first is studied in \cite{HSS-delannoy}, and all four are treated in \cite{CHS-delannoy}.

\item Finite classical groups give rise to oligomorphic groups (see \cref{ssec:setup_classcial_grps}). Over a field of characteristic zero, these are studied in depth in \cite{HS-classical} where their measures are classified \cite[\S\S6,8,9,10]{HS-classical}. We study their tensor ideals in Part \ref{sec:classical}. 

For example, the pair $(GL_{\infty}(\mathbb{F}_q) = \bigcup_n GL_n(\mathbb{F}_q), \F_q^{\infty})$ is oligomorphic \cite{HS-oligomorphic}, corresponding to the  general linear group (parabolic version) in \cref{ssec:setup_classcial_grps}. The Karoubi completions of the categories $\Perm(G,\mu)$ in this case coincide with the categories studied in \cite{Knop} and \cite{EAHGL}. 
%Other examples include a Levi version of the general linear group and the symplectic, orthogonal, and unitary analogues. 

\end{enumerate}

\subsection{Restriction functors}\label{sec:restriction}

The following proposition is a special case of the restriction construction in \cite[Proposition~3.4]{HS-classical}.

Let $U\in \Sub(G)$. There is a faithful SM functor
\[
 \Res_U^G:\Perm(G,\mu)\longrightarrow\Perm(U,\mu|_U)
\]
defined as follows:
\begin{enumerate}[label=\textup{(\roman*)}]
\item a finitary smooth $G$-set $X$ is regarded as the same set with the restricted $U$-action;
\item a $G$-invariant function $a:Y\times X\to\kk$ is regarded as the same function, now merely $U$-invariant.
\end{enumerate}
The restriction functor extends to the additive Karoubi completions: $\Res_U^G:\Perm(G,\mu)^{\Kar}\longrightarrow\Perm(U,\mu|_U)^{\Kar}$.

\subsection{The relative setting}\label{sec:relative}
Let $G$ be a (pro-) oligomorphic group and let $\mathcal{E}$ be a stabilizer class: that is, a collection of open subgroups in $G$ satisfying the following conditions (see \cite[Section 2.6]{HS-oligomorphic}):
\begin{enumerate}
    \item $G\in \mathcal{E}$,
    \item $\mathcal{E}$ is stable under conjugation,
    \item $\mathcal{E}$ is stable under finite intersections,
    \item $\mathcal{E}$ forms a neighborhood base of $1_G$.
\end{enumerate}

\begin{example}\label{ex:Omega_stabilizer_class}
    Let $(G, \Omega)$ be an oligomorphic pair as in \cref{ssec:olig_defs} and consider the induced oligomorphic structure on $G$. The pointwise stabilizers of finite subsets in $\Omega^n$, $n\geq 0$ form the $\Omega$-stabilizer class $\mathcal{E}_{\Omega}$.
\end{example}

We may now consider {an} $\mathcal{E}$-relative version of the theory described in the previous subsections, by replacing open subgroups with subgroups lying in $\mathcal{E}$, while smooth $G$-sets will be replaced by $\mathcal{E}$-smooth $G$-sets (ones whose stabilizers lie in $\mathcal{E}$). One can then define \emph{$\mathcal E$-relative measures} on $G$; in particular, the $\mathcal E$-relative measure of a map between two $\mathcal{E}$-smooth $G$-sets with a transitive base. Given such a measure $\mu$, one can define the category $\Perm(G, \mathcal{E}, \mu)$ in precisely the same manner as before. Clearly, any {``}absolute'' measure on $G$ will be a relative measure as well, for any stabilizer class.

The same construction as in \cref{ssec:olig_defs}, with acting subgroups and point stabilizers in $\mathcal E$, defines the relative universal measure ring
$\Theta(G,\mathcal E)$. This ring represents $\mathcal E$-relative measures (\cite[Remark~3.8]{HS-oligomorphic}).

Let us mention one useful link between the relative and the absolute setting.

\begin{definition}\label{def:saturating}
Let $G$ be a (pro-) oligomorphic group with a $\kk$-valued measure $\mu$ and let $\mathcal{E}$ be a stabilizer class of $G$. We say that a stabilizer class $\mathcal{E}$ is \emph{saturating} if for any $U\in \Sub(G)$, there exists $V\in \mathcal{E}${, $V\subseteq U$,} such that $\mu(U/V)\neq 0$.
\end{definition}
This is a variation of the notion of a \emph{large stabilizer class} (see {\cite[Section 3.4(d)]{HS-oligomorphic}}): there, one replaces the condition $\mu(U/V)\neq 0$ by the condition $|U/V|<\infty$.
\begin{example}
    Assume $\mathrm{char}~\kk = 0$. Then any large stabilizer class is saturating. 
\end{example}

\begin{proposition}\label{prop:saturating_stab_class}
    Let $G$ be a (pro-) oligomorphic group with a $\kk$-valued measure $\mu$. If $\mathcal{E}$ is a saturating stabilizer class of open subgroups in $G$, then $\Perm(G, \mathcal{E}, \mu)^{\Kar}\cong \Perm(G, \mu)^{\Kar}$.
\end{proposition}

\begin{proof}
We have an obvious embedding (fully faithful functor) $I:\Perm(G, \mathcal{E}, \mu) \to \Perm(G, \mu)$ inducing a fully faithful functor $\Perm(G, \mathcal{E}, \mu)^{\Kar}{\to} \Perm(G, \mu)^{\Kar}$. To show that the latter is an equivalence, we just need to show that for any $U\in \Sub(G)$, $\mathcal{C}(G/U)$ lies in the Karoubi envelope of the essential image of $I$.

    Let $U\leq G $ be an open subgroup. Since $\mathcal{E}$ is a saturating stabilizer class, we may choose $V\in \mathcal{E}$ such that $V\subseteq U$ and $\mu(U/V)\neq 0$. In other words, $\mu(G/V\to G/U)\neq 0$. By \cref{lem:fibre-splitting}, this means that $\mathcal{C}(G/U)$ is a direct {summand} of $\mathcal{C}(G/V)$, as required. 
\end{proof}

\begin{corollary} \label{cor:tensor-generator} If $(G,\Omega)$ is an oligomorphic pair and $\mathcal E_\Omega=\{G_S \ | \ S\subset\Omega\text{ is finite}\}$ is saturating for $\mu$, then the object $\mathcal{C}(\Omega)$ generates $\Perm(G, \mu)^{\Kar}$ as a Karoubian additive rigid SM category. 
\end{corollary}

\begin{proof} It suffices to prove that for any open subgroup $U \le G$, the object $\mathcal{C}(G/U)$ is a direct summand of some object generated by tensor powers of $\mathcal{C}(\Omega)$. Since the tensor product is induced by the Cartesian product of $G$-sets, the $n$-th tensor power is\[\mathcal{C}(\Omega)^{\otimes n} \cong \mathcal{C}(\Omega^n).\] Now $\Omega^n$ decomposes into a finite disjoint union of transitive $G$-sets $\Omega^n = \bigsqcup_{i=1}^k O_i$. For any orbit $O_i$, let $x = (x_1, x_2, \dots, x_n) \in O_i$ a point and define $S = \{x_1, \dots, x_n\} \subseteq \Omega$. Then $O_i \cong G/G_S$, yielding the decomposition \[ \mathcal{C}(\Omega)^{\otimes n} \cong \bigoplus_{i=1}^k \mathcal{C}(G/G_{S_i}).\] Therefore, for any finite subset $S \subseteq \Omega$, the object $\mathcal{C}(G/G_S)$ is a direct summand of $\mathcal{C}(\Omega)^{\otimes |S|}$. By the saturation condition, for any $U$ there exists a pointwise stabilizer $G_S \le U$ such that {$\mu(U/G_S)\neq 0$}. By Lemma \ref{lem:fibre-splitting}, $\mathcal{C}(G/U)$ is then a direct summand of $\mathcal{C}(G/G_S)$.
\end{proof}

\begin{example}\label{ex:saturating}
Assume $\operatorname{char}~\kk=0$. The following are examples of saturating
pointwise stabilizer classes.
\begin{enumerate}
\renewcommand{\labelenumi}{\textnormal{(\roman{enumi})}}
\item \emph{Finite groups.}
The stabilizer class {$\{\{1\},G\}$} is saturating.

\item \emph{The infinite symmetric group.}
For $S_\infty$ this holds for every measure by \cite[Proposition~14.1]{HS-oligomorphic}.

\item \emph{The ordered line.}
For $G=\operatorname{Aut}(\mathbb R,<)$ acting on $\mathbb R$, every open subgroup is itself
the pointwise stabilizer of a finite subset
\cite[Proposition~16.1]{HS-oligomorphic}, so saturation holds for every measure.

\item \emph{Regular measures.}
If $\mu$ is regular, every stabilizer class is saturating.

\item \emph{Finite products.}
Suppose $(G_i,\Omega_i)$ has a saturating pointwise stabilizer class for
$\mu_i$, for $i=1,2$. Then $G=G_1\times G_2$ acting on
$\Omega_1\sqcup\Omega_2$ has a saturating pointwise stabilizer class for
$\mu=\mu_1\boxtimes\mu_2$ (see Lemma \ref{lem:finite-index-product-subgroup}). This can be iterated for finitely many factors.
\end{enumerate}
\end{example}

%%%%%%%%%%%%%%%

\subsection{Traces in \texorpdfstring{$\Perm(G,\mu)$}{permutation categories}}

Let $G$ be an oligomorphic group and $\mu$ a $\kk$-valued measure on $G$.

Let $X\in \Set{G}$. The trace of $h \in \mathrm{End}(\mathcal{C}(X))$ in $\Perm(G,\mu)$ is evaluated by integrating along the diagonal
\[\mathrm{Tr}(h) = \int_X h(x, x) \, d\mu(x). \] Given an object $(\mathcal{C}(X),e)$ in $\mathcal{T}=\Perm(G,\mu)^{\Kar}$, an endomorphism $f \in \mathrm{End}_{\T}((\mathcal{C}(X), e))$ is a function $f: X \times X \to {\kk}$ such that $f = e \circ f \circ e$. The trace of $f$ on the subobject $(\mathcal{C}(X), e)$ is identical to its trace as a morphism on $\mathcal{C}(X)$
\[ \mathrm{Tr}_{(\mathcal{C}(X),e)}(f) = \mathrm{Tr}_{\mathcal{C}(X)}(f) = \int_X f(x, x) \, d\mu(x).\]

Let $X,Y\in \Set{G}$ and $f\in \cC(Y \times X), ~g\in \cC( X \times Y)$ two $G$-invariant Schwartz functions. 
The trace of the composition defines a symmetric bilinear pairing $B$ between the two $\Hom$-spaces
\begin{align*}
B(g, f)=\mathrm{Tr}(g \circ f) &= \int_X (g \circ f)(x, x) \, d\mu(x)= \int_X \int_Y g(x, y) f(y, x) \, d\mu(y) \, d\mu(x). 
\end{align*}

For any orbit $O \subseteq Y \times X$, its transpose $O^T = \{(x, y) \mid (y, x) \in O\}$ is a $G$-orbit in $X \times Y$. Then
\begin{align*}
B(\one_{O^T}, \one_O) &= \int_X \int_Y \one_{O^T}(x, y) \one_O(y, x) \, d\mu(y) \, d\mu(x) \\
&= \int_{X \times Y} \one_{O^T}(x, y) \, d(\mu \times \mu)(x, y) = \mu(O^T) = \mu(O)
\end{align*}
For distinct orbits $O_i \neq O_j$ in $Y\times X$, the pointwise product $\one_{O_i^T}(x, y) \one_{O_j}(y, x)$ is identically zero, so $B(\one_{O_i^T}, \one_{O_j}) = 0$.

\begin{corollary} \label{cor:trace-diagonal} Let $O,O'$ be two orbits in $Y \times X$.  Then \[ \operatorname{Tr}\bigl(\one_{O'^{T}}\circ\one_O\bigr)
 =\begin{cases}
      \mu(O),&\text{ if } O'=O,\\
      0,&\text{ if } O'\ne O.
   \end{cases}\]
\end{corollary}

%%%%%%%%

\subsection{Negligible morphisms in permutation categories}

Let $G$ be an oligomorphic group and $\mu$ a $\kk$-valued measure on $G$. We denote by $\mathcal{N}$ the tensor ideal of negligible morphisms (\cref{def:negligible}) in either $\Perm(G, \mu)$ or $\Perm(G, \mu)^{\Kar}$.

\begin{corollary}
Let $f = \sum_{i=1}^n c_i \one_{O_i} \in \Hom(\mathcal{C}(X), \mathcal{C}(Y))$. The morphism $f$ is negligible if and only if $c_i = 0$ for all orbits $O_i\subset Y\times X$ where $\mu(O_i) \neq 0$. Therefore \[ \mathcal{N}(\mathcal{C}(X), \mathcal{C}(Y)) = \Span_\kk \{ \one_O \mid O \text{ is a } G\text{-orbit in } Y \times X \text{ and } \mu(O) = 0 \text{ in } \kk \} \]
\end{corollary}

\begin{proof}
Because the basis vectors $\{\one_{O_i^T}\}$ span $\Hom(\mathcal{C}(Y), \mathcal{C}(X))$ and the pairing is orthogonal on orbits, $f$ is in the radical of the pairing if and only if $B(\one_{O_i^T}, f) = 0$ for all $i$. 
$$ B(\one_{O_i^T}, f) = c_i B(\one_{O_i^T}, \one_{O_i}) = c_i \mu(O_i) $$
For this to be zero for all $i$, we must have $c_i = 0$ whenever $\mu(O_i) \neq 0$ in $\kk$.
\end{proof}
\begin{corollary}
  The thick tensor ideal in $\Perm(G, \mu)$ associated with $\mathcal{N}$ is given by $\cC(X)$, $X\in \Set{G}$ where all the $G$-orbits of $X$ have measure $0$.  
\end{corollary}

\begin{proof}
   For $X\in \Set{G}^{tran}$, we have, by the above corollary: $$\cC(X) \in \mathcal{N}\; \Longleftrightarrow \; \id_{\cC{(X)}}\in \mathcal{N}(\mathcal{C}(X), \mathcal{C}(X)) \; \Longleftrightarrow \;\mu(\Delta X)=0 \; \Longleftrightarrow \; \mu(X)=0$$ where $\Delta X\subset X\times X$ denotes the diagonal orbit.
\end{proof}
%Correspondingly, the thick tensor ideal in $\Perm(G, \mu)^{\Kar}$ associated with $\mathcal{N}$ consists of all direct sums of indecomposable objects of categorical dimension $0$.
\begin{corollary}

The ideal $\mathcal{N}$ is trivial if and only if for every $X\in \Set{G}^{tran}$, we have $\mu(X) \neq 0$.

\end{corollary}

\begin{corollary}
    \label{lem:negl-basis}
The morphism space $\Hom_{\Perm(G, \mu)/ \mathcal{N}}(\mathcal{C}(X),\mathcal{C}(Y))$ has a $\kk$-basis in bijection with the transitive $G$-orbits of $Y \times X$ that have non-zero measure.
\end{corollary}

%%%%%%%%%%%%%%

%%%%%%%%%%%%%%%%%%

%%%%%%%%%%%%%%%%%

\part{Products and wreath products of oligomorphic groups}

In this section we prepare the groundwork for the examples of oligomorphic groups arising from direct products and wreath products, since details have been missing in the literature so far.

\section{Wreath products of a finite group with \texorpdfstring{$S_\infty$}{the infinite symmetric group}}\label{sec:wreath-products_intro}

In this section $\kk$ is a field, $\mathrm{char}~\kk=0$. 

Let $\Gamma$ be a finite group. Denote its full wreath product with the infinite symmetric group by
$G=\Gamma^{\Omega}\rtimes S_\infty$, where $\Omega=\{1,2,\ldots\}$. For the measures $\mu$ defined below, the categories $\Perm(G,\mu)^{\Kar}$  coincide with the finite-group wreath-product interpolation categories defined in \cite{Knop} and further studied in \cite{Nyobe-Savage} (under the name \emph{group partition category}) and in \cite{FHH}.

Let $d=|\Gamma|$. The natural faithful $G$-set is
$D=\Gamma\times\Omega$, with action given by
\[
  (\gamma,\sigma)\cdot(h,i)
  =(\gamma_{\sigma(i)}h,\sigma(i)).
\]
The subsets $B_i=\Gamma\times\{i\}$, for $i\in\Omega$, are called
\emph{blocks}. The blocks are the fibres of the obvious map $D\to\Omega$ (projection on the second factor) and form a
$G$-invariant partition, since $(\gamma,\sigma)B_i=B_{\sigma(i)}$.
Thus every translate of a block is equal to it or disjoint from it.
We can identify the set of blocks with $\Omega$ and $G$ acts on this set through
the quotient map $G\to S_\infty$.

For a finite subset $A\subset\Omega$, let $U_A$ be the pointwise stabilizer
of $\Gamma\times A$. Then
\[
  U_A\cong\Gamma^{\Omega\setminus A}\rtimes
       \operatorname{Sym}(\Omega\setminus A).
\]
These subgroups form a neighbourhood basis of the identity. Write
$U_n=U_{\{1,\ldots,n\}}$, with $U_0=G$.

\begin{lemma}\label{lem:wp-open-subgroups}
Every open subgroup $V\leq G$ has the form
\[
  V=H\times U_A,
  \qquad \text{ with } \qquad H\leq\Gamma^A\rtimes\operatorname{Sym}(A),
\]
for some finite $A\subset\Omega$. 

In particular, $[V:U_A]<\infty$,
and the stabilizer class
$\mathcal E=\{U_A \ | \ A\subset\Omega\text{ finite}\}$ is large.
\end{lemma}

\begin{proof}
We first claim that for finite subsets $A,B\subset\Omega$, the subgroup $\langle U_A,U_B\rangle$ generated by $U_A, U_B$ is $U_{A\cap B}$. One inclusion is obvious: $\langle U_A,U_B\rangle\subseteq U_{A\cap B}$.

Let us prove the second inclusion. The two base subgroups $\Gamma^{\Omega\setminus A}, \Gamma^{\Omega\setminus B}$ generate
$\Gamma^{\Omega\setminus(A\cap B)}$. For the permutation subgroups $\operatorname{Sym}(\Omega\setminus A), \operatorname{Sym}(\Omega\setminus B)$,
put $C=A\cap B$. A transposition of $\Omega\setminus C$ already
fixes at least one of the sets $A, B$ pointwise, unless it exchanges
$a\in A\setminus B$ with $b\in B\setminus A$. 
In the latter case,
choose $z\notin A\cup B$ and use the identity
\[
  (a\,b)=(a\,z)(b\,z)(a\,z).
\]
The outer transpositions fix $B$ pointwise and the middle one
fixes $A$ pointwise. Thus the two permutation subgroups generate
every finitely supported permutation fixing $C$.

Let us show that the same holds for any permutation (not necessarily finitely supported) which fixes $C$ pointwise. We let $\sigma$ be any permutation fixing $C$ pointwise and put
$A'=A\setminus C$. The partial bijection
$\sigma(A')\to A'$, $\sigma(a)\mapsto a$, extends to a permutation
of the finite set $A'\cup\sigma(A')$. Extending it by the identity elsewhere
gives a finitely supported permutation $\tau$ fixing $C$.
Then $\tau\sigma$ fixes $A$ pointwise. Both $\tau$ and
$\tau\sigma$ lie in the subgroup generated by the two pointwise
stabilizers, so $\sigma=\tau^{-1}(\tau\sigma)$ lies in the same subgroup as well. This proves $\langle U_A,U_B\rangle= U_{A\cap B}$.

We now choose $A$ of minimal cardinality among those finite subsets of $\Omega$ with $U_A\leq V$. For $v\in V$,
let $\bar v \in \operatorname{Sym}(\Omega)$ be its permutation of $\Omega$. Since
\[
  vU_Av^{-1}=U_{\bar v(A)},
\]
both $U_A$ and $U_{\bar v(A)}$ are contained in $V$. The identity
proved above therefore gives
\[
  U_{A\cap\bar v(A)}
  =\langle U_A,U_{\bar v(A)}\rangle
  \leq V.
\]
As $A\cap\bar v(A)\subseteq A$, the minimality of $|A|$ forces
$A\cap\bar v(A)=A$. Thus $A\subseteq\bar v(A)$, and equality
follows because these finite sets have the same cardinality.
Hence every element of $V$ preserves $A$ setwise.

We put $F_A=\Gamma^A\rtimes\operatorname{Sym}(A)$, viewed as the
subgroup acting trivially on $\Gamma\times(\Omega\setminus A)$.
An element of $G$ preserving $A$ setwise acts independently on
$\Gamma\times A$ and $\Gamma\times(\Omega\setminus A)$.
Its actions on these two parts belong to $F_A$ and $U_A$,
respectively. Since these subgroups act on disjoint parts,
they commute and have trivial intersection. The setwise
stabilizer of $A$ is therefore $F_A\times U_A$, and we obtain
\[
  U_A\leq V\leq F_A\times U_A.
\]

We identify $U_A$ with $\{1\}\times U_A$ and define
\[
  H=\{h\in F_A \ | \ (h,1)\in V\}.
\]
Then $H\leq F_A$. If $(h,u)\in V$, the containment
$\{1\}\times U_A\leq V$ gives
\[
  (h,1)=(h,u)(1,u^{-1})\in V,
\]
so $h\in H$ and hence $V\subseteq H\times U_A$.
Conversely, for $h\in H$ and $u\in U_A$, both $(h,1)$ and
$(1,u)$ belong to $V$, so $(h,u)\in V$.
We therefore obtain
\[
  V=H\times U_A.
\]
Since $F_A$ is finite, so is $H$, and
$[V:U_A]=|H|<\infty$.

\end{proof}

As wreath products are not treated in \cite{HS-oligomorphic}, we need to establish whether they admit any measures. Instead of trying to construct such measures directly, we use Snowden's results on the realizability of Knop's tensor envelopes in the oligomorphic world \cite{Knop,Snowden-Knop}. In \cite{Knop,Nyobe-Savage} interpolation categories on $G$ have been constructed which depend on a parameter $t \in \kk$. We use these existing constructions to define measures based on Knop's degree functions. 

Let $\mathcal A$ be the opposite of the category of finite free right $\Gamma$-sets. To be in keeping with Snowden's convention on regular
categories, we use $\mathcal A^\varnothing$, obtained by adjoining
a formal initial object as in \cite[Definition~2.2]{Knop}. Its principal objects are the original objects of
$\mathcal A$. An embedding $i:X\hookrightarrow Y$ of such sets becomes a
regular epimorphism $i^{\mathrm{op}}:Y\twoheadrightarrow X$ in
$\mathcal A$. A degree function is a function on regular epimorphisms in $\mathcal{A}^\varnothing$ satisfying certain conditions (see \cite[Definition 3.1]{Knop}). By \cite[\S8, Example~2, equation~(8.15)]{Knop},
the degree functions are parametrized by $\delta\in\kk$, the
degree of $(\varnothing\hookrightarrow\Gamma)^{\mathrm{op}}$.

\begin{definition}
   For a finite free right $\Gamma$-set $X$,
define
\[
  \Emb_\Gamma(X,D)
  =\{f:X\to D\mid f \text{ is injective and }
       \Gamma\text{-equivariant}\}.
\] 
The group $G=\operatorname{Aut}_\Gamma(D)$ acts on this set by
postcomposition.
\end{definition}

%Any element of this set specifies a finite substructure of $D$ together with an equivariant identification with $X$.

\begin{lemma}
    We have: $  \Emb_\Gamma(X,D) \in \Set{G}^{tran}$.
\end{lemma}
\begin{proof}
    Given two embeddings
$f_1,f_2:X\to D$, the map
\[
  f_2\circ f_1^{-1}:f_1(X)\longrightarrow f_2(X)
\]
is an equivariant bijection. Both complements
$D\setminus f_1(X)$ and $D\setminus f_2(X)$ consist of countably
infinitely many free $\Gamma$-orbits. By matching these orbits
and choosing an equivariant bijection on each matched pair,
we extend this map to an element $g\in G$. Then
$g\circ f_1=f_2$, which proves transitivity.
\end{proof}

Fix an embedding $i:X\to D=\Gamma\times \Omega$ with image $\Gamma\times A$ ($A\subset\Omega$ finite). The stabilizer of $i$ is $U_A$, so we get
\[
  \operatorname{Emb}_\Gamma(X,D)\cong G/U_A.
\]

Snowden's construction
\cite[Construction~6.7 and Proposition~6.9]{Snowden-Knop} describes the measures of the $G$-sets $\Emb(X, D)$ for every finite free right $\Gamma$-set $X$, depending on our choice of
$\delta \in \kk$. Therefore for every $\delta\in \kk$, this construction
 gives us a measure
relative to the stabilizer class
$\mathcal E=\{U_A\mid A\subset\Omega\text{ finite}\}$. By Lemma~\ref{lem:wp-open-subgroups} and
\cite[\S3.4(d)]{HS-oligomorphic}, it extends uniquely to an absolute $\kk$-valued
measure on $G$.

\begin{lemma}
    These are the only measures on $G$.
\end{lemma}
\begin{proof}

To see that these are all the measures, we will show that any measure on $G$ is obtained from some degree function on regular epimorphisms in $\mathcal{A}^\varnothing$. This is done using \cite[Proposition~6.11 b)]{Snowden-Knop}, which gives an explicit criterion (in the general setting of Knop categories) as to when a measure on $G$ is obtained from a degree function.

In our case, per \cite[Proposition~6.11 b)]{Snowden-Knop}, we need to check that the following condition holds for any measure $\mu$ on $G$: Given two finite free right $\Gamma$-sets $X, Y$ and a $\Gamma$-equivariant embedding $i:X\to Y$, the fibers of the projection $$\Hom_\Gamma(Y, D)\to\Hom_{\Gamma}(X, D)$$ should all have the same $\mu$-measure.

Let $a=|X/\Gamma|$ be the number of $\Gamma$-orbits in $X$ and $r=|Y/\Gamma|-a$. Then $i$ corresponds to adjoining $r$ free right $\Gamma$-orbits to $X$. Identifying the $\Gamma$-sets $X\cong \Gamma\times \{1, \ldots, a\}$, $Y\cong \Gamma\times \{1, \ldots,a+r\}$, the above projection becomes
\[
  \pi:D^{a+r}\longrightarrow D^a,\qquad
  (p_1,\ldots,p_{a+r})\longmapsto(p_1,\ldots,p_a).
\]
For any $G$-orbit $Q\subseteq D^a$, we have: $\pi^{-1}(Q)=Q\times D^r$.
We may now rewrite the required condition as follows:
The measure $\mu\left(\pi:Q\times D^r\to Q\right)$ should not depend on the choice of the $G$-orbit $Q\subseteq D^a$.

Yet the map $\pi:Q\times D^r\to Q$ is a base change of
$D^r\to\bullet_G$. Base-change invariance and multiplicativity
therefore give
\[
  \mu(Q\times D^r\longrightarrow Q)=\mu(D^r)=\mu(D)^r,
\]
which is independent of the choice of $Q$, as required. 
Thus every relative measure satisfies the criterion of
\cite[Proposition~6.11]{Snowden-Knop} and comes from the unique
degree function with $\delta=\mu(D)$. 

\end{proof}
Consider the $G$-set $\Omega$ on which $\Gamma$ acts trivially. Since the $G$-map $D\to\Omega$ has
fibres of cardinality $d$, we may alternatively parametrize the
absolute measures by
\[
  t=\mu(\Omega),\qquad \delta=\mu(D)=td.
\]
We write $\mu_t$ for the resulting measure ($t\in\kk$).

The comparison of categories
\cite[Theorem~7.5]{Snowden-Knop},
\cite[Theorem~9.5]{Nyobe-Savage}, together with the finite covers
in Lemma~\ref{lem:wp-open-subgroups}, gives
\[
  \Perm(G,\mu_t)^{\Kar}
  \simeq \operatorname{Par}(\Gamma,td)^{\Kar},
\]
where the right-hand side is the additive Karoubi envelope of
the group partition category in the sense of \cite{Nyobe-Savage}.

\begin{lemma}\label{lem:wreath-regular} The measure $\mu_t$ is regular iff $t \notin \mathbb{Z}_{\geq 0}$. All measures are quasi-regular.
\end{lemma}

\begin{proof}

By \cite[Corollary~9.6(b)]{Nyobe-Savage}, the category $\operatorname{Par}(\Gamma,td)^{\Kar}$ has a trivial negligible ideal exactly when
$td\notin d\mathbb Z_{\geq0}$.
Corollary~\ref{cor:trace-diagonal} and
Lemma~\ref{lem:ideals-kar} therefore imply
\[
  \mu_t\text{ is regular}
  \quad\Longleftrightarrow\quad t\notin\mathbb Z_{\geq0}.
\]

Let $A\subset \Omega$ be a finite subset and set $a=|A|$. The $U_A$-orbits in $\Omega$ are the $a$
fixed points and $\Omega\setminus A$. Hence
$(\mu_t|_{U_A})(\Omega\setminus A)=t-a$. 

Recall that $  U_A\cong\Gamma^{\Omega\setminus A}\rtimes
       \operatorname{Sym}(\Omega\setminus A)$ so $U_A$ is isomorphic to $G$ as an oligomorphic group. The measures on $U_A$ are therefore parametrized again by elements of $\kk$, and the measure
$\mu_t|_{U_A}$ therefore corresponds to parameter $t-a$. In particular, the above argument implies:
\begin{equation}\label{eq:wp-restricted-regularity}
      \mu_t|_{U_A}\text{ is regular}
  \quad\Longleftrightarrow\quad t-a\notin\mathbb Z_{\geq0}.
\end{equation}

Hence for $t\in \Z_{\geq 0}$, we may choose $a>t$ so that $\mu_t|_{U_A}$ would be regular for $A=\{1, \ldots, a\}$. Thus all measures $\mu_t$ on $G$ are quasi-regular. The singular measures are the $\mu_m$ with
$m\in\mathbb Z_{\geq0}$. 
%Since $\mu_m|_{U_{m+1}}$ is regular, (as the corresponding parameter is $-1$), every $\mu_t$ is quasi-regular.

\end{proof}

\begin{remark}
    Let $t \in \Z_{\geq 0}$. In this case the group partition categories admit a semisimple quotient by \cite{Nyobe-Savage}, and we conclude that any nilpotent endomorphism in $\Perm(G, \mu_t)$ has trace zero. The quasi-regularity of $\mu_t$ then implies that the wreath product categories admit an abelian envelope by \cite[§3.4]{HS-classical}, as previously proven with different methods in \cite{FHH}.
\end{remark}

%%%%%%%%%%%%%%%%%%%%%%%%%

\section{Products and categorical decompositions}

In this section we study representations of products of two oligomorphic groups.

\subsection{External tensor product of categories}

Let $\mathcal{T}_1$, $\mathcal{T}_2$ be two $\kk$-linear monoidal categories.

We define the \emph{naive external tensor product}, denoted $\mathcal{T}_1 \boxtimes_{\mathrm{naive}} \mathcal{T}_2$, as follows: the objects are formal pairs $(X_1, X_2)$ for $X_1 \in \mathcal{T}_1$ and $X_2 \in \mathcal{T}_2$, denoted as $X_1 \boxtimes X_2$. For objects $X_1 \boxtimes X_2$ and $Y_1 \boxtimes Y_2$, the $\Hom$-space is the tensor product of the underlying $\kk$-modules
\[ \Hom(X_1 \boxtimes X_2, Y_1 \boxtimes Y_2) = \Hom_{\mathcal{T}_1}(X_1, Y_1) \otimes_\kk \Hom_{\mathcal{T}_2}(X_2, Y_2)\]

and the composition of morphisms is defined by the bilinear extension of the rule
\[(f_1 \otimes_\kk f_2) \circ (g_1 \otimes_\kk g_2) = (f_1 \circ g_1) \otimes_\kk (f_2 \circ g_2).\]

This is in an obvious way a $\kk$-linear monoidal category where the tensor product on objects is given by
\[(X_1 \boxtimes X_2) \otimes_{\boxtimes} (Y_1 \boxtimes Y_2) = (X_1 \otimes_1 Y_1) \boxtimes (X_2 \otimes_2 Y_2)\]

and on morphisms by 
\[(f_1 \otimes_\kk f_2) \otimes_{\boxtimes} (g_1 \otimes_\kk g_2) = (f_1 \otimes_1 g_1) \otimes_\kk (f_2 \otimes_2 g_2).\] The monoidal unit is $\triv_{\boxtimes} = \triv_{\mathcal{T}_1} \boxtimes \triv_{\mathcal{T}_2}$.

We define additionally \[\mathcal{T}_1 \boxtimes \mathcal{T}_2 = (\mathcal{T}_1 \boxtimes_{\mathrm{naive}} \mathcal{T}_2)^{\Kar}.\]

\begin{lemma}\label{lem:external-kar}
There is a canonical monoidal equivalence \[ \mathcal{T}_1\boxtimes\mathcal{T}_2 \simeq (\mathcal{T}_1)^{\Kar}\boxtimes (\mathcal{T}_2)^{\Kar}.\]
Equivalently, \[(\mathcal{T}_1\boxtimes_{\mathrm{naive}}\mathcal{T}_2)^{\Kar}
\simeq \bigl(\mathcal{T}_1^{\Kar}\boxtimes_{\mathrm{naive}}\mathcal{T}_2^{\Kar}\bigr)^{\Kar}. \]
\end{lemma}

\begin{proof}
The fully faithful embeddings $\mathcal{T}_i\hookrightarrow (\mathcal{T}_i)^{\Kar}$ induce a fully faithful monoidal functor $\mathcal{T}_1\boxtimes_{\mathrm{naive}}\mathcal{T}_2 \to \mathcal{T}_1^{\Kar}\boxtimes_{\mathrm{naive}} \mathcal{T}_2^{\Kar}$. Passing first to finite direct sums and then to images of idempotents gives a fully faithful monoidal functor
\[ (\mathcal{T}_1\boxtimes_{\mathrm{naive}}\mathcal{T}_2)^{\Kar} \longrightarrow
\bigl(\mathcal{T}_1^{\Kar}\boxtimes_{\mathrm{naive}}\mathcal{T}_2^{\Kar}\bigr)^{\Kar}.
\] To prove essential surjectivity, it is enough to treat an object $$(X,e)\boxtimes(Y,f) \in \mathcal{T}_1^{\Kar}\boxtimes_{\mathrm{naive}}\mathcal{T}_2^{\Kar}$$ with $(X,e)\in \mathcal{T}_1^{\Kar}$ and $(Y,f)\in\mathcal{T}_2^{\Kar}.$ Here $X$ and $Y$ are finite direct sums of objects of $\mathcal{T}_1$ and $\mathcal{T}_2$ respectively, and $e\in \End_{\mathcal{T}_1}(X)$, $f\in \End_{\mathcal{T}_2}(Y)$ are idempotents. The object $(X,e)\boxtimes(Y,f)$ is the image of the idempotent $e\boxtimes f$ on $X\boxtimes Y$, and distributivity expresses $X\boxtimes Y$ as a finite direct sum of objects $X_i\boxtimes Y_j$ where $X_i\boxtimes Y_j\in \mathcal{T}_1\boxtimes_{\mathrm{naive}}\mathcal{T}_2$. Therefore $(X,e)\boxtimes(Y,f)$ lies in the essential image of $(\mathcal{T}_1\boxtimes_{\mathrm{naive}}\mathcal{T}_2)^{\Kar}$, as required.
\end{proof}

Suppose now that $\mathcal{T}_1$, $\mathcal{T}_2$ are Karoubian $\kk$-linear monoidal categories with a non-trivial ideal $\mathcal{I}$ in $\mathcal{T}_1$ and a non-trivial ideal $\mathcal{J}$ in $\mathcal{T}_2$. We define two new ideals in $\mathcal{T}_1\boxtimes\mathcal{T}_2$, given by
\[\mathcal{K}_1  = \mathcal{I} \boxtimes \mathcal{T}_2, \;\; \mathcal{K}_2  = \mathcal{T}_1 \boxtimes \mathcal{J}
\]

Formally, $\mathcal{K}_1$ is the ideal in $\mathcal{T}_1\boxtimes\mathcal{T}_2$ generated by all morphisms of the form $f \boxtimes g$, where $f \in \mathcal{I}$ and $g$ is any morphism in $\mathcal{T}_2$, and similarly for $\mathcal{K}_2$. The two ideals $\mathcal{K}_1$ and $\mathcal{K}_2$ are non-zero, proper, and distinct. Similarly one can form the ideal $\mathcal{I}\boxtimes \mathcal{J}$ in $\mathcal{T}_1\boxtimes\mathcal{T}_2$, which is distinct from both $\mathcal{K}_1$ and $\mathcal{K}_2$.

\begin{lemma} $\mathcal{K}_1$, $\mathcal{K}_2$ and $\mathcal{I}\boxtimes \mathcal{J}$ define distinct non-trivial tensor ideals in $\mathcal{T}_1 \boxtimes \mathcal{T}_2$.
\end{lemma}

\begin{proof}
First, we show that these ideals are non-trivial: $\mathcal{I}, ~\mathcal{J}\neq 0$ so $\mathcal{K}_1, ~\mathcal{K}_2,~\mathcal{I}\boxtimes \mathcal{J} \neq 0$, and $\id_{\triv\boxtimes\triv} \notin \mathcal{K}_1, \mathcal{K}_2, \mathcal{I}\boxtimes \mathcal{J}$, so none of them is $\mathcal{T}_1 \boxtimes \mathcal{T}_2$. Furthermore, these ideals are distinct: if $f_1\in \mathcal{I}$, $f_2\in \mathcal{J}$ are non-zero then $\id_{\triv}\boxtimes f_2\in \mathcal{K}_2$ but not in $\mathcal{K}_1$, $\mathcal{I}\boxtimes \mathcal{J}$, and $f_1\boxtimes \id_{\triv}\in \mathcal{K}_1$ but not in $\mathcal{K}_2$, $\mathcal{I}\boxtimes \mathcal{J}$.

Next, let $\mathcal{K}^0=\mathcal{I}\boxtimes_{\mathrm{naive}} \mathcal{J}$ in $\mathcal{T}_1\boxtimes_{\mathrm{naive}}\mathcal{T}_2$ be the ideal of morphisms spanned by $f_1\otimes_{\kk} f_2$ where $f_1\in \mathcal{I}, f_2\in\mathcal{J}$. The ideal $\mathcal{K}^0$ is a tensor ideal by the bilinearity of composition and of the tensor product. We first extend $\mathcal{K}^0$ to the additive envelope of $\mathcal{T}_1\boxtimes_{\mathrm{naive}}\mathcal{T}_2$ and obtain an ideal $\mathcal{K}^{\oplus}$ in the additive envelope. The morphisms in the additive envelope can be represented as matrices whose entries are morphisms in $\mathcal{T}_1\boxtimes_{\mathrm{naive}}\mathcal{T}_2$; such a matrix belongs to $\mathcal{K}^{\oplus}$ if and only if all its entries lie in $\mathcal{K}^0$. We now extend $\mathcal{K}^{\oplus}$ to an ideal $\mathcal{K}$ in  the additive Karoubi envelope $\mathcal{T}_1 \boxtimes \mathcal{T}_2$ by setting \[\mathcal{K}((X,e),(Y,f))=f~\mathcal{K}^{\oplus}(X,Y)~e.\] \cref{lem:ideals-kar} shows that $\mathcal{K}$ is a tensor ideal. An analogous construction applies to $\mathcal K_1$ and $\mathcal K_2$.
\end{proof}

Hence a monoidal category with only one non-trivial tensor ideal $\mathcal{N}$ cannot admit a product decomposition into two  monoidal categories with proper non-trivial tensor ideals.

%%%%%%%%%%%%%%%%%%

\subsection{Products of oligomorphic groups}\label{sec:products}
Let $G_1, G_2$ be two (pro-) oligomorphic groups with respective measures. We consider the group $G = G_1\times G_2$ with the product topology; it is a pro-oligomorphic group. Let $\mathcal{E}_{\times}$ be the product stabilizer class in $G_1\times G_2$, consisting of subgroups $U_1\times U_2$ with $U_i\leq G_i$ open.

By \cite[Section 3.4(c)]{HS-oligomorphic}, we have homomorphisms of measure rings (see \cref{ssec:olig_defs}):
$$ \Theta(G_1)\otimes \Theta(G_2) \xrightarrow{\sim} \Theta(G_1{\times} G_2, \mathcal{E}_{\times}) \to \Theta(G_1\times G_2)$$
where the first homomorphism is given by $[X]\otimes [Y]\mapsto [X\times Y]$. The first {isomorphism gives a correspondence} between relative measures on $(G_1{\times} G_2, \mathcal{E}_{\times})$ and products of measures  $ \mu_i$ on $G_i$: namely,
{given} a measure $\mu$ on $G_1\times G_2$, and finitary smooth $G_i$-set $X_i$ for $i=1,2$, we let
\[\mu_1(X_1)=\mu(X_1\times\bullet_{G_2}),\qquad
\mu_2(X_2)=\mu(\bullet_{G_1}\times X_2).\]
This gives us a measure $\mu_i$ for each $G_i$.
Conversely, let $\mu_i$ be measures on $G_i$. For each $i=1,2$, fix a finitary smooth $G_i$-set $X_i$. We define
\[(\mu_1\boxtimes \mu_2) \bigl(X_1\times X_2\bigr)=\mu_1(X_1)\mu_2({X_2}).\]

\begin{lemma}\label{lem:finite-index-product-subgroup}
Let $G_1,G_2$ be pro-oligomorphic groups and let $H\leq G_1\times G_2$ be open. Then there are open subgroups $N_i\leq G_i$ such that
\[N_1\times N_2\leq H,\qquad [H:N_1\times N_2]<\infty.\]
\end{lemma}

\begin{proof}
Since $H$ is open, it contains $U_1\times U_2$ for some open subgroups $U_i\leq G_i$. Put $H_i=\operatorname{pr}_i(H)$ and
\[N_1=\{g_1\in H_1\mid(g_1,1)\in H\}=\operatorname{pr}_1(\Ker \operatorname{pr}_2),\qquad
N_2=\{g_2\in H_2\mid(1,g_2)\in H\}=\operatorname{pr}_2(\Ker \operatorname{pr}_1).\]
Because $U_i\leq N_i$, each $N_i$ is open. If $h_i\in H_i$, choose $(h_1,h_2)\in H$ projecting to it. Conjugating $(n_1,1)$ by $(h_1,h_2)$ shows that $h_1n_1h_1^{-1}\in N_1$, and similarly for $N_2$; hence $N_i\triangleleft H_i$.

Goursat's lemma (see \cite[Theorem 2.1]{Bauer-Sen-Zvengrowski}) identifies $H$ with the fiber product of the quotient maps $H_i\to H_i/N_i$ over an isomorphism
\[ H_1/N_1\cong H_2/N_2=:Q.\] Since $H_i$ is open in the pro-oligomorphic group $G_i$, it is itself pro-oligomorphic. Since $N_i$ is open and normal,
\[|Q|=|H_i/N_i|=|N_i\backslash H_i/N_i|<\infty\] by Roelcke precompactness. The map $H\to Q$ has kernel $N_1\times N_2$, and therefore \[[H:N_1\times N_2]=|Q|<\infty.\]
\end{proof}
The following statements are direct consequences of \cref{lem:finite-index-product-subgroup} and \cite[Section 3.4(d)]{HS-oligomorphic}.
\begin{corollary}\label{cor:prod_is_large_stab_class}
Assume $\mathrm{char}~\kk=0$. Then the class $\mathcal E_{\times}$ is a large stabilizer class in the sense of \cite[Section 3.4(d)]{HS-oligomorphic}. In particular, $ \Theta(G_1\times G_2, \mathcal E_{\times})\otimes \mathbb{Q} \cong \Theta(G_1\times G_2)\otimes \mathbb{Q} $.
 
\end{corollary}
Combining this result with \cref{prop:saturating_stab_class}, we obtain:
\begin{corollary}\label{cor:prod_saturating}
       Assume $\mathrm{char}~\kk = 0$. Then any relative measure $\mu$ on $(G_1\times G_2, \mathcal{E}_{\times})$ extends uniquely to an absolute $\kk$-valued measure $\overline{\mu}$ on $G_1\times G_2$. 
       
       Furthermore, the functor $\Perm(G_1\times G_2, \mathcal{E}_{\times}, \mu) \to \Perm(G_1\times G_2, \overline{\mu})$ extends to an equivalence
       $$\Perm(G_1\times G_2, \mathcal{E}_{\times}, \mu)^{\Kar} \xrightarrow{\sim} \Perm(G, \overline{\mu})^{\Kar}.$$ 
\end{corollary}

We call this the \emph{product measure} and denote it by $\mu_1\boxtimes\mu_2$. It satisfies
\[ (\mu_1\boxtimes\mu_2)(X\times Y)=\mu_1(X)\mu_2(Y). \]

\begin{theorem}\label{thm:products}

Let $G_1, G_2$ be oligomorphic groups. 
Let $\mu$ be an $\cE_{\times}$-relative measure on $G_1\times G_2$, and let $$\mu_1=\mu(-\times \bullet_{G_2}),\mu_2=\mu(\bullet_{G_1}\times -)$$ be the measures on $G_1, G_2$ defined before. There exists a monoidal functor
\[
F:\Perm(G_1,\mu_1)\boxtimes_{\mathrm{naive}}\Perm(G_2,\mu_2)
  \longrightarrow \Perm(G_1\times G_2,\cE_{\times},\mu)
\]
inducing a monoidal equivalence
\[ F^{\Kar}:\bigl(\Perm(G_1,\mu_1)\boxtimes_{\mathrm{naive}}\Perm(G_2,\mu_2)\bigr)^{\Kar} \xrightarrow{\sim} \Perm(G_1\times G_2,\cE_{\times},\mu)^{\Kar}. \]
\end{theorem}

Equivalently, using \cref{lem:external-kar}, we obtain an equivalence
\[ \Perm(G_1,\mu_1)^{\Kar}\boxtimes \Perm(G_2,\mu_2)^{\Kar}
 \simeq \Perm(G_1\times G_2,\cE_{\times},\mu)^{\Kar}. \]
\begin{proof}
Let $\mathcal{T} = \Perm(G_1, \mu_1) \boxtimes_{\mathrm{naive}} \Perm(G_2, \mu_2)$ and define the functor $F: \mathcal{T} \to \Perm(G, \cE_{\times},\mu)$ as follows. Let $X_i, Y_i \in \Set{G_i}$ and let \[ F(\mathcal{C}(X_1) \boxtimes \mathcal{C}(X_2)) = \mathcal{C}(X_1 \times X_2).\]

To define the functor $F$ on morphisms, recall that the $\Hom$-space in $\mathcal{T}$ is 
\begin{align*}
    \Hom_{\mathcal{T}}(\mathcal{C}(X_1) \boxtimes \mathcal{C}(X_2), \mathcal{C}(Y_1) \boxtimes \mathcal{C}(Y_2)) &= \Hom_{G_1}(\mathcal{C}(X_1),\mathcal{C}(Y_1)) \otimes_\kk \Hom_{G_2}(\mathcal{C}(X_2), \mathcal{C}(Y_2))\\
&=\cC(Y_1 \times X_1)^{G_1} \otimes_\kk \cC(Y_2 \times X_2)^{G_2}.
\end{align*}  
%So any morphism from $\mathcal{C}(X_1) \boxtimes \mathcal{C}(X_2)$ to $\mathcal{C}(Y_1) \boxtimes \mathcal{C}(Y_2)$ is a finite $\kk$-linear combination of elementary tensors $A_1 \otimes A_2$, where $A_1 \in\cC(Y_1 \times X_1)^{G_1}$ is a $G_1$-invariant Schwartz function and $A_2 \in\cC(Y_2 \times X_2)^{G_2}$ is a $G_2$-invariant Schwartz function.

In the target category $\Perm(G, \mu)$, the $\Hom$-space is
\[ \Hom_G(\mathcal{C}(X_1 \times X_2), \mathcal{C}(Y_1 \times Y_2)) = \cC((Y_1 \times Y_2) \times (X_1 \times X_2))^{G_1 \times G_2}.\]

We define the functor $F$ on an elementary tensor $A_1 \otimes A_2$ as the function on $(Y_1 \times Y_2) \times (X_1 \times X_2)$ given by the product of their evaluations
$$F(A_1 \otimes A_2)\big((y_1, y_2), (x_1, x_2)\big) = A_1(y_1, x_1) \cdot A_2(y_2, x_2)$$ and extend this $\kk$-linearly to all morphisms in $\mathcal{T}$. This is clearly invariant under $G = G_1 \times G_2$, $\kk$-valued and of finitary support.

\medskip

We first verify that this assignment is functorial. For $A = A_1 \otimes A_2$ and $B = B_1 \otimes B_2$ we have
\begin{align*} \int_{Y_1 \times Y_2} B_1(z_1, y_1) B_2(z_2, y_2) A_1(y_1, x_1) A_2(y_2, x_2) \, d\mu(y_1, y_2) = \left( \int_{Y_1} B_1 A_1 \, d\mu_1 \right) \left( \int_{Y_2} B_2 A_2 \, d\mu_2 \right)\end{align*} by the product measure identity and Fubini's theorem (\cite[Corollary 4.9]{HS-oligomorphic}). As the tensor product in $\Perm(G,\mu)$ is the Cartesian product of $G$-sets, we conclude that $F$ is a monoidal functor.

\medskip

We claim that the map $F$ on the $\Hom$-sets induces an isomorphism
\[\cC(Y_1 \times X_1)^{G_1} \otimes_\kk\cC(Y_2 \times X_2)^{G_2} \xrightarrow{\sim}\cC((Y_1 \times Y_2) \times (X_1 \times X_2))^G.\]

Let $Z_i$ {be} a $G_i$-orbit in ${Y_i\times X_i}$ for $i=1,2$ and let $\one_{Z_i}$ be the respective indicator functions. Then $Z_1\times Z_2$ is a transitive $G$-subset of ${(Y_1\times Y_2)\times(X_1\times X_2)}$ and $$F(\one_{Z_1}\otimes \one_{Z_2}) = \one_{Z_1\times Z_2}$$ with the right hand side denoting the indicator function on $Z_1\times Z_2$. Since the corresponding indicator functions span the spaces  $\cC(Y_1 \times X_1)^{G_1},  \cC(Y_2 \times X_2)^{G_2}, \cC((Y_1 \times Y_2) \times (X_1 \times X_2))^G$ respectively, we conclude that $F$ induces the required isomorphism and so $F$ is a fully faithful functor.

\medskip

Finally, let us show that the functor $F^{\Kar}$ is essentially surjective. 

The category $\Perm(G_1\times G_2,\cE_{\times},\mu)^{\Kar}$ is generated, under taking finite direct sums and direct summands, by the objects $\cC(G/(U_1\times U_2))= \cC(G_1/U_1\times G_2/U_2)$ for $U_i\in \Sub({G_i})$. Such objects lie in the essential image of $F^{\Kar}$, so $F^{\Kar}$ is essentially surjective.

\end{proof}

%%%%%%%%%%%%%%%%%%%%%%

\section{Semidirect products with a finite group and equivariantization}\label{subsec:equivariantization}
Let $G$ be an oligomorphic group, and $\Gamma$ be a finite group acting on $G$ by continuous group automorphisms. Let $\kk$ be a ring. We will consider the semidirect product $\oG=G\rtimes \Gamma$. 

\subsection{Twisting of \texorpdfstring{$G$-sets}{G-sets}}

Given $X\in \Set{G}$ and $\gamma\in \Gamma$, we define $X^{\gamma}$ to be a $G$-set whose underlying set is $X$, but the action of $G$ is twisted by the automorphism $\gamma$ of $G$: that is, for $g\in G$, and $x\in X^{\gamma}$, $g\cdot_\gamma x=\gamma^{-1}(g).x$. Since the automorphisms $\gamma: G\to G$ are continuous, they are also open and so they send open subgroups to open subgroups. Thus $X^\gamma$ is again a smooth finitary $G$-set. 

It is immediate that for any $X, Y\in \Set{G}$ and $\gamma\in \Gamma$ we have a $G$-isomorphism $(X\times Y)^\gamma\cong X^\gamma\times Y^\gamma$. For a function $\underline{f}\in \cC(Y\times X)^G$, we denote by $\underline{f}^{\gamma}\in \cC(Y^\gamma\times X^\gamma)$ the same function on the underlying sets. For $f\in \Hom_G(\cC(X), \cC(Y))$ corresponding to $ \underline{f}\in \cC(Y\times X)^G$, we denote by $f^\gamma\in \Hom_G(\cC(X^\gamma), \cC(Y^\gamma))$ the map corresponding to $\underline{f}^{\gamma}$.

\begin{definition}
    Let $\mu$ be an (absolute) $\kk$-valued measure on $G$. We say that $\mu$ is \emph{$\Gamma$-invariant} if for any open subgroup $U$ of $G$, $X\in \Set{U}$ and $\gamma\in \Gamma$ we have $\mu_{U}(X) =\mu_{\gamma(U)}(X^\gamma)$.  
\end{definition}
\subsection{Semidirect product}
Let $\oG=G\rtimes \Gamma$. This is again an oligomorphic group, and $G$ is a normal subgroup of $\oG$. Therefore \cite[\S3.4(e)]{HS-oligomorphic} implies $\Theta(G)_{\Gamma} = \Theta(\oG)$ where $(-)_{\Gamma}$ denotes the coinvariants. 

In other words, given a $\Gamma$-invariant measure $\mu$ on $G$, we have a unique corresponding measure $\overline{\mu}$ for $\oG$ (and all the measures for $\oG$ are obtained in this way). The measure $\overline{\mu}$ is defined as follows: for any open subgroup $\overline{U}\leq \overline{G}$ and any $\overline{X}\in \Set{\oU}$, we let \[\label{eq:equivariantization-measure}
 \overline{\mu}_{\overline{U}}(\overline{X})=\mu_{\overline{U}\cap G}\bigl(\overline{X}|_{\overline{U}\cap G}\bigr).
\]

\begin{lemma}\label{lem:semidirect_restriction_on_sets}
    Let $\overline{X}$ be a transitive $\oG$-set, and let $X\in G\backslash  \overline{X}$. Then for each $\gamma\in \Gamma$, $\gamma X\cong X^{\gamma}$ is a $G$-orbit in $\overline{X}$ and we have an isomorphism of sets $\overline{X}|_G\cong \bigcup_{\gamma\in \Gamma} X^{\gamma}$.
\end{lemma}
 The subgroup $\Gamma \leq \oG$ acts on $\overline{X}$ by permuting the above $G$-orbits.
\begin{proof}
    Let $x\in X\subseteq \overline{X}$ and $\gamma\in \Gamma$. The transitive $G$-set $G.(\gamma.x)$ is clearly isomorphic to $X^\gamma$, so $X^{\gamma} $ form $G$-orbits in $\overline{X}$ and the group $\Gamma$ permutes these orbits. Since $\overline{X}$ is a transitive $\oG$-set, we conclude that $\overline{X}|_G\cong \bigcup_{\gamma\in \Gamma} X^{\gamma}$.
\end{proof}

The inclusion of $G \leq \oG$  gives a faithful SM functor
\[
 \operatorname{Res}^{\oG}_{G}:\Perm(\oG,\overline{\mu})\longrightarrow\Perm(G,\mu).
\]
We will use a subscript when denoting spaces of Schwartz functions in $\Perm(\oG,\overline{\mu}),~\Perm(G,\mu)$ to indicate which group is acting on the space.

The functor $\operatorname{Res}^{\oG}_{G}$ sends $\mathcal C_{\oG}(\oX)$ to $\mathcal C_{G}(\oX|_{G})$, and a morphism $\mathcal C_{\oG}(\oX)\to \mathcal C_{\oG}(\overline{Y}) $ represented by a $\oG$-invariant Schwartz function $a:\overline{Y}\times \oX\to\kk$ to the morphism $\mathcal C_{G}(\oX)\to \mathcal C_{G}(\overline{Y}) $ represented by the same function (now regarded as a $G$-invariant function). 

\subsection{Equivariantization}

\medskip

Let $\mathcal T$ be a $\kk$-linear symmetric monoidal category with an action of a finite group $\Gamma$ by monoidal functors: that is, for every $\gamma\in \Gamma$ there is a monoidal functor $T_{\gamma}:\mathcal T\to \mathcal T$ so that $T_1\cong \id $ and $T_\sigma T_\gamma \cong T_{\sigma \gamma} $ for any $\gamma, \sigma \in \Gamma$.
\begin{definition}[{\cite[Definition~2.7.2]{EGNO}}]\label{def:equivariantization}
The \emph{equivariantization} $\mathcal T^{\Gamma}$ has objects $(M,u)$, where $M\in\mathcal T$ and $u$
is a collection of maps \[
 u_\gamma:T_\gamma M\xrightarrow{\sim}M\; \text{ for } \gamma \in \Gamma, \qquad \text{ such that } \qquad
 u_1=\id_M,\qquad u_{\gamma\tau}=u_\gamma T_\gamma(u_\tau).
\]
A morphism $f:(M,u)\to(N,v)$ is a morphism $f:M\to N$ in $\mathcal T$ satisfying
\[\forall \; \gamma~ \in ~\Gamma, \qquad
 f u_\gamma=v_\gamma T_\gamma(f).
\]
\end{definition}

As in \cite[Definition~2.7.2]{EGNO} the category $\mathcal T^{\Gamma}$ is monoidal, and the obvious forgetful functor $\For:\mathcal T^{\Gamma}\to\mathcal T$ which sends $(M,u)$ to $M$ is a SM functor.

In our setting, we define the functor $T_{\gamma}: \Perm(G,\mu)\to \Perm(G,\mu)$ via $\cC(X)\mapsto \cC(X^{\gamma})$ and $f\mapsto f^{\gamma}$.
\begin{lemma}
Let $\gamma\in \Gamma$. The functors $T_{\gamma}$, $\gamma\in \Gamma$ are SM functors satisfying $T_1\cong \id $ and $T_\sigma T_\gamma \cong T_{\sigma \gamma} $ for any $\gamma, \sigma \in \Gamma$.
\end{lemma}
\begin{proof}
    The fact that $T_{\gamma}$ is SM follows from the fact that $(X\times Y)^\gamma\cong X^\gamma\times Y^\gamma$.
    The equivalence $T_{1}\cong \id $ is also straightforward. To prove the remaining statement, consider the maps $(X^\gamma)^\sigma \to X^{\sigma \gamma}$ given by the identity on the underlying set. These are $G$-maps (by the definition of group action of $\Gamma$ on $G$), so they define a natural isomorphism $T_\sigma T_\gamma \to T_{\sigma \gamma} $.
\end{proof}
We thus obtain an equivariantization $\Perm(G,\mu)^\Gamma$ of $\Perm(G,\mu)$ and denote the forgetful functor by $\For:\Perm(G,\mu)^\Gamma \to \Perm(G,\mu)$.

\subsection{Equivalence of categories}
 We will describe now the category $\Perm(\oG,\overline{\mu})$ as an equivariantization of the category $\Perm(G,{\mu})$ with respect to the $\Gamma$-action.

\begin{proposition}\label{prop:equivariantization-product}
Assume $\kk$ is a field such that $\operatorname{char}\kk\nmid |\Gamma| $.
We have a $\kk$-linear SM equivalence
\[
 \mathcal R:\Perm(\oG,\overline{\mu})^{\Kar}\xrightarrow{\sim}\left(\Perm(G, \mu)^{\Gamma}\right)^{\Kar}.
\]
\end{proposition}
\begin{remark}
    We have an equivalence $\left(\Perm(G, \mu)^{\Gamma}\right)^{\Kar}\cong \left(\Perm(G, \mu)^{\Kar}\right)^{\Gamma}$ due to \cite[Lemma~2.16]{Sun-equivariantization}. 
\end{remark}
\begin{proof} 

Our aim is to construct a SM functor 
\[ \mathcal R_0:\Perm(\oG,\overline{\mu})\longrightarrow\Perm(G, \mu)^{\Gamma} \]
such that $\operatorname{Res}^{\oG}_{G}\cong\For\circ \mathcal{ R}_0$ and $\mathcal{ R}_0$ will yield the desired equivalence $\mathcal{R}$ once it is extended to the Karoubi envelope.

Recall that \cite[\S7.7, Propositions~7.20 and~7.21(a)]{HS-oligomorphic} attach to any $\oG$-equivariant morphism $f:\oX\to \overline Y$ two maps $f_*: \mathcal{C}_{\oG}( \overline X) \to \mathcal{C}_{\oG}( \overline Y)$ and $f^*:\mathcal{C}_{\oG}( \overline Y) \to \mathcal{C}_{\oG}( \overline X)$. Let $\oX$ be a finitary smooth $\oG$-set. Regard $\gamma\in \Gamma$ as an element in $\oG$. Then the bijection
\[
 b_\gamma^\oX: \oX^{\gamma}\longrightarrow \oX,\qquad
 x\longmapsto\gamma x,
\]
is $G$-equivariant (cf. \cref{lem:semidirect_restriction_on_sets}), because for $g\in G$ one has
\[
 b_\gamma^\oX(g\cdot_{\gamma}x)
 =\gamma\bigl((\gamma^{-1}.g).x\bigr)
 =g(\gamma x)
 =g\,b_\gamma^\oX(x).
\]

Recall that $ T_\gamma\mathcal C_{G}(\oX)=\cC_G(\oX^{\gamma})$. We set
\[
 u_\gamma^\oX=(b_\gamma^\oX)_*:
 T_\gamma\mathcal C_{G}(\oX)
 \xrightarrow{\sim}\mathcal C_{G}(\oX), \;\;  (u_\gamma^\oX\psi)(x)=\psi(\gamma^{-1}x) \; \text{ for } \;  {\psi\in \cC_G(\oX^{\gamma})}.  
\]
The map $u_\gamma^\oX$ is invertible with inverse $((b_\gamma^\oX)^{-1})_*$. We now define

\[
 \begin{aligned}
 \mathcal R_0:\Perm(\oG,\overline{\mu})&\longrightarrow\Perm(G, \mu)^{\Gamma},\\
 \mathcal C_{\oG}(\oX)&\longmapsto
   (\mathcal C_{G}(\oX),(u_\gamma^\oX)_{\gamma\in \Gamma})
 \end{aligned}
\]
which acts as restriction on morphisms. We first need to check that this is well-defined, i.e. that $ (\mathcal C_{G}({\oX}),(u_\gamma^\oX)_{\gamma\in \Gamma})$ define indeed equivariant objects and that restriction of morphisms from $\Perm(\oG,\overline{\mu})$ obey the rules of Definition \ref{def:equivariantization}.

For $\gamma,\tau\in \Gamma$, consider the composite
\[
 T_\gamma T_\tau \oX
 \xrightarrow{\,T_\gamma(b_\tau^\oX)\,}T_\gamma \oX
 \xrightarrow{\,b_\gamma^\oX\,}\oX
\]
which sends $x$ first to $\tau x$ and then to $\gamma\tau x$. Since $T_\gamma T_\tau \cong T_{\gamma\tau}$, this equals $b_{\gamma\tau}^\oX$. Since the formation of $(-)_*$ respects composition and commutes with twisting, we get
\[ u_1^\oX=\id,\qquad
 u_\gamma^\oX\circ T_\gamma(u_\tau^\oX)=u_{\gamma\tau}^{{\oX}}.
\]
Hence $(\mathcal C_{G}(\oX),(u_\gamma^\oX)_{\gamma\in \Gamma})$
is an object of $\Perm(G, \mu)^{\Gamma}$.

Let $f:\mathcal C_{\oG}(\overline X)\to\mathcal C_{\oG}(\overline Y)$ be represented by a $\oG$-invariant Schwartz function $f(y,x)$. Composition with the above relabeling maps gives
\[
 (f\circ u_\gamma^\oX)(y,x)
 =f(y,\gamma x)
 =f(\gamma^{-1}y,x)
 =(u_\gamma^{\overline Y}\circ T_\gamma(f))(y,x),
\]
where the middle equality follows from $\oG$-invariance. Hence the restricted morphism is a morphism in $\Perm(G, \mu)^{\Gamma}$.

The functoriality of $\mathcal{R}_0$ follows from that of ordinary restriction $\Res^{\oG}_G$, since the underlying objects and morphisms are unchanged. Finally,
$b_\gamma^{\overline X\times \overline Y}=b_\gamma^\oX\times b_\gamma^{\overline Y}$ and
$b_\gamma^\bullet=\id$, so the tensor identifications of restriction respect the equivariant structures. Therefore $\mathcal R_0$ is a SM functor. By construction, the diagram
\[
\begin{tikzcd}[column sep=large,row sep=large]
 \Perm(\oG,\overline{\mu}) \arrow[r,"\mathcal R_0"]
 \arrow[dr,swap,"\operatorname{Res}^{\oG}_{G}"]
 & \Perm(G, \mu)^{\Gamma}\arrow[d,"\operatorname{For}"]\\
 &\Perm(G, \mu)
\end{tikzcd}
\]
commutes. We claim that the functor $\mathcal R_0$ is fully faithful. A function $a\in \cC(\overline Y\times \overline X)^G$ defines a morphism $\mathcal{R}_0\cC_{\oG}(\oX)\to \mathcal{R}_0\cC_{\oG}(\overline Y)$ if and only if 
\[
 a(\gamma y,\gamma x)=a(y,x)\qquad(\gamma\in \Gamma).
\]
Since $G$ and $\Gamma$ generate $\oG$, this is just $\oG$-invariance of $a$. Thus $$\mathcal R_0~:~\Hom_{\oG}(\cC_{\oG}(\oX),\cC_{\oG}(\overline Y))\;\longrightarrow\; \Hom_{\Perm(G, \mu)^{\Gamma}}(\mathcal{R}_0\cC_{{\oG}}(\oX),\mathcal{R}_0\cC_{{\oG}}(\overline Y)) $$ is a bijection and $\mathcal R_0$ is fully faithful.

The universal property of the Karoubi completion now extends $\mathcal R_0$ to a fully faithful SM functor
\[
 \mathcal R=\mathcal R_0^{\Kar}:\Perm(\oG,\overline{\mu})^{\Kar}
 \longrightarrow(\Perm(G, \mu)^{\Gamma})^{\Kar}.
\]

It remains to show that $\mathcal{R}$ is essentially surjective.

Let $I:\Perm(G, \mu)\to\Perm(G, \mu)^{\Gamma}$ be the additive induction functor defined by
\[
 I(M)=(\bigoplus_{\gamma\in \Gamma}T_\gamma M, u), \qquad I(f)=\bigoplus_\gamma T_\gamma(f)\]
 where $u_{\tau}: T_{\tau}\left(\bigoplus_{\gamma\in \Gamma}T_\gamma M\right)\to \bigoplus_{\gamma\in \Gamma}T_\gamma M$ is the obvious isomorphism. 
This functor is both left and right adjoint to $\For$, and extends to a functor $I:\Perm(G, \mu)^{\Kar}\to\left(\Perm(G, \mu)^{\Gamma}\right)^{\Kar}$.

The standard induction--forgetful splitting \cite[Lemma~2.14(1),(3)]{Sun-equivariantization} makes every object $(M,u)\in \Perm(G, \mu)^{\Gamma}$ a retract of $I(M)$, and by definition, every object of $\left(\Perm(G, \mu)^{\Gamma}\right)^{\Kar}$ is a {retract} of some object $(M,u)\in \Perm(G, \mu)^{\Gamma}$. So it suffices to show that every $I(\cC_G(X))$, $X\in \Set{G}$, is a retract of an object in the image of $\mathcal{R}_0$.

Consider the induced set
\[
 \overline{X}=\oG\times_{G}X
   =(\oG\times X)/\bigl((gh,x)\sim(g,hx),\ h\in G\bigr).
\]
This is a finitary smooth $\oG$-set (where $\oG$ acts on the first coordinate \cite[\S2.4]{HS-oligomorphic}). When we restrict we get
\[
 \overline{X}|_{G}\cong\coprod_{\gamma\in \Gamma} X^\gamma,
\]

We will denote by $[\gamma,x]$ the copy of an element $x\in X$ in the $X^\gamma$. Multiplication by $\tau\in \Gamma$ sends $[\gamma,x]$ to $[\tau\gamma,x]$. We obtain
\[
 I\bigl(\mathcal C_{G}(X)\bigr)\cong\mathcal R_0\bigl(\mathcal C_{\oG}(\overline{X})\bigr).
\]
Thus $I(\cC_G(X))$, $X\in \Set{G}$, is a retract of an object in the image of $\mathcal R_0$, as required.
\end{proof}
\begin{remark}
The functor $\mathcal{R}_0$ need not be an equivalence. 

Indeed, let $G=\{1\}$ and $\Gamma=\Z/3$, so $\oG=\Z/3$. In this case $\mu, \overline{\mu}$ are just counting measures. Then $\Perm(\oG, \overline{\mu})$ is the category $\Perm(\Z/3)$ of permutation modules for $\Z/3$, while $\Perm(G, {\mu})^{\Gamma}=\left(\Vect\right)^{\Z/3}\cong \Rep_{\kk}(\Z/3)$. The functor $\mathcal{R}_0:\Perm(\oG, \overline{\mu})\to \Perm(G, {\mu})^{\Gamma}$ is then $\Perm(\Z/3)\to \Rep_{\kk}(\Z/3)$, which is not essentially surjective even when $\mathrm{char}~\kk=0$. 

    This example shows that the assumption that $\operatorname{char}\kk\nmid |\Gamma| $ cannot be dropped either: if $\mathrm{char}~\kk=3$, then the Karoubi envelope of $\Perm(\Z/3)$ is not $\Rep_{\kk}(\Z/3)$.
\end{remark}

%%%%%%%%%%%%%%%%%%%
\subsection{Main example: wreath product with \texorpdfstring{$S_d$}{the symmetric group}}
    A natural example for the above setting is the symmetric group $S_d$ acting on a product of $d$ copies of an oligomorphic group. Assume $\mathrm{char} ~\kk=0$.
    
 Let $G_0$ be a pro-oligomorphic group with an absolute $\kk$-valued measure $\mu_0$, and set
\[
\nu=\mu_0^{\boxtimes d},\qquad \mathcal T =\Perm(G_0,\mu_0)^{\Kar}
\]
The symmetric group $S_d$ acts on $G_0^d$ by permuting the factors. For $\sigma\in S_d$, write $\alpha_\sigma(g)=\sigma g\sigma^{-1}$ for the corresponding automorphism of $G_0^d$.
 Let $\mathcal T'=\Perm(G_0^d,\nu)^{\Kar}$. We have an equivalence of categories, due to \cref{lem:external-kar}: $
E:\mathcal T^{\boxtimes d}~\xrightarrow{\sim} ~  \mathcal T'
$.

The general construction above allows us to define an action of $S_d$ on $\Perm(G_0^d,\nu)$, given by $\kk$-linear SM autoequivalences
\[
 T_\sigma:\Perm(G_0^d,\nu)\longrightarrow\Perm(G_0^d,\nu).
\] Similarly, we may define the action of $S_d$ on the category $\mathcal T^{\boxtimes d}$ by monoidal $\kk$-linear SM autoequivalences \[  T'_\sigma:\mathcal T^{\boxtimes d}\longrightarrow \mathcal T^{\boxtimes d}, \qquad
 T'_\sigma (M_1\boxtimes\cdots\boxtimes M_d)
 = M_{\sigma^{-1}(1)}\boxtimes\cdots\boxtimes M_{\sigma^{-1}(d)},
\]
for each $\sigma \in S_d$.

\begin{lemma}\label{lem:twisting}
There exists a natural transformation $T_\sigma\circ E~\xrightarrow{\sim} ~E\circ T'_{\sigma}$. 

\end{lemma}

\begin{proof}
 For finitary smooth $G_0$-sets $X_1,\ldots,X_d$, the coordinate permutation
\[
 \begin{aligned}
 c_\sigma:(X_1\times\cdots\times X_d)^{\sigma}
     &\xrightarrow{\sim}
       X_{\sigma^{-1}(1)}\times\cdots\times X_{\sigma^{-1}(d)},\\
 (x_1,\ldots,x_d)&\longmapsto
       (x_{\sigma^{-1}(1)},\ldots,x_{\sigma^{-1}(d)})
 \end{aligned}
\]
is $G_0^d$-equivariant: the $i$-th coordinate of
$c_\sigma(g\cdot_{T_\sigma}x)$ is $g_i x_{\sigma^{-1}(i)}$ for $g\in {G_0^d}$.
Under the product construction of \cref{thm:products}, this corresponds to the permutation of the external tensor factors. The statement for arbitrary $M_i\in\mathcal T$ follows from extending to finite sums and retracts.
    
\end{proof}
Thus the functors $T'_{\sigma}, ~\sigma\in S_d$ define an action of $S_d$ on {$\mathcal T^{\boxtimes d}$ by $\kk$-linear SM autoequivalences, and $E$ is an $S_d$-equivariant equivalence}.

Now let us consider the semidirect product $G_0^d\rtimes S_d$.
The product measure $\nu=\mu_0^{\boxtimes d}$ on $G_0^d$ is $S_d$-invariant, so it extends uniquely to a measure $\overline{\nu}$ on $G_0^d\rtimes S_d$ by \cite[\S3.4(e)]{HS-oligomorphic}.
\cref{prop:equivariantization-product} now implies:
\begin{corollary}
      We have $\kk$-linear SM equivalences
\[
\Perm(G_0^d\rtimes S_d,\overline{\nu})^{\Kar}~\cong~(\mathcal T')^{S_d} ~\cong~ (\mathcal T^{\boxtimes d})^{S_d}.
\]
  
\end{corollary}

\part{Classification of tensor ideals}

 \section{Classification of tensor ideals}

\begin{definition}

Given $2$ transitive $G$-sets $X, Y\in \mathbf{S}(G)^{tran}$, we write $X\succeq_\mu Y$ if there exist $G$-maps $f_X:Z\to X$, $f_Y:Z\to Y$ in $\mathbf{S}(G)$ such that $\mu(f_Y)\neq0$, where $\mu(f_Y)$ is the measure of any fiber of $f_Y$ (recall that these fibers all have the same measure).

\end{definition}
\begin{remark}

\mbox{}

    \begin{itemize}[itemindent=0pt,leftmargin=*]
    \item Given a $G$-map $Y\to X$ of transitive $G$-sets, we have: $X\succeq_\mu Y$.
        \item {Since} $G$ is a maximal open subgroup, $\bullet_G$ is a maximal transitive $G$-set.
    \end{itemize}
\end{remark}
\begin{lemma}\label{lem:succeq_mu}
    Let $X, Y\in \mathbf{S}(G)^{tran}$. The following conditions are equivalent:
    \begin{enumerate}
        \item $X\succeq_\mu Y$.
        \item There exists $Z \in \mathbf{S}(G)^{tran}$ and $G$-maps $f_X:Z\to X$, $f_Y:Z\to Y$ such that $\mu(f_Y)\neq 0$.
        \item There exists $Z\subseteq X\times Y$ in $\mathbf{S}(G)^{tran}$ such that $\mu(\pi_Y\rvert_Z)\neq 0$, where $\pi_Y: X\times Y\to Y$ is the projection map.
    \end{enumerate}
\end{lemma}
\begin{proof}
Clearly, $(3) \implies (2)\implies (1)$ so we will prove that $(1)\implies (2)\implies  (3)$.

    Indeed, let $Z\in {\mathbf{S}(G)}$ with $G$-maps $f_X:Z\to X$, $f_Y:Z\to Y$ such that $\mu(f_Y)\neq0$. 
            
    Let $Z_1, \ldots, Z_k$ be the $G$-orbits in $Z$, so that $Z=Z_1 \sqcup \ldots{\sqcup} Z_k$. Then for every $i$, we have $G$-maps $f_X\rvert_{Z_i}:Z_i\to X$, $f_Y\rvert_{Z_i}:Z_i\to Y$ and $\mu(f_Y)= \sum_{i=1}^k\mu(f_Y\rvert_{Z_i})$. In order to have $\mu(f_Y)\neq 0$, at least one of the summands must be non-zero. Hence $(1)\implies (2)$. 
    
    Now, let $Z\in \mathbf{S}(G)^{tran}$ with $G$-maps $f_X:Z\to X$, $f_Y:Z\to Y$ such that $\mu(f_Y)\neq0$.  We have a $G$-map $f:Z\to X\times Y$ so that $f_Y=\pi_Y\circ f$. Let $Z'\subseteq X\times Y$ be the image of $f$; this is a $G$-transitive set. Then $f_Y=\pi_Y\rvert_{Z'}\circ f $ and thus $\mu(\pi_Y\rvert_{Z'})\mu( f:Z\to Z') = \mu(f_Y)\neq 0$. We conclude that $\mu(\pi_Y\rvert_{Z'})\neq 0$. Hence $(2)\implies (3)$.
\end{proof}
\begin{remark}\label{rmk:order_on_subgroups}
      Let $U, V\leq G$ be open subgroups. We say that $V \preceq_\mu U$ if and only if there exists $g\in G$ such that $\mu(V/(V\cap gUg^{-1}))\neq 0$. 
      
      By the above Lemma, we have: $V \preceq_\mu U $ iff $G/V\preceq_\mu G/U$. 
\end{remark}

We defined a relation on transitive $G$-sets, and this relation is reflexive. Let us show that it is transitive.

\begin{lemma}
    Let $X, Y, Y'$ be transitive $G$-sets. 
    If $X\succeq_\mu Y$, $Y \succeq_\mu Y'$ then $X\succeq_\mu Y'$.
  
\end{lemma}

\begin{proof}

  %  Let us prove the corresponding statement for open subgroups. Note that in the definition of the relation $\preceq_\mu$, it is enough to show that there exists an open subgroup $W<V\cap U^g$ for which $\mu(V/W)\neq 0$ (since $\mu(V/W) = \mu(V/V\cap gUg^{-1})\mu(V\cap gUg^{-1}/W)$).    Assume that $ V\preceq_\mu U$ and $V'\preceq_\mu V$. Let $g, h\in G$ be such that $\mu(V/V\cap gUg^{-1})\neq 0$, $\mu(V'/V'\cap hVh^{-1})\neq 0$.     Let $$ W=V'\cap hVh^{-1} \cap (hg)U(hg)^{-1} = V'\cap h(V\cap gUg^{-1})h^{-1} \subseteq V' \cap (hg)U(hg)^{-1}.$$ We need to show that $\mu(V'/W)\neq 0$.
Assume that $X\succeq_\mu Y$, $Y \succeq_\mu Y'$.
We have a commutative diagram in $\mathbf{S}(G)$:
$$\xymatrix{&Z'\times_Y Z \ar_{\pi_{Z'}}[d]\ar^-{\pi_Z}[r] &Z \ar_{f_Y}[d]\ar^{f_X}[r] \ar[r] &X\\
&Z' \ar_{f'_{Y'}}[d]\ar^{f'_Y}[r] &Y & \\
&Y' & &}$$
so that $Z, Z'$ are transitive and $\mu(f_Y), \mu(f'_{Y'})\neq 0$.
Then $\mu(f'_{Y'}\circ \pi_{Z'}) = \mu(f'_{Y'})\mu( \pi_{Z'})$. Now, the fiber of $  \pi_{Z'} $ is isomorphic as a $U$-set to the fiber of $f_Y$ for a common subgroup of definition $U\in \Sub(G)$. So  $\mu( \pi_{Z'}) \neq 0$. Thus $\mu(f'_{Y'}\circ \pi_{Z'})\neq 0$, hence $X\succeq_\mu Y'$.
      
\end{proof}

\begin{definition}
    Let $X, Y\in \mathbf{S}(G)^{tran}$. We say that $X, Y$ are \emph{$\mu$-comparable}, and write $X\sim_{\mu}Y$ if $X\succeq_\mu Y, Y\succeq_\mu X$. 
    
\end{definition}

\begin{definition}
Let $\mathcal{P}^G_\mu$ be the poset of equivalence classes of $\mu$-comparable transitive $G$-sets, equipped with the induced partial order $[X] \preceq_\mu [Y]$. 

\end{definition}
\begin{remark}
    The partial order on $\mathcal{P}^G_\mu$ is given explicitly by: $[X] \succeq_\mu [Y]$ if there exists a $G$-map $Z\rightarrow X$ where $Z\in [Y]$.
\end{remark}

\begin{notation}
    Let $X\in \mathbf{S}(G)$. We denote by $\eta_X:\triv \to \mathcal{C}(X)$ the $G$-map corresponding to the constant function $1$ on $X$. 
\end{notation}

\begin{definition} 

\mbox{}

\begin{itemize}[itemindent=0pt,leftmargin=*]
    \item A subset $D \subseteq \mathcal{P}^G_\mu$ is \emph{downward-closed} if $[Y] \in D$ and $[Z] \preceq_\mu [Y]$ implies $[Z] \in D$. 
    \item Let $\Down(\mathcal{P}^G_\mu)$ be the poset of downward-closed subsets of $\mathcal{P}^G_{\mu}$. The poset $\Down(\mathcal{P}^G_\mu)$ is ordered by subset inclusion $\subseteq$.
\end{itemize}
\end{definition}

\begin{theorem}\label{thm:classification}
    The poset of tensor ideals of $\Perm(G, \mu)$ is in bijection with the poset of downward-closed subsets of $\mathcal{P}^G_\mu$. A tensor ideal $\mathcal{I}\subset\Perm(G, \mu)$ corresponds to $D\subseteq \mathcal{P}^G_\mu$ if {and} only if $\mathcal{I}$ is generated by $\eta_X$ where $[X]\in D$.

    Explicitly, for $D\in \Down(\mathcal{P}^G_\mu)$, the corresponding ideal $\mathcal{I} $ is given by $$\mathcal{I}(\cC(X), \cC(Y))=\Span\{\one_O~:~ O\in G\backslash (Y\times X) \text{ and }[O]\in D\}, \qquad X, Y\in \Set{G}$$
    where $\one_O\in \Hom_{\Perm(G, \mu)}(\cC(X), \cC(Y))$ corresponds to the indicator function on the orbit $O$ under the isomorphism $\Hom_{\Perm(G, \mu)}(\cC(X), \cC(Y))\cong \cC(Y\times X)^G$.
\end{theorem}
\begin{example}\label{ex:largest-class}
    The class $[\bullet_G]$ contains transitive $G$-sets which have non-zero measure. It is maximal with respect to the order $\preceq_\mu$ and the downward-closed subset it generates is the entire $\mathcal{P}^G_\mu$. The corresponding tensor ideal is the entire category $\Perm(G, \mu)$.
\end{example}

\begin{example}\label{ex:negligible} Under this correspondence $\mathcal{N}$ corresponds to the downward-closed set $$\mathcal{P}^G_\mu\setminus\{[\bullet_G]\}=\{[X]\in \mathcal{P}^G_{\mu} ~:~\mu(X)=0\}.$$
 Indeed, $\mathcal{N}$ is the unique maximal proper tensor ideal, and its downward-closed subset is the complement of the greatest equivalence class. 
\end{example}

\begin{proof}[Proof of \cref{thm:classification}]
    The proof of the theorem is done in several steps:
    \begin{enumerate}[itemindent=0pt,leftmargin=*]
        \item First of all, consider the additive functor $\Hom(\triv, -): \Perm(G, \mu)\to \Vect$.  By \cref{thm:Coulembier} (see \cite[Theorem 3.1.1]{Coulembier-ideals}), there exists a bijection between the poset of tensor ideals of $\Perm(G, \mu)$ and the poset of additive subfunctors of $\Hom(\triv, -)$. The bijection is given by $\mathcal{I} \longmapsto \mathcal{I}(\triv, -)$.
        \item To understand the poset of subfunctors of $\Hom(\triv, -)$, we use \cref{lem:subfunc_Hom_from_1_gen_by_eta} below: such an additive subfunctor $H\subseteq \Hom(\triv, -)$ is uniquely determined by the collection of $\eta_X$ (here $X\in \mathbf{S}(G)^{tran}$) for which $H(\mathcal{C}(X))\neq 0$; that is, $H(\mathcal{C}(X))\neq 0$ if and only if $H(\mathcal{C}(X))=\Span\{\eta_X\}$. 
        In other words, the maps $\eta_X$ for which $H(\mathcal{C}(X))\neq 0$ generate the corresponding tensor ideal.
        \item For an additive subfunctor $H\subseteq \Hom(\triv, -)$, let $$D_{H}=\{[X]\in \mathcal{P}^G_\mu~: ~H(\mathcal{C}(X))\neq 0\} = \{[X]\in \mathcal{P}^G_\mu~: ~H(\mathcal{C}(X)) =\Span\{\eta_X\}~ \} $$ By \cref{lem:subfunc_Hom_from_1_functoriality} below, the subset $D_{H}$ is downward-closed, that is, $D\in \Down(\mathcal{P}^G_\mu)$.
        \item For any $D\in \Down(\mathcal{P}^G_\mu)$ we define an additive subfunctor $H_D\subseteq \Hom(\triv, -)$. This functor sends $A\in \Perm(G, \mu)$ to $$H_D(A)=\Span\{\varphi\circ \eta_X~|~ \varphi:\mathcal{C}(X) \rightarrow A,\; [X]\in D  \}$$
        In other words, we consider all morphisms into $A$ from objects $\mathcal{C}(X)$ where $[X]\in D$. Composing such a morphism with $\eta_X$ we obtain an element of $\Hom(\triv, A)$, and the subspace spanned by these elements is denoted by $H_D(A)$. This is clearly an additive subfunctor of $\Hom(\triv, -)$.
        \item For completeness of exposition, let us check explicitly that for any $D\in \Down(\mathcal{P}^G_\mu)$ and any additive subfunctor $H \subseteq \Hom(\triv, -)$ we have: $D_{H_D}=D$ and $H_{D_H}{=} H$. 
        
        Let us prove the first assertion. Clearly, $D\subseteq D_{H_D}$. Given $[Y]\in D_{H_D}$ we have: $H_D(\mathcal{C}(Y))= \Span\{\eta_Y\}$. This means that there exists a morphism $\varphi:\mathcal{C}(X)\to \mathcal{C}(Y)$ in $\Perm(G, \mu)$ such that $[X]\in D$, and because $\Hom(\one,\mathcal{C}(Y))$ is one-dimensional, at least one composite $\varphi \circ \eta_X$ is non-zero. Rescaling it produces $\eta_Y=\varphi\circ \eta_X $. By \cref{lem:subfunc_Hom_from_1_functoriality} below, we conclude that $[X]\succeq_\mu [Y]$ so $[Y]\in D$. Thus $D= D_{H_D}$.

        We now prove that $H_{D_H}{=} H$. For any $X\in \mathbf{S}(G)^{tran}$, we have:
        $$H_{D_H}(\mathcal{C}(X))= \Span\{\eta_X\} \;\;\Longleftrightarrow \;\; [X]\in D_H \;\;\Longleftrightarrow \;\; H(\mathcal{C}(X))= \Span\{\eta_X\} $$
        so $ H_{D_H} {=} H$.

        Under rigidity, a morphism $\mathcal{C}(X)\to \mathcal{C}(Y)$ is identified with a morphism $\one\to \mathcal{C}(Y\times X)$. The decomposition of $Y\times X$ into $G$-orbits identifies the latter space with $\bigoplus_O \kk\eta_O$. Hence the ideal corresponding to $D$ consists of the linear combinations of $\one_O$ for which $[O]\in D$.

    \end{enumerate}
\end{proof}
\begin{lemma}\label{lem:subfunc_Hom_from_1_gen_by_eta}
    Let $X\in \mathbf{S}(G)^{tran}$. Then $\eta_X$ spans $\Hom_{\Perm(G, \mu)}(\triv, \mathcal{C}(X))$.
\end{lemma}
\begin{proof}
   The space $\Hom_{\Perm(G, \mu)}(\triv, \mathcal{C}(X))$ is given by the $G$-invariant Schwartz functions on $\bullet_G\times X = X$. Such functions are all multiples of the constant function, implying the required statement.
\end{proof}
\begin{lemma}\label{lem:subfunc_Hom_from_1_functoriality}
Let $X, Y\in \mathbf{S}(G)^{tran}$.
The following conditions are equivalent:
\begin{enumerate}
    \item $X \succeq_{\mu} Y$,
    \item There exists a morphism $\varphi:\mathcal{C}(X)\to \mathcal{C}(Y)$ in $\Perm(G, \mu)$ such that $\varphi\circ \eta_X=\eta_Y$.
\end{enumerate}

\end{lemma}
\begin{proof}
Denote the $G$-orbits in $Y\times X$ by $Z_1, \ldots, Z_k$. Let $p_i^X:Z_i\to X$, $p_i^Y:Z_i\to Y$ be the corresponding projections. By definition, $X\succeq_\mu Y$ if and only if there exists $i$ such that $\mu(p_i^Y)\neq 0$.

Recall that $\Hom(\mathcal{C}(X), \mathcal{C}(Y))$ has a basis $\{\one_{Z_i}\}_{i=1}^k$ where $\one_{Z_i}:\mathcal{C}(X)\to \mathcal{C}(Y)$ is given as follows: for $f\in \mathcal{C}(X)$, define $ \one_{Z_i}(f)(y)=\int_{x\in X}1_{Z_i}(y,x)f(x)dx$, with $1_{Z_i}(y,x)=1$ if $(y,x)\in Z_i$ and $0$ otherwise. 

Hence the constant function $1_X\in \cC(X)$ is sent by the morphism $\one_{Z_i}$ to the constant function $\one_{Z_i}(1_X)= \mu(p_i^Y)1_Y$. This means that $\one_{Z_i}\circ \eta_X$ sends the constant function $1$ on $\bullet_G$ to $\mu(p_i^Y)\cdot1_Y$. Since $\Hom(\triv, \mathcal{C}(Y))=\Span(\eta_Y)$, we conclude that $\one_{Z_i}\circ \eta_X$ is a multiple of $\eta_Y$, and it is a non-zero multiple if and only if $\mu(p_i^Y)\neq 0$.
This completes the proof of the lemma.

\end{proof}
The following lemma gives a slightly different description the bijection between the posets $\mathcal{TID}(\Perm(G, \mu))$ and $\Down(\mathcal{P}^G_\mu)$.
\begin{lemma}\label{cor:identities_gen_tens_ideal}
    Let $D\in  \Down(\mathcal{P}^G_\mu)$. The tensor ideal in $\Perm(G, \mu)$ generated by $\{\eta_X ~:~{[X]}\in D\}$ coincides with the tensor ideal generated by $\{\id_{\mathcal{C}(X)} ~:~{[X]}\in D\}$. 
\end{lemma}
\begin{proof}
    Let $\mathcal{I}$ be a tensor ideal in $\Perm(G, \mu)$. Clearly, if $\id_{\mathcal{C}(X)}\in \mathcal{I}$ then $\eta_X\in \mathcal{I}$. 

    Now assume that $\eta_X\in \mathcal{I}$. To show that $\id_{\mathcal{C}(X)}\in \mathcal{I}$, it is enough to show that $coev_{\mathcal{C}(X)}\in {\mathcal{I}}$.

    Now $\Delta_{\mathcal{C}(X)}\circ \eta_{X}:\triv \to \mathcal{C}(X\times X)$ corresponds to $coev_{\mathcal{C}(X)}: \triv \to \mathcal{C}(X)\otimes \mathcal{C}(X)^*$ under the isomorphism $\mathcal{C}(X\times X) \cong \mathcal{C}(X)\otimes \mathcal{C}(X)^*$. 
    Clearly $\Delta_{\mathcal{C}(X)}\circ \eta_{X}:\triv \to \mathcal{C}(X\times X)$ is in $\mathcal{I}$, and we are done.

\end{proof}
The following is an immediate corollary of \cref{thm:classification}:
\begin{corollary}\label{cor:principal-class}
    For $X\in \mathbf{S}(G)^{tran}$, the principal tensor ideal
generated by $\id_{\mathcal{C}(X)}$ coincides with the principal tensor ideal
generated by $\eta_X$. This corresponds to the downward-closed subset 
\[
  \downarrow[X]
  =\{[Y]\in \mathcal{P}_{\mu}^{G}:[Y]\preceq_\mu[X]\} \in \Down(\mathcal{P}^G_\mu) .
\]
Thus principal tensor ideals in $\Perm(G, \mu)$ generated by the identity of a single transitive permutation object are in bijection with $\mathcal{P}^G_\mu$.
\end{corollary}

\begin{corollary}\label{cor:tensor_ideals_are_thick} The decategorification map \[ Ob: \TID(\Perm(G, \mu)) \to \ID(\Perm(G, \mu)), \;\; \mathcal{I}\mapsto \{X\in Ob(\Perm(G, \mu)), \; \id_X\in \mathcal{I}\}\] is a bijection.
\end{corollary}

\begin{proof}
To show that this is a bijection, we need to show that every tensor ideal $\mathcal{I}$ in $\Perm(G, \mu)$ is generated by $\{\id_A~:~ \id_A\in\mathcal{I}\}$. This follows directly from \cref{cor:identities_gen_tens_ideal}.
\end{proof}

\begin{remark} The parametrization is reminiscent of the following situation. Let $(W,S)$ be a Coxeter system and $\mathcal{H}$ its associated Hecke category over a field of characteristic zero. In $\mathcal{H}$, the Grothendieck ring is the Hecke algebra, equipped with the canonical Kazhdan-Lusztig  basis $\{C_w\}_{w \in W}$. In $\Perm(G, \mu)$, the Grothendieck semiring is generated by the transitive $G$-sets $\{\mathcal{C}(X)\}$.
 \begin{itemize}
        \item In $\mathcal{H}$, the pre-order $x \le_{LR} y$ is defined by requiring that $C_x$ appears as a direct summand (with a shift) in $h \cdot C_y \cdot h'$ for some objects $h, h'$. The thick two-sided tensor ideals of the Karoubian category $\mathcal{H}$ are in bijection with the downward-closed sets of the two-sided cells (the equivalence classes of $\le_{LR}$).
\item In $\Perm(G, \mu)$, the pre-order $\mathcal{C}(X) \preceq_\mu \mathcal{C}(Y)$ is defined by requiring that $\mathcal{C}(X)$ is a direct summand of $\mathcal{C}(Y) \otimes V$ for some object $V$ (equivalently, $\id_{\mathcal{C}(X)}$ factors through $\mathcal{C}(Y) \otimes V$). The tensor ideals of morphisms in the skeletal category $\Perm(G, \mu)$ are in bijection with the downward-closed sets of the $\mu$-equivalence classes. By \cref{cor:tensor_ideals_are_thick} this coincides with the thick tensor ideals.
\end{itemize}

Note that in $\mathcal{H}$ we have to distinguish between left, right and two-sided ideals since $\mathcal{H}$ is not braided.
\end{remark}

%%%%%%%%%%%%%%

%%%%%%%%%%%%%%%%

\section{Structural properties of the lattice of tensor ideals}
\subsection{Minimal ideals and regular measures}
Let $G$ be an oligomorphic group and $\mu$ a measure on $G$.

\begin{remark}
    For a quasi-regular measure $\mu$ that satisfies (Nil), Harman and Snowden constructed in \cite[§3.4]{HS-classical} the monoidal abelian envelope $\uRep(G, \mu)$ of $\Perm(G, \mu)$ (cf. \cite{HS-oligomorphic}).

    The objects $\mathcal{C}(G/U)$, where $\mu\rvert_U$ is regular, are projective (and injective) objects in $\uRep(G,\mu)$. In particular, $\uRep(G,\mu)$ has enough projective objects (this is one of the fundamental ideas behind the construction of $\uRep(G, \mu)$ in \cite{HS-oligomorphic}). The thick tensor ideal of projective-injective objects in $\uRep(G,\mu)$ is generated, under taking direct sums and direct summands, precisely by the objects $\mathcal{C}(G/U)$ when $\mu\rvert_U$ is regular. 
\end{remark}

\begin{proposition}\label{prop:quasiregular}
Let $(G, \mu)$ be an oligomorphic group and let $U\in \Sub(G)$. The following are equivalent:
 
 \begin{enumerate}
    \item $\mu\rvert_U$ is regular,
        \item There exists $V\subseteq U$ such that $\mu\rvert_V$ is regular and $\mu(U/V)\neq 0$, 
     \item $\mu$ is quasi-regular, and $[G/U]$ is the least element of $\mathcal{P}_{\mu}^G$.
 
 \end{enumerate}
If these conditions hold, the poset ${\Down(\mathcal{P}_{\mu}^G)}\setminus \{\emptyset\}$ has a least element, namely the downward-closed singleton $\{[G/U]\}\subseteq \mathcal{P}_{\mu}^G$. 

\end{proposition}

\begin{proof}[Proof of \cref{prop:quasiregular}]

\mbox{}

\begin{enumerate}[itemindent=0pt,leftmargin=*]
    \item[(1) $\implies$ (2).] This is obvious (set $V=U$).
       \item[(2) $\implies$ (3).] 
    Assume that $\mu\rvert_V$ is regular.
Let $X\in \Set{G}^{tran}$ and choose $x\in X$. Consider the orbit
\[
  G.(V,x)\subseteq G/V\times X.
\]
The fiber of $G.(V,x)\to G/V$ over $V$ is the $V$-orbit $V.x$, which is a transitive $V$-set.  Since $\mu|_{V}$ is regular, its measure is nonzero.  Hence
\[
  G/V\preceq_\mu X
\] and $[G/V]$ is the least element in the poset $\mathcal{P}_{\mu}^G$. Now, $\mu(U/V)\neq 0$ implying that $[G/U]=[G/V]$, so $[G/U]$ is the least element of $\mathcal{P}_{\mu}^G$.

\item[(3) $\implies$ (1).] 
Assume that $[G/U]$ is the least element of $\mathcal{P}_{\mu}^G$ and $\mu$ is quasi-regular. Let $W'\in \Sub(G)$ such that $\mu\rvert_{W'}$ is regular. By our assumption, $[G/U] \preceq_{\mu} [G/W']$ so by \cref{rmk:order_on_subgroups}, there exists an open $W\leq U$ such that $\mu(U/W)\neq 0$ and a $G$-conjugate of $W$ lies in $W'$. This immediately implies that $\mu\rvert_W$ is regular as well.

To prove that $\mu\rvert_U$ is regular, let $V\leq U$ be an open subgroup; we need to show that $\mu(U/V)\neq 0$. Since $W$ is regular, we have: $\mu\bigl(W/(V\cap W)\bigr)\neq 0$. Hence
$$\mu(U/V)\,\mu\bigl(V/(V\cap W)\bigr)=\mu\bigl(U/(V\cap W)\bigr)=\mu(U/W)\,\mu\bigl(W/(V\cap W)\bigr)\neq0 $$ 
proving that $\mu(U/V)\neq 0$.
\end{enumerate}
\end{proof}
\begin{remark}
    The existence of a least element in $\mathcal{P}_{\mu}^G$ does not guarantee that $\mu$ is quasi-regular (hence the two conditions in (3)). Indeed, in the second Delannoy category $\Perm(\Aut(\mathbb{R}, <), \mu_{2})$ (notation of \cite{CS-delannoy}) there is just one non-trivial tensor ideal, by \cite[Proposition 5.14]{CS-delannoy} or by \cref{thm:classification}. Yet the measure $\mu_2$ is not quasi-regular, see \cite[§1]{CS-delannoy}. 
\end{remark}
\begin{corollary}\label{cor:q_reg_implies_minimal_ideal}
    Let $(G, \mu)$ be an oligomorphic group with a quasi-regular measure. Then $\Perm(G, \mu)$ has a unique smallest nonzero ideal, which is a principal ideal.

    If the measure $\mu$ is regular, then this ideal is the trivial (maximal) ideal.
\end{corollary}
\begin{remark}
     In particular, $[G/U]=[G/U']$ for any open subgroups $U, U'\subseteq G$ such that $\mu|_U$, $\mu|_{U'}$ {are} regular. Assume that the measure $\mu$ satisfies property (P) of \cite[Definition 7.10]{HS-oligomorphic}; then for any such $U$, $\id_{\mathcal{C}(G/U)}$ generates the splitting ideal in $\Perm(G, \mu)^{\Kar}$ in the sense of \cite[Definition 2.2]{BEO}.

\end{remark}

\subsection{Useful criterion for a unique non-trivial tensor ideal}
\begin{corollary}\label{cor:two-class}
Assume that $\mu$ is quasi-regular but not regular, and that for every $X\in \Set{G}^{tran}$ at least one of the following holds:
\begin{enumerate}[label=\textup{}]
\item (A) 
$ \mu(X)\neq0$;
\item (B) for some $x\in X$ there exists an open subgroup ${W_*}\leq \Stab(x)$ such that
\[
 \mu(\Stab(x)/W_*)\neq0
\]
and the restricted measure $\mu|_{W_*}$ is regular.
\end{enumerate}
Then $\mathcal{P}^G_\mu$ is the two-element chain
\[
 \mathsf N\lneq[\bullet_G],
\]
where $[\bullet_G]$ consists of the transitive sets of nonzero global measure and $\mathsf N = [X] $ for any $X$ satisfying condition (B). 

In particular, in this setting $\Perm(G, 
\mu)$ has a unique nonzero proper tensor ideal corresponding to the downward-closed subset $\{\mathsf N\}$. This is the ideal of negligible morphisms $\N$.
\end{corollary}

\begin{proof} Let $X$ be a transitive $G$-set and let $V=\Stab(x)$ for some $x\in X$.

If $X$ satisfies condition (A) then by \cref{ex:largest-class}, we have $[X]=  [\bullet_G]$, which is the greatest element in $ \mathcal{P}^G_\mu$.

If $X$ satisfies condition (B) then by \cref{prop:quasiregular}, $[X]$ is the least element in $ \mathcal{P}^G_\mu$. Since $\mu$ is quasi-regular, there exists $V\in \Sub(G)$ such that $\mu|_{V}$ is regular; for such $V$, $G/{V}$ satisfies condition (B), so $ \mathcal{P}^G_\mu$ has a least element, which we denote by $\mathsf N$.

It remains to check that $ \mathsf N\neq[\bullet_G]$. Indeed, by our assumption, $\mu$ is not regular, so there exists $V\leq G$ so that $\mu(G/V){=} 0$. By the above argument, $V$ cannot satisfy condition (A) and so it must satisfy condition (B). Thus $\mathsf N=[G/V]\neq[\bullet_G]$, as required.
\end{proof}

\begin{question} In case $\mathrm{char}~\kk=0$, does the quasi-regularity of the measure $\mu$ imply that $\mathcal{P}_{\mu}^G$ is finite? In positive characteristic, computations with GPT seem to show that the circular Delannoy category has quasi-regular measures where $\mathcal{P}_{\mu}^G$ is infinite.
\end{question}

%%%%%%%%%%%%%%%%%

\subsection{The relative setting}\label{ssec:rel_setting_ideals}
Let $G$ be a (pro-) oligomorphic group and let $\mathcal{E}$ be a stabilizer class (see \cite[Section 2.6]{HS-oligomorphic}). Let $\mu$ be a relative $\mathcal{E}$-measure on $G$. 

We may define the relative posets $\mathcal{P}_{\cE, \mu}^G$ and $\Down(\mathcal{P}_{\cE, \mu}^G)$. The relative versions of all the results in this section hold, including \cref{thm:classification} and \cref{cor:two-class}. 

If $\mathcal{E}$ is a stabilizer class and $\mu$ is an absolute (not just $\mathcal{E}$-relative) measure on $G$, we have an obvious map between posets of tensor ideals $\mathcal{TID}(\Perm(G,  \mu)) \to   \mathcal{TID}(\Perm(G, \mathcal{E}, \mu))$ and an obvious embedding of posets $\mathcal{P}_{\cE, \mu}^G \to  \mathcal{P}_{\mu}^G$ inducing a map $\Down(\mathcal{P}_{ \mu}^G) \to  \Down(\mathcal{P}_{\cE,\mu}^G)$.
The following statement is straightforward, but extremely useful:
\begin{lemma}\label{lem:saturating_class_poset_equiv}
  If $\mathcal{E}$ is a saturating stabilizer class then the map $\mathcal{P}_{\cE, \mu}^G \to  \mathcal{P}_{\mu}^G$ is a bijection.  

  In particular, if $\mathrm{char}~\kk=0$ and $\mathcal{E}$ is a large stabilizer class then the map $\mathcal{P}_{\cE, \mu}^G \to  \mathcal{P}_{\mu}^G$ is a bijection.
\end{lemma}
 In particular, $\mathcal{TID}(\Perm(G,  \mu)) \to   \mathcal{TID}(\Perm(G, \mathcal{E}, \mu))$ is a bijection, which is compliant with \cref{prop:saturating_stab_class}.

%One may view \cref{cor:two-class} as a special case of this statement, with the stabilizer class $\mathcal{E}$ taken to be the collection $$\mathcal{E}=\{U \in \Sub(G)~|~ \mu(G/U)\neq 0 \text{ or } \mu\rvert_U \text{ is regular}\}.$$

%%%%%%%%%%%%%%%%%%%%%%%%%%%%%

%%%%%%%%%%%%%%%%%%%%%%%%%%%%%%

\section{First examples}

%%%%%%%%%%%%%%%%%%%%%%%%%

\subsection{Wreath products with the infinite symmetric group}\label{subsec:wreath-products-iedals}

Let $\kk$ be a field of
characteristic zero, $\Gamma$ be a finite group and $G=\Gamma^{\Omega}\rtimes S_\infty$ be the wreath product of \cref{sec:wreath-products_intro} with $d=|\Gamma|$ and $\Omega=\{1,2,\ldots\}$. Let $D=\Gamma\times \Omega$ and for a finite subset $A\subset \Omega$, denote by $U_A$ the pointwise stabilizer
of $\Gamma\times A\subset D$, and set $U_{m}=U_{\{1, \ldots, m\}}$. The measures $\mu_t$ on $G$ are parameterized by $t\in \kk$ (see \cref{sec:wreath-products_intro}). Recall that by \cref{lem:wreath-regular}, $\mu_t$ is regular when $t\notin \Z_{\geq 0}$; in that case, $\Perm(G, \mu_t)$ has only trivial tensor ideals.

\begin{theorem}\label{thm:wp-tensor-ideals}
Let $m\in\mathbb Z_{\geq0}$. Then
\[
  \mathcal P^G_{\mu_m}
  =\{\mathsf N\prec[\bullet_G]\},
  \qquad \mathsf N=[G/U_{m+1}].
\]

Therefore $\Perm(G,\mu_m)$ has one non-trivial
tensor ideal, namely $\mathcal N$, generated by
  $\mathrm{id}_{\cC(G/U_{m+1})}$.
\end{theorem}

\begin{proof}
We verify the hypotheses of Corollary~\ref{cor:two-class}. Let $X=G/V \in \Set{G}^{tran}$. By
Lemma~\ref{lem:wp-open-subgroups}, we can choose $U_A\leq V$ with finite
index. We order $A=\{a_1,\ldots,a_n\}$. The map
\[
  gU_A\longmapsto\bigl(g(1,a_1),\ldots,g(1,a_n)\bigr)
\]
identifies $G/U_A$ with the ordered $n$-tuples of points of $D$
in pairwise distinct blocks (the action is transitive, and the
base tuple has stabilizer $U_A$).
Forgetting the $\Gamma$-coordinates maps this set to the ordered
configurations of $n$ distinct indices in $\Omega$, with fibres
of cardinality $d^n$. After $j$ indices have been chosen, they
are fixed by their stabilizer and the remaining index set has
measure $t-j$ (by restriction invariance and additivity of the measure). The successive forgetful maps therefore give the falling
factorial $(t)_n=t(t-1)\cdots(t-n+1)$, with $(t)_0=1$.
The Tower Formula \eqref{eq:tower-formula} now implies the following formulas for the measure:
\[
  \mu_t(G/U_A)=d^n(t)_n,
  \qquad
  \mu_t(G/V)=\frac{d^n(t)_n}{[V:U_A]}.
\]
If $n\leq m$, this value at $t=m$ is nonzero, so $X$ satisfies
condition~\textup{(A)} of Corollary~\ref{cor:two-class}.
If $n>m$, the restriction $\mu_m|_{U_A}$ is regular by Equation
\eqref{eq:wp-restricted-regularity} of \cref{lem:wreath-regular}, and
$\mu_m(V/U_A)=[V:U_A]\ne0$. Thus condition~\textup{(B)} of \cref{cor:two-class}
holds with $x=V\in G/V$ and $W_*=U_A$.

Since $\mu_m$ is quasi-regular but not regular,
Corollary~\ref{cor:two-class} gives the two-element poset and the
unique nonzero proper tensor ideal $\mathcal N$.
Proposition~\ref{prop:quasiregular} identifies its least class
with $[G/U_{m+1}]$, since $\mu_m|_{U_{m+1}}$ is regular.
The statement about the generator of $\mathcal{N}$ follows from Corollary~\ref{cor:principal-class}.
\end{proof}

\begin{example}\label{ex:ideals_symmetric}
    For $\Gamma=\{1\}$, we recover the classification of tensor ideals in $\Perm(S_{\infty}, \mu_t)$, which is due to Comes and Ostrik (see \cite[Theorem 3.15]{CO-abelian}):
    \begin{itemize}
        \item when $t\notin \Z_{\geq 0}$ the measure $\mu_t$ is regular and $\Perm(G, \mu_t)$ has only trivial tensor ideals,
        \item when $t\in \Z_{\geq 0}$ the measure $\mu_t$ is quasi-regular but not regular, and $\Perm(G, \mu_t)$ has exactly one non-trivial tensor ideal, the negligible ideal $\mathcal{N}$.
    \end{itemize}
\end{example}
%%%%%%%%%%%%%%%%%%%%

\subsection{Finite groups in positive characteristic}\label{sec:finite-groups}
Let $p$ be a prime integer and $G$ a finite group.

Let $\kk$ be a field of characteristic $p$, and let $   \Rep_{\kk}(G)$ denote the category of finite-dimensional representations of $G$ (that is, finite-dimensional $\kk G$-modules). In this setting, there exists a unique $\kk$-valued measure on $G$-sets, given by $\mu(X) = |X|\cdot 1_{\kk}\in \mathbb{F}_p\subseteq \kk$. The full subcategory $\Perm(G, \mu)^{\Kar} \subseteq \Rep_{\kk}(G)$ then contains precisely the so-called $p$-permutation modules.

For $H\leq G$, write $X_H=G/H$, and choose a Sylow $p$-subgroup
$P_H\in\Syl_p(H)$. Write $Q\leq_G R$ when $Q$ is contained in a
$G$-conjugate of $R$.

\begin{theorem}\label{thm:E-p-subgroups}
Let $\kk$ be a field of characteristic $p$.
Let $G$ be a finite group, and let $\mu$ be the unique $\kk$-valued measure on $G$-sets. One has
\[
  G/H\preceq_\mu G/K
  \quad\Longleftrightarrow\quad
  P_H\leq_G P_K.
\]
Consequently, we have an isomorphism of posets
\[
\quad
  \cP_\mu^G\cong \Sub_p(G)/G,
\]
where $\Sub_p(G)/G$ is the poset of $G$-conjugacy classes of $p$-subgroups,
ordered by subconjugacy.
\end{theorem}

\begin{proof}
    Recall that $G/H\preceq_\mu G/K$ if and only if $H\preceq_\mu K$; equivalently, if and only if there exists $g\in G$ such that $\mu(H/gKg^{-1}\cap H)\neq 0$. In our setting, this is equivalent to the condition that $p$ does not divide the index $[H:gKg^{-1}\cap H]$. The latter is equivalent to $P_H \leq_G P_{gKg^{-1}\cap H}$ which implies that $ P_{gKg^{-1}\cap H} \leq_G P_K$. Thus $G/H\preceq_\mu G/K$ implies $P_H\leq_G P_K$.

    Vice versa, if $P_H\subseteq g P_Kg^{-1}$ for some $g\in G$, then $[H:gKg^{-1}\cap H]$ divides $[H:P_H]$. The latter is not divisible by $p$, so neither is $[H:gKg^{-1}\cap H]$. Hence $G/H\preceq_\mu G/K$. 
\end{proof}
This reproves \cite[Corollary 15.3(ii)]{TW} on the lattice of subfunctors of the fixed point Mackey functor (see \cite{Yo} for the link between the fixed point Mackey functor and the functor $\Hom_G(\triv, -)$ on permutation representations).

Let us describe explicitly the thick tensor ideal $I_D$ in $\Rep_{\kk}(G)$ corresponding to any downward-closed subset $D \subseteq \Sub_p(G)/G$. 
The corresponding thick tensor ideal in $\Perm(G, \mu)$ would just be the intersection of $I_D$ with the full subcategory $\Perm(G, \mu) \subseteq \Rep_{\kk}(G)$.

We start with a well-known definition.
\begin{definition}
    Let $M\in  \Rep_{\kk}(G)$, and let $P\leq G$ be a $p$-subgroup. We say that $M$ is \emph{relatively $P$-projective} if $M$ is a direct summand of $\Ind_P^G N$ for some $N\in \Rep_\kk(P)$.

    We denote by $\mathcal{I}_P$ the full subcategory of $\Rep_{\kk}(G)$ consisting of relatively $P$-projective modules.
\end{definition}
An equivalent characterization of relative projectivity is (see \cite[Proposition 2.3.2]{Benson}): $M$ is relatively $P$-projective if the natural map of $G$-modules $M\to \Ind_P^G \Res_P^G M\cong  M\otimes \kk[G/P]$ splits.
\begin{example}\label{ex:rel_proj}
\mbox{}
\begin{itemize}[itemindent=0pt,leftmargin=*]
\item The relatively $\{1_G\}$-projective representations are precisely the projective $\kk G$-modules.
\item Let $P\leq G$ be a Sylow subgroup. Any $M\in \Rep_{\kk}(G)$ is relatively $P$-projective. 
    \item  Let $U\leq G$ be a subgroup. Then $\kk[G/U]$ is relatively $P$-projective for any $p$-subgroup $P\leq G$ such that $P_U\leq_G P$. Moreover, there exists an indecomposable direct summand $M$ of $\kk[G/U]$ such that $M$ is relatively $P$-projective if and only if $P_U\leq_G P$ (see \cite[Chapter 4, \S 8]{NagaoTsushima}).
\end{itemize}

\end{example}
 The following elementary statements are also well-known.
\begin{lemma}
    Let $P',P\leq G$ be two $p$-subgroups of $G$.
    \begin{enumerate}
  
        \item If $P \leq P'$ then $\mathcal{I}_P \subseteq \mathcal{I}_{P'}$.
                \item If $P, P'$ are conjugate then $\mathcal{I}_P = \mathcal{I}_{P'} $.
        \item $\mathcal{I}_P$ is a Karoubi additive subcategory which is closed under tensoring with finite-dimensional $\kk G$-modules; that is, $\mathcal{I}_P$ defines a thick tensor ideal.

    \end{enumerate}
\end{lemma}
\begin{proof}
\mbox{}

    \begin{enumerate}[itemindent=0pt,leftmargin=*]
    
        \item If $P \leq P'$ and $M$ is relatively $P$-projective, then $M$ is a direct summand of $\Ind^G_P N$ for some $N\in \Rep_{\kk}(P)$. This means that  $M$ is a direct summand of $\Ind^G_{P'} \Ind^{P'}_P N$, so $M$ is relatively $P'$-projective.
                \item If $P, P'$ are conjugate then $\kk[G/P]\cong \kk[G/P']$ as $G$-modules. Therefore, the natural map $M\to   M\otimes \kk[G/P]$ splits if and only if the map $M\to   M\otimes \kk[G/P']$ splits.
        \item The property of being Karoubi additive is clear, and the tensoring property is proved in \cite[Corollary 2.3.5]{Benson}.
 
    \end{enumerate}
\end{proof}
Now let $D \subseteq \Sub_p(G)/G$ be a downward-closed subset; that is, we may identify $D$ with a subset of $\Sub_p(G)$ closed under conjugation and inclusion.
\begin{proposition}
    Let $I_D$ be the thick tensor ideal in $\Rep_{\kk}(G)$ given by the full subcategory $\left(\bigcup_{P\in D} \mathcal{I}_P\right)^{\oplus} $. Then $I_D\cap \Perm(G, \mu)$ is precisely the thick tensor ideal in $\Perm(G, \mu)$ corresponding to $D$ under \cref{thm:classification}.
\end{proposition}
\begin{proof}
    By \cref{thm:classification}, the thick tensor ideal in $\Perm(G, \mu)$ corresponding to $D$ is generated as an additive Karoubi category by the permutation representations $\kk[G/U]$ where $U\leq G$ is a subgroup such that $P_U\in D$. 
    By \cref{ex:rel_proj}, if $P_U\in D$ then $\kk[G/U] \in I_D\cap \Perm(G, \mu)$. 
    
    It remains to check that whenever $\kk[G/U] \in I_D$ for some $U\leq G$, we have: $P_U\in D$. Indeed, assume $\kk[G/U]\in I_D$ and let $M$ be the direct summand of $\kk[G/U]$ mentioned in \cref{ex:rel_proj}. Then $M \in I_D$. Now, $M$ is indecomposable, so $M \in \mathcal{I}_P$ for some $ P_U\leq_G P\in D$. Since $D$ is downward closed we conclude that $P_U\in D$.
\end{proof}

%%%%%%%%%%%%%%%%%%%%%%%%%%%%

%%%%%%%%%%%%%%%%%%%%%%%%%

\subsection{Products of oligomorphic groups}

\begin{corollary}\label{cor:product-poset}
    Let $\mathrm{char}~\kk=0$. Let $G_1,\ldots,G_d$ be pro-oligomorphic groups with $\kk$-valued measures $\mu_i$. Put
$G=\prod_{i=1}^dG_i$, $\mu=\boxtimes_{i=1}^d\mu_i$, and
$\mathcal D_i=\Perm(G_i,\mu_i)^{\Kar}$. Then
\[
 \mathcal P^G_\mu\cong\prod_{i=1}^d\mathcal P^{G_i}_{\mu_i},\qquad
 \TID\!\left(\boxtimes_{i=1}^d\mathcal D_i\right)
 \cong\TID\!\left(\Perm(G,\mu)^{\Kar}\right)
 \cong\Down\!\left(\prod_{i=1}^d\mathcal P^{G_i}_{\mu_i}\right).
\]
\end{corollary}

\begin{proof}

We only need to prove that $\mathcal P^G_\mu\cong\prod_{i=1}^d\mathcal P^{G_i}_{\mu_i}$, as the second statement follows from \cref{thm:classification}. By \cref{lem:saturating_class_poset_equiv} and \cref{cor:prod_is_large_stab_class},
we have an isomorphism of posets $\mathcal{P}^G_{\mathcal{E}_{\times},\mu}\cong \mathcal P^G_\mu$. It remains to check that $\mathcal P^G_{\mathcal{E}_{\times},\mu}\cong\prod_{i=1}^d\mathcal P^{G_i}_{\mu_i}$, which follows directly from the definition of product measure.

\end{proof}

\begin{example}\label{ex:boolean-product-ideals}
Take $\mathrm{char}~\kk=0$, $G_i=S_\infty$ and $t_i\in\mathbb Z_{\geq0}$ for every $i$. \cref{ex:ideals_symmetric} gives
$\mathcal P^{G_i}_{\mu_{t_i}}=\{\mathsf N_i<[\bullet_{G_i}]\}$. By Corollary \ref{cor:product-poset}
\[
 \mathcal P^{S_\infty^{\times d}}_{\boxtimes_i\mu_{t_i}}
 \cong B_d,\qquad
 \TID\!\left(\boxtimes_i\Perm(S_\infty,\mu_{t_i})^{\Kar}\right)
 \cong\Down(B_d),
\]
where $B_d=(2^{\{1,\ldots,d\}},\subseteq)$ is the Boolean lattice. It is well-known that this contains $M(d)$ downward-closed subsets (which parametrize the tensor ideals including the zero ideal and the category itself)
where $M(d)$ is the $d$-th \emph{Dedekind number}. These numbers grow very rapidly (for $d=1,2,3,4,5$ these numbers are
$3,6,20,168,7581$, respectively).
    
\end{example}

%%%%%%%%%%%%%%%%%%%%%%

\subsection{Equivariantizations}\label{subsec:ideals_equivariantization}
Let $G$ be an oligomorphic group, and $\Gamma$ be a finite group acting on $G$ by continuous group automorphisms. Let $\kk$ be a field with $\operatorname{char} ~\kk \nmid |\Gamma|$, and let $\mu$ be an (absolute) $\kk$-valued measure on $G$ which is $\Gamma$-invariant.

    \begin{lemma}\label{lem:twist-order}
    Let $X, Y\in \Set{G}^{tran}$ and $\gamma\in \Gamma$. If $X\succeq_{\mu} Y$ then $X^{\gamma}\succeq_{\mu} Y^{\gamma}$. 
    \end{lemma}
    \begin{proof}
        Assume $X\succeq_{\mu} Y$. By definition, there exists $Z\in \Set{G}$ with $G$-maps $$f_X:Z\to X, ~f_Y:Z\to Y$$ such that $\mu(f_Y)\neq0$. Consider the triple $(Z^{\gamma}, ~f_X^{\gamma}:Z^{\gamma}\to X^{\gamma}, ~f_Y^{\gamma}:Z^{\gamma}\to Y^{\gamma})$. Since $\mu$ is $\Gamma$-invariant, we have: $\mu(f_Y^{\gamma})=\mu(f_Y)\neq 0$. Thus $X^{\gamma}\succeq_{\mu} Y^{\gamma}$. 
    \end{proof}
    We may therefore now define $[X]^{\gamma}=[X^{\gamma}]$ for any $\gamma\in \Gamma$ and $[X]\in \mathcal{P}_{\mu}^G$.
    \begin{definition}
        A downward-closed subset $D\subseteq \mathcal{P}^G_\mu $ is called \emph{$\Gamma$-stable} if $[X]\in D$ $\implies$ $[X]^\gamma\in D$ for every $\gamma\in \Gamma$. We denote by $\Down(\mathcal{P}^G_\mu)^{\Gamma}$ the poset of downward-closed $\Gamma$-stable subsets of $ \mathcal{P}^G_\mu $.
    \end{definition}

   Let $\overline{\mu}$ be the measure on $G\rtimes \Gamma$ corresponding to $\mu$, as before.

   \begin{lemma}\label{lem:orbitwise-order}
Let $\oX,\overline{Y}\in \Set{\oG}^{tran}$ and let $X\in G\backslash\oX$. Then
\[
 \oX\succeq_{\overline{\mu}}\overline{Y}
 \quad\Longleftrightarrow\quad
 X\succeq_{\mu} Y \text{ for some } Y\in G\backslash\overline{Y}.
\]
   \end{lemma}
  \begin{proof}
Let $X\in G\backslash\oX$ and let $\pi_{\overline{X}}: \oX\times \overline{Y} \to  \overline{X}$, $\pi_{\overline{Y}}: \oX\times \overline{Y} \to  \overline{Y}$ be the projections on $\oX, ~ \overline{Y}$ respectively. For a fixed $\overline{Z}\in \overline{G}\backslash( \oX\times \overline{Y})$, let $\pi$ be the restriction of $\pi_{\overline{Y}}$ to ${\overline{Z}}$. Then for any $Y\in G \backslash\overline{Y}$, the fiber of $\pi$ over a point $y\in Y$ is $$\bigsqcup_{Z\in G\backslash\overline{Z}, ~ \pi(Z)=Y} \pi^{-1}(y)\cap Z = \bigsqcup_{Z\in G\backslash\overline{Z}, ~ \pi(Z)=Y} \pi^{-1}\rvert_Z(y).$$

Hence we obtain: 
\begin{equation}\label{eq:fibre-sum}
\overline{\mu}(\pi) = \sum_{Z \in G\backslash\overline{Z}, ~ \pi(Z)=Y}\mu\bigl(\pi|_{Z}:Z\to Y\bigr)
\end{equation}
  
  \begin{enumerate}[itemindent=0pt,leftmargin=*]
    \item[$(\Longrightarrow)$] Assume $\oX\succeq_{\overline{\mu}}\overline{Y}$. By our assumption and \cref{lem:succeq_mu}, there exists $\overline{Z}\in \overline{G}\backslash( \oX\times \overline{Y})$ such that $\overline{\mu}(\pi_{\overline{Y}}\rvert_{\overline{Z}})\neq 0$. Let $\pi: \overline{Z}\to \overline{Y}$ be the projection as before. Then \eqref{eq:fibre-sum} implies the existence of $Z' \in G\backslash\overline{Z} $ such that $\mu(\pi|_{Z'})\neq 0$. Let $Y'=\pi(Z')$ and let $X'=\pi_{\oX}(Z')$. Then $X'\in G\backslash\oX$, so $X'=X^{\gamma}$ for some $\gamma\in \Gamma$ by \cref{lem:semidirect_restriction_on_sets}. Now let $Y=(Y')^{\gamma^{-1}}$ and $Z=(Z')^{\gamma^{-1}}$. Then $\mu(\pi: Z\to Y)=\mu(Z'\to Y')\neq 0$, and $X=\pi_{\oX}(Z)$. Thus $ X\succeq_{\mu} Y  $.
      \item[$(\Longleftarrow)$] Let $Y\subseteq\overline Y|_G$ be a $G$-orbit with $X\succeq_\mu Y$. By \cref{lem:succeq_mu}, there is a $G$-orbit $Z\subseteq X\times Y$ such that $\mu(Z\to Y)\neq 0$. Let $\overline{Z}=\{\oG.z~:~z\in Z\}\subset \oX\times \overline{Y}$, and $\pi: \overline{Z}\to \overline{Y}$ be the projection as before. Clearly $\overline{Z}$ is transitive, and its $G$-orbits are the sets $\gamma Z=Z^\gamma $, $\gamma\in\Gamma$, with $\pi( \gamma Z)= \gamma Y$. Let $\Gamma_{Y}=\{\gamma\in\Gamma:\gamma Y =Y\}$ be the stabilizer of $Y$ under the action of $\Gamma$ on $G\backslash\overline{Y}$, and let $S$ be the orbit of $Z$ under the action of $\Gamma_{Y}$. Since $\operatorname{char}~\kk\nmid|\Gamma|$, we conclude that $\operatorname{char}~\kk\nmid|S|$. Now \eqref{eq:fibre-sum} gives:
      $$ \overline{\mu}(\pi) =  \sum_{Z'\in S}\mu\bigl(\pi|_{Z'}:Z'\to Y\bigr)=|S|\cdot \mu(Z\to Y)\neq 0.$$
      The second equality follows from the fact that $\gamma Y=Y$ for $\gamma\in \Gamma_Y$ and $\mu$ is $\Gamma$-invariant. So $\oX\succeq_{\overline\mu}\overline Y$ by \cref{lem:succeq_mu}.
      %$$ \overline{\mu}(\pi)=|S|\cdot \mu(Z\to Y)=\frac{|\Gamma|}{|\Gamma^Y|}\cdot \mu(Z\to Y) \neq 0$$

    \end{enumerate}

\end{proof}
   
\begin{proposition}\label{prop:equivariantization-ideals}
    We have bijections between posets 
    $$\mathcal{TID}(\Perm(G, \mu)^{\Gamma})~\cong~ \mathcal{TID}(\Perm(G\rtimes \Gamma, \overline{\mu}))~\cong~ \Down(\mathcal{P}^{G\rtimes \Gamma}_{\overline{\mu}})~\cong ~\Down(\mathcal{P}^G_\mu)^{\Gamma}.$$
    The bijection $\Down(\mathcal{P}^{G\rtimes \Gamma}_{\overline{\mu}}) ~\xrightarrow{\sim}~ \Down(\mathcal{P}^G_\mu)^{\Gamma}$ is given explicitly by 
    \[\overline{D}~\longmapsto~ D=\left\{[X]\in \mathcal{P}^G_\mu~:~\mathrm{Mor}_G(X, \overline{X}\rvert_G)\neq \emptyset  \text{ for some } [\overline{X}]\in \overline{D}, ~ X\in [X]\right\}.\]
\end{proposition}
\begin{proof}
    The first bijection follows from \cref{prop:equivariantization-product}, together with the bijection between tensor ideals in a category and in its Karoubi completion. The second is just the statement of classification of the tensor ideals in $\Perm(G\rtimes \Gamma, \overline{\mu})$, so we only need to prove that the explicit map $\Down(\mathcal{P}^{G\rtimes \Gamma}_{\overline{\mu}}) ~\xrightarrow{}~ \Down(\mathcal{P}^G_\mu)^{\Gamma}$ is well-defined and is a bijection. 

    Given $\overline{D}\in \Down(\mathcal{P}^{G\rtimes \Gamma}_{\overline{\mu}})$, let 
    \[\Phi(\overline{D})=\left\{[X]\in \mathcal{P}^G_\mu~:~\mathrm{Mor}_G(X, \overline{X}\rvert_G)\neq \emptyset  \text{ for some } [\overline{X}]\in \overline{D}, ~ X\in [X]\right\}.\]
{\it The subset $\Phi(\overline{D})\subseteq \mathcal{P}^G_\mu$ is downward-closed:} given $[X]\in \Phi(\overline{D})$, fix $\overline{X}\in \Set{G\rtimes \Gamma}$ such that $\mathrm{Mor}_G(X, \overline{X}\rvert_G)\neq \emptyset$. If $Y\preceq_{\mu} X$ then there exists $Z\in [Y]$ such that $\mathrm{Mor}_G(Z, X)\neq \emptyset$. Thus $\mathrm{Mor}_G(Z,  \overline{X}\rvert_G)\neq \emptyset$ and so $[Y]=[Z] \in \Phi(\overline{D})$.

Moreover, as a downward-closed subset of $\mathcal{P}^G_{\mu}$, $\Phi(\overline{D})$ is clearly generated by the set $$\left\{[X]\in \mathcal{P}^G_\mu~:~ X\in G\backslash (\overline{X}\rvert_G)  \text{ for some } [\overline{X}]\in \overline{D}\right\}.$$

{\it The set $\Phi(\overline{D})$ is $\Gamma$-stable:} due to \cref{lem:semidirect_restriction_on_sets}, for $\overline{X} \in \Set{G\rtimes \Gamma}$ and $X\in \Set{G}$, there exists a $G$-map $X\to \overline{X}\rvert_G$ if and only if some $G$-orbit $O$ in $\overline{X}\rvert_G$ has $\mathrm{Mor}_G(X, O)\neq \emptyset$. By \cref{lem:semidirect_restriction_on_sets}, for every $\gamma\in \Gamma$, $\mathrm{Mor}_G(X^{\gamma},O^{\gamma})\neq \emptyset$. 

Thus the correspondence $\overline{D}\mapsto \Phi(\overline{D})$ is well-defined. Let us check that it is a bijection. We define the inverse map $\Down(\mathcal{P}^G_\mu)^{\Gamma} ~\xrightarrow{}~ \Down(\mathcal{P}^{G\rtimes \Gamma}_{\overline{\mu}})$ to be given by 
$$ D~\in ~\Down(\mathcal{P}^G_\mu)^{\Gamma}\; \longmapsto \; \Psi(D)=\left\{[\overline{X}]~\in~ \mathcal{P}^{G\rtimes \Gamma}_{\overline{\mu}}\,:\; \exists \, X\in G\backslash (\overline{X}\rvert_G) \text{ such that }  [X]\in D   \right\}.$$
Since $D$ is $\Gamma$-stable, and using \cref{lem:semidirect_restriction_on_sets}, ``there exists'' may be replaced by ``for any''.

\emph{$\Psi(D)$ is well defined and downward-closed:} let $\oX$ be a transitive $\oG$-set with a $G$-orbit $X$ such that $[X]\in D$, and let $\overline Y$ be a transitive $\oG$-set with $\overline Y\preceq_{\overline\mu}\oX$. By \cref{lem:orbitwise-order}, $X\succeq_\mu Y'$ for some $G$-orbit $Y'$ of $\overline Y|_G$. As $D$ is downward-closed, $[Y']\in D$. Applied to $\overline Y\sim_{\overline\mu}\oX$, this shows that the defining condition of $\Psi(D)$ only depends on the class $[\oX]$; in general, it shows that $\Psi(D)$ is downward-closed. 

\emph{Showing that $\Psi\circ\Phi=\id$:} Let $\overline D\in\Down(\mathcal{P}^{G\rtimes \Gamma}_{\overline{\mu}})$. If $[\oX]\in\overline D$ and $X$ is a $G$-orbit of $\oX|_G$, then $[X]\in \Phi(\overline D)$, so $[\oX]\in\Psi\Phi(\overline D)$. 

Conversely, let $[\oX]\in\Psi\Phi(\overline D)$, and let $X$ be a $G$-orbit of $\oX|_G$ with $[X]\in\Phi(\overline D)$. Then there exist $[\overline Y]\in\overline D$ and a $G$-orbit $Y'$ of $\overline Y|_G$ with $Y'\succeq_\mu X$. By \cref{lem:orbitwise-order}, applied to the pair $\overline Y,\oX$ and the $G$-orbit $Y'\subseteq\overline{Y}|_G$, we get $\overline Y\succeq_{\overline\mu}\oX$. Since $\overline D$ is downward-closed, $[\oX]\in\overline D$.

\emph{Showing that $\Phi\circ\Psi=\id$.} Let $D\in\Down(\mathcal{P}^G_\mu)^\Gamma$. For any $[\oX]\in\Psi(D)$ and any $X\in G\backslash(\oX|_G)$, we have: $ [X]\in D$. Now, $\Phi(\Psi(D))$ is generated by such classes $[X]$ and $D$ is downward-closed, so $\Phi\Psi(D)\subseteq D$. 

Conversely, let $X\in \Set{G}^{tran}$ such that $[X]\in D$, and let $\oX=\oG\times_G X$ be the induced $\oG$-set considered in the proof of \cref{prop:equivariantization-product}. It is transitive, and $X$ is one of the $G$-orbits of $\oX|_G\cong\coprod_{\gamma\in\Gamma}X^\gamma$. Hence $[\oX]\in\Psi(D)$ and $[X]\in \Phi\Psi(D)$. Thus $D\subseteq \Phi\Psi(D)$.

Therefore $\Phi$ and $\Psi$ are mutually inverse, inclusion-preserving bijections, as required.

\end{proof}

\begin{example}

Let $G_0=S_\infty$, $\mathrm{char} ~\kk=0$ and $t\in\mathbb Z_{\geq0}$. Consider the measure $\mu_t$ on $S_{\infty}$ (see \cref{subsec:main-examples}). 

Set $\nu=\mu_t^{\boxtimes d}$ and let $\overline{\nu}$ be the measure for $S_\infty^{\times d}\rtimes S_d$ corresponding to $\nu$. Then $\Perm(S_\infty^{\times d}, \mu_t^{\boxtimes d})^{S_d} $ has exactly $d+2$ tensor ideals
\[
 0=\mathcal I_{d+1}\subsetneq\mathcal I_d\subsetneq\cdots
 \subsetneq\mathcal I_1=\mathcal N\subsetneq\mathcal I_0=\Perm(S_\infty^{\times d}, \mu_t^{\boxtimes d})^{S_d}.
\]
Indeed, in this case $\mathcal{P}^{G_0^{\times d}}_\nu$ is the Boolean lattice $B_d=(2^{\{1,\ldots, d\}}, \subseteq)$, and the $S_d$-stable subsets are the ones which are invariant under the action of $S_d$ on $2^{\{1,\ldots, d\}}$. Each such $D\in \Down(B_d)^{S_d}\setminus\{\emptyset\}$ corresponds to choosing an integer $0\leq r\leq d$ so that the sets $A\in 2^{\{1,\ldots, d\}}$ belonging to $D$ would have at most $r$ elements.

The same is true for $\Perm(S_\infty^{\times d}\rtimes S_d, \overline{\nu})$, since $\left(\Perm(S_\infty^{\times d}, \mu_t^{\boxtimes d})^{S_d} \right)^{\Kar}\cong \Perm(S_\infty^{\times d}\rtimes S_d, \overline{\nu})^{\Kar}$ by \cref{prop:equivariantization-product}.
\end{example}

%%%%%%%%%%%%%%%%%%%%%%%

%%%%%%%%%%%%%%%%%%%%%%%%%%

\part{Tensor ideals for classical linear groups}\label{sec:classical}

We assume throughout this section that $\mathrm{char}~\kk=0$. We will compute the lattice of tensor ideals for all classical linear (finitary) groups and $\kk$-valued measures $\mu$, focusing on the singular measures. This section is heavily based on \cite{HS-classical}, but we will avoid model-theoretic language to simplify the exposition.

\section{The setup: classical linear oligomorphic groups}\label{ssec:setup_classcial_grps}

Let $q\in \Z_{>0}$ be a power of a prime, 
and let $\mathrm{Vec}^{{\infty}}_{\F_q}$ denote the category of all ({not} necessarily finite dimensional) vector spaces over $\F_q$. Let $\mathcal{A}$ be a category with a faithful functor $\mathcal{A}\to \mathrm{Vec}^{{\infty}}_{\F_q}$. Using this forgetful functor, we may consider objects in $\mathcal{A}$ as $\F_q$-vector spaces, and maps in $\mathcal{A}$ as $\F_q$-linear maps. The objects in $\mathcal{A}$ which are finite-dimensional over $\F_q$ will be referred to in this section as \emph{finite structures}. 

We will denote by $\mathcal{A}^{mono}$ the wide subcategory whose morphisms are only monomorphisms. For $A, B\in \mathcal{A}$, we will denote: $$\Emb(A, B)=\Hom_{\mathcal{A}^{mono}}(A, B).$$ 
The full subcategory of $\mathcal{A}$ (respectively, $\mathcal{A}^{mono}$) on all the finite structures will be denoted by $\mathcal{A}_f$ (respectively, $\mathcal{A}^{mono}_f
$).

In all the cases below, we will fix a category $\mathcal{A}$ as above, and a countable-dimensional $\F_q$-vector space $\V \in \mathcal{A}$.
We will denote by $G=\Aut_{\mathcal{A}}(\V) {\subsetneq} \Emb(\V, \V)$ the corresponding group of automorphisms; the elements of $G$ can be represented by $\F_q$-matrices with countably many rows and columns.

\medskip

The pair $(G, \V)$ will be an oligomorphic pair, with the $\V$-stabilizer class $\mathcal{E}_\V$ (see \cref{ex:Omega_stabilizer_class}) being a large stabilizer class in $G$. 

For a finite structure $A \in \mathcal{A}^{mono}_f$, the group $G$ acts on the set $X_A=\Emb(A, \V)$. In our examples, this will always be a transitive $G$-set. By \cite[Proposition 2.7]{HS-classical} the functor $ A \mapsto X_A$ will be an antiequivalence between $\mathcal{A}^{mono}_f$ and the category of transitive $\V$-smooth $G$-sets (see \cref{sec:relative} for definitions).

Throughout this paper, the term \emph{classical groups} will refer to one of the five oligomorphic pairs $(G, \V)$ listed below: the general linear {(parabolic and Levi)}, the symplectic, the orthogonal and the unitary. The specific category $\mathcal{A}$ in each of these cases will be described below.

\medskip

For each $A\in \mathcal{A}_f$ in the cases below, we will denote by $A_\mathrm{rad}$ the \emph{radical} part of $A$, as defined below. In each case, there exists a non-degenerate complement $A_{\mathrm{nd}}$ of $A_\mathrm{rad}${,} so that $A=A_\mathrm{nd}\oplus A_\mathrm{rad}$ and this decomposition is orthogonal with respect to the relevant form on $A$, if present. We will denote by $\rk(A)$ an appropriate \emph{rank} of $A_{\mathrm{nd}}$, again defined case-by-case, but independent of the choice of $ A_{\mathrm{nd}}$. If $A_{\rad}=0$ then we write $A=A_{\nd}$ and say that $A$ is \emph{non-degenerate}.

\subsection{The classical groups}\label{sec:classical-cases}

\subsubsection{The general linear case - parabolic version}

In this case, we take $\mathcal{A}$ to be $\mathrm{Vec}_{\F_q}$. Up to isomorphism, the finite structures are parametrized by non-negative integers $a$. The finite structure corresponding to $a$ is an $a$-dimensional $\F_q$-vector space. 

For $A\in \mathcal{A}$, we set $A_{\mathrm{nd}}=A$, $A_{\mathrm{rad}}=\{0\}$ and $\rk(A)=\dim_{\F_q} A_{\mathrm{nd}}$.

The group $G$ consists of {all linear} automorphisms of the countable-dimensional space $\V$. The $\V$-stabilizer class then consists of parabolic-type subgroups of $G=GL_{\infty}(\F_q)$ (pointwise stabilizers of a finite-dimensional linear subspace of $\V$).

\subsubsection{The general linear case - Levi version}
This part is based on \cite[\S6.1]{HS-classical}.

In this case, we take $\mathcal{A}$ to be the category of $\Z/2$-graded $\F_q$-vector spaces $A=A_{\bar 0}\oplus A_{\bar 1}$
with a bilinear pairing $A_{\bar 0}\times A_{\bar 1}\to\F_q$. The morphisms are grading-preserving, pairing-preserving linear maps.

For $A\in \mathcal{A}$, we set $(A_{\mathrm{rad}})_{\bar 0},~ (A_{\mathrm{rad}})_{\bar 1}$ to be the radical summands of $A_{\bar 0}, A_{\bar 1}$ respectively, and $A_{\mathrm{nd}}$ to be the complement to $A_{\mathrm{rad}}$ - a maximal subspace of $A$ on which the pairing is perfect. We denote $\rk(A)=\dim_{\F_q}(A_{\mathrm{nd}})_{\bar 0}={\dim_{\F_q}}(A_{\mathrm{nd}})_{\bar 1}$.

Up to isomorphism, the finite structures are parametrized by triples of non-negative integers $(a,b,c)$. The  finite structure corresponding to a triple $(a,b,c)$ is given by $V_{a,b,c}$ where $(V_{a,b,c})_{\mathrm{nd}} $ has dimension $(a\mid a)$ and carries a perfect pairing, while the dimensions of the $0$-th and the $1$-st grade of $(V_{a,b,c})_{\mathrm{rad}}$ are respectively $b$ and $c$.

In the countable dimensional space $\V=\V_{\bar 0}\oplus\V_{\bar 1}$ we fix an explicit basis:
\[
\V_{\bar 0}=\bigoplus_{i\geq1}\F_q e_i,
 \qquad
\V_{\bar 1}=\bigoplus_{i\geq1}\F_q f_i
 \qquad \text{ so that } \qquad
 \langle e_i,f_j\rangle=\delta_{ij}.
\] 
We write: $X_{a,b,c}=\Emb(V_{a,b,c}, \V)$.
The group $G$ consists of grading-preserving automorphisms of $\V$ preserving the pairing. The $\V$-stabilizer class then consists of Levi-type subgroups of $G$.

\begin{example}
\mbox{}
\begin{itemize}[itemindent=0pt,leftmargin=*]
    \item For $A=V_{1,0,0}=\F_q\oplus \F_q$, an element $i\in X_{1,0,0}$ is determined by choosing {vectors} $v_i\in \V_{i}$ for $i=0,1$ so that $ \langle v_0,v_1\rangle=1$. Then $\Stab(i)\subseteq G$ is the subgroup of automorphisms of $\V$ which fix both $v_0, v_1$.  
\item For $A=V_{0,1,0}$, an element $i\in X_A=\Emb(A, \V)$ is determined by choosing a non-zero vector in $\V_{\bar 0}$. Then $\Stab(i)\subseteq G$ is the subgroup of automorphisms of $\V$ which fix that vector. 
\end{itemize}

\end{example}

\subsubsection{The symplectic case}
This part is based on \cite[\S8.1]{HS-classical}. We assume that $q$ is odd.

In this case, we take $\mathcal{A}$ to be the category of $\F_q$-vector spaces equipped with an alternating bilinear form. The morphisms are isometries.

For $A\in \mathcal{A}$, we set $A_{\mathrm{rad}}=Rad(A)$ and $A_{\mathrm{nd}}$ to be a maximal nondegenerate symplectic subspace. We denote $\rk(A)=\frac{1}{2}\dim_{\F_q} A_{\mathrm{nd}}$.

Up to isomorphism, the finite structures are parametrized by pairs of non-negative integers $(a,b)$. The finite structure corresponding to a pair $(a,b)$ is $V_{a,b}$ such that $2a=\dim_{\F_q} (V_{a,b})_{\mathrm{nd}}$ and $b=\dim_{\F_q} (V_{a,b})_{\mathrm{rad}}$. 

In the countable dimensional space $(\V, \langle, \rangle)\in \mathcal{A}$ we fix an explicit basis $\{e_i,f_i\}_{i\geq1}$ such that
\[
 \langle e_i,e_j\rangle=\langle f_i,f_j\rangle=0,
 \qquad \langle e_i,f_j\rangle=\delta_{ij}.
\]
 We write $X_{a,b}=\Emb(V_{a,b},\mathbf V)$.

\subsubsection{The orthogonal case}\label{sec:orthogonal-structure}
This part is based on\cite[\S9.1]{HS-classical}. We assume that $q$ is odd.

In this case, we take $\mathcal{A}$ to be the category of $\F_q$-vector spaces with a symmetric bilinear form (such spaces are called quadratic spaces in \cite{HS-classical}). The morphisms are isometries.

For $A\in \mathcal{A}$, we set $A_{\mathrm{rad}}=Rad(A)$ and $A_{\mathrm{nd}}$ to be a maximal subspace with a nondegenerate induced form. We denote $\rk(A)=\dim_{\F_q} A_{\mathrm{nd}}$.

Up to isomorphism, the finite structures are parametrized by triples $(a,\epsilon, b)$ where $a,b$ are non-negative integers and $\epsilon\in\F_q^\times/(\F_q^\times)^2$. The finite structure corresponding to a triple $(a,\epsilon, b)$ is $V_{a,\epsilon,b}$ such that $a=\dim_{\F_q} (V_{a,\epsilon, b})_{\mathrm{nd}}$, $b=\dim_{\F_q}(V_{a,\epsilon,b})_{\mathrm{rad}}$, and $\epsilon\in\F_q^\times/(\F_q^\times)^2$ is the \emph{modified discriminant}. If $a>0$ and $\delta$ is the ordinary discriminant of the induced form on $(V_{a,\epsilon, b})_{\mathrm{nd}}$, define the modified discriminant to be $\epsilon=(-1)^{\lfloor\frac{a}{2}\rfloor}\delta$ \cite[\S9.1]{HS-classical}. 
 
 For $a=0$ the discriminant label is irrelevant.

In the countable dimensional space $(\V, \langle, \rangle)\in \mathcal{A}$ we fix an explicit basis $\{e_i,f_i\}_{i\geq 1}$ such that 
\[
 \langle e_i,e_j\rangle=\langle f_i,f_j\rangle=0,
 \qquad \langle e_i,f_j\rangle={\delta_{ij}}.
\]
 We write $X_{a,\epsilon,b}=\Emb(V_{a,\epsilon,b},\mathbf V).$

\subsubsection{The unitary case}
This part is based on\cite[\S10.1]{HS-classical}. Let $\F = \F_{q^2}$ be a quadratic extension of $\F_q$ with $q$ odd and let $\tau$ be its nontrivial Galois automorphism.

In this case, we take $\mathcal{A}$ to be the category of $\F $-vector spaces with a Hermitian form (also called a symmetric sesquilinear form)
\[
 \langle ax,by\rangle=\tau(a)b\langle x,y\rangle,
 \qquad \tau(\langle x,y\rangle)=\langle y,x\rangle.
\]
Morphisms are $\F $-linear isometries. 

For $A\in \mathcal{A}$, we set $A_{\mathrm{rad}}=Rad(A)$ and $A_{\mathrm{nd}}$ to be a maximal subspace with non-degenerate induced form. We denote $\rk(A)=\dim_{\F } A_{\mathrm{nd}}$.

Up to isomorphism, the finite structures are parametrized by pairs of non-negative integers $(a,b)$. The finite structure corresponding to a pair $(a,b)$ is $V_{a,b}$ such that $a=\dim_{{\F}} (V_{a,b})_{\mathrm{nd}}$ and $b=\dim_{\F } (V_{a,b})_{\mathrm{rad}}$. 

In the countable dimensional space $(\V, \langle, \rangle)\in \mathcal{A}$ we fix an explicit basis $\{e_i\}_{i\geq1}$ such that $\langle e_i,e_j \rangle = \delta_{i,j}$.
 We write $ X_{a,b}=\Emb(V_{a,b},\mathbf V)$.

\subsection{Measures and counting polynomials}\label{sec:counting}
Let us fix a classical oligomorphic pair $(G, \V)$ as above.

We begin with the non-orthogonal case. In \cite[§§6.7, 8.5, 10.4]{HS-classical}, Harman and Snowden defined $\kk$-valued measures $\mu_t$ on $G$ parametrized by a parameter $t\in \kk$. In the general linear (parabolic) and symplectic case, the parameter $t$ is just the value of the corresponding measure $\mu_t$ on the $G$-set $\V$; in the general linear (Levi) version, the parameter $t$ is the measure of $\V_{\bar 0}$ (equivalently, of $\V_{\bar 1}$); in the unitary case, the parameter $t$ is given by $\mu_t(X_{0,1})-\mu_t(X_{1,0})=t-1$, so $\mu_t(\V)=t^2$.

The measures $\mu_t$ above exhaust all the measures for the oligomorphic pair $(G, \V)$.

To describe such a measure explicitly, one needs to explain how to compute $\mu_t(X_A)$ for $A\in \mathcal{A}_f$. This turns out to be a polynomial in $t$, as seen below. Note that the polynomial itself depends only on the finite structure $A$ and not on the measure!
\begin{lemma}[\cite{HS-classical}]\label{lem:measure_is_poly}
    For every $A\in \mathcal{A}_f$, there exists a polynomial $P_A\in \kk[x]$ such that for any $t\in \kk$, $ \mu_t(X_A)=P_A(t)$.
\end{lemma}
Let us briefly explain how these polynomials are constructed.

Fix $A\in \mathcal{A}_f$, and let $k$ denote the number of integer parameters needed to describe a finite structure in $\mathcal{A}$ (so $k=1$ for the parabolic general linear group, $k=3$ for the Levi general linear group, and $k=2$ in the remaining cases). 

By \cite{HS-classical}, one can find explicit polynomials $Q_A$ over $\Z[1/q]$ in variables $x_1, x_2, \ldots$ so that for any $M\in  \mathcal{A}_f$, we have: $$\#\Emb(A,M)=Q_A(q^a, q^b, \ldots)$$ where $a,b,\ldots$ denote the $\F_q$-dimensions of the non-degenerate and radical parts of $M$. All these polynomials are given by explicit recursions. 
\begin{remark}
    To keep our notation aligned with \cite{HS-classical} in the unitary case, we will choose $Q_A$ so that $\#\Emb(A,M)=Q_A((-q)^a, (-q)^b)$. 
\end{remark}

The polynomial $P_A$ is then obtained by setting $P_A(x)= Q_A(x, 1, \ldots, 1).$
\medskip

Finally, let us consider the case of the orthogonal group. In this case, the above statement needs to be modified: by \cite[§9]{HS-classical}, one can find explicit polynomials $Q^*_A$ over $\Z[1/q]$ in variables $x_1, x_2$ so that for any $V_{a,1,b}\in  \mathcal{A}_f$, we have: $$\#\Emb(A,V_{a,1,b})=Q_A^*({q^{\lfloor a/2\rfloor}}, q^b)$$ where $* \in \{even, ~odd\}$ is the parity of $a$. Again, these polynomials are given by explicit recursions. We will denote $P_A^*(x)=Q_A^*(x, 1)$ with $* \in \{even, ~odd\}$. In fact, if we denote by $\chi: \F_q^{\times} \to \{\pm 1\}$ the quadratic character, then $ \#\Emb(A,V_{a,\eps,b})=Q_A^*({\chi(\eps)q^{\lfloor a/2\rfloor}, q^b})$ for any $\eps\in {\F_q^{\times}/(\F_q^{\times})^2}$.

The measures for the orthogonal oligomorphic pair form a family $\{\mu^{even}_t, \mu^{odd}_t\}_{t\in \kk}$, with $\mu^*_t(X_A)=P_A^*(t)$ for any $t\in \kk$ and any $*$. Additionally, $\mu^{even}_0=\mu^{odd}_0$ (see \cite[§9.6]{HS-classical}). 

\subsection{Saturation property}\label{ssec:saturation}
The structural stabilizer class $ \mathcal{E}_{\V}=\{G_S~:~S\subset {\V} \text{ is finite} \}$ for all classical oligomorphic pairs $(G, \V)$ is saturating in the sense of \cref{def:saturating}, for every
$\kk$-valued measure. This is due to \cite[Proposition~15.2]{HS-oligomorphic} \cite[Propositions~6.8, 8.8, 9.10, and the paragraph
immediately preceding Corollary~10.10]{HS-classical}.

\subsection{Valuations} \label{sec:valuations}
Let $(G, \V)$ be a classical Levi oligomorphic pair (we exclude the parabolic $GL$ case), and $\mathcal{A}_f$ be the corresponding category of finite structures. 

For any $A\in \mathcal{A}_f$, let $\rk(A)$ denote the rank of the non-degenerate part of $A$ as before (that is, if $A=V_{a,\ldots}$ then $\rk(A)=a$). To understand the importance of the rank $
\rk(A)$ in the case $t=0$, we will use the following lemma.
\begin{lemma}\label{lem:valuation}

For any $A\in \mathcal{A}_f$,
\[
 P_A(T)=T^{\rk(A)}R_A(T),\qquad\text{ where } \qquad R_A(0)\neq0.
\]
(in the orthogonal case, we obtain $P_A^{*}(T)=T^{\rk(A)}R_A^{*}(T)$ for $*\in \{even, odd\}$).
\end{lemma}

\begin{proof} This follows from the explicit formulas for the polynomials $P_A$ given in \cite[Proposition~6.3, Proposition~8.3, Proposition~9.5, Proposition~10.4]{HS-classical}. In the orthogonal case, the computation is made separately for $P_A^{even}$ and $P_A^{odd}$ and gives the same order of $T$ in both branches.

\end{proof}

Recall that {(see \cref{sec:singular-non-zero})}, an embedding $i\in \Emb(A, B)$ induces a $G$-map $ p_i:X_B\longrightarrow X_A.$
Let $Q_{B/A,i}(T)$ be the polynomial describing the measure of the fiber of $p_i$ (that is, for any $t\in \kk$, we have $Q_{B/A,i}(t)=\mu_t(p_i)$; in the orthogonal case, we fix either the even or the odd branch). 

\begin{lemma}\label{lem:zero-fibre}
For an embedding $i\in \Emb(A, B)$
\[
 Q_{B/A,i}(0)\ne0
 \quad\Longleftrightarrow\quad
 \rk(B)=\rk(A).
\]
\end{lemma}

\begin{proof}
Recall that by Equation \eqref{eq:polynomial-ratio}, we have: $ P_B(T)=P_A(T)Q_{B/A,i}(T)
$ (after choosing a branch in the orthogonal case).

The polynomial $P_A$ is nonzero, and the above formula holds in the integral domain $\mathbb Q[T]$. Taking orders of vanishing at zero gives
\[
 \ord_0 P_B=\ord_0 Q_{B/A,i}+\ord_0 P_A.
\]
By Lemma \ref{lem:valuation} this is the same as
\[
 \rk(B)=\ord_0 Q_{{B}/A,i}+\rk(A).
\]
Thus $Q_{B/A,i}(0)$ is nonzero exactly when $\rk(B)=\rk(A)$.
\end{proof}
\subsection{Singular parameters}\label{ssec:singular_param}
The measures $\mu_t$ for the classical oligomorphic pairs $(G, \V)$ are known to be regular  for generic values of $t$.

We will say that the value $t\in \kk$ for a given classical oligomorphic pair $(G, \V)$ is \emph{singular} if $\mu_t$ is not a regular measure. The singular parameters for each classical oligomorphic pair are listed in the table below. \[
\begin{array}{c|c|c}
\text{family of oligomorphic groups}&\text{measure}&\text{singular values of the parameter } t\\ \hline
\text{general linear (parabolic)} &\mu_{t}&q^n \text{ for } n\geq 0\\
\text{general linear (Levi)} &\mu_{t}&0,~q^n \text{ for } n\geq 0\\
\text{symplectic} &\mu_{t}&0,~q^{2n} \text{ for }  n\geq 0\\
\text{unitary}&\mu_{t}&0, ~(-q)^n \text{ for }  n\geq 0\\
\text{orthogonal}&\mu^{even}_{t}&0,~ 1,~ \pm q^r \text{ for } r\geq 1\\
&\mu^{odd}_{{t}}& 0, ~\pm q^r \text{ for } r\geq 0
\end{array}
\]

%%%%%%%%%%%%%%%%%%%%

\section{Tensor ideals at singular non-zero parameters}\label{sec:singular-non-zero}
Our goal in this section is to study the poset of ideals in the permutation categories for classical oligomorphic groups. The problem is only interesting when the values of the parameter $t$ are {singular}: when $\mu_t$ is regular, the category $\Perm(G, \mu_t)$ possesses only trivial ideals, per \cref{cor:q_reg_implies_minimal_ideal}. 

Moreover, it is proved in \cite{HS-classical} that for any $t\neq 0$, the categories $\Perm(G, \mu_t)^{\Kar}$ coincide for the parabolic and the Levi version of the general linear group. In particular, it means that the posets of ideals in the corresponding permutation categories coincide, so it is enough to study the Levi version.

In this subsection, we prove the following theorem.

\begin{theorem} \label{thm:classify-classical-nonzero} Let $(G, \V)$ be a classical oligomorphic pair. Suppose that the measure parameter $t$ is singular but non-zero. Then $\Perm(G, \mu_t)$ has exactly one non-trivial ideal: the negligible ideal $\mathcal{N}$.
\end{theorem}

%We will prove this statement using Corollary \ref{cor:two-class}. 
The crucial step in the proof is that every nonzero root of the counting polynomial $P_A(T)$ is simple; this statement is wrong when $t=0$, in which case there appear additional tensor ideals as we will show in \cref{sec:ideals_at_zero}.

Before we prove \cref{thm:classify-classical-nonzero}, we provide some auxiliary results. Fix $(G, \V)$ to be a classical oligomorphic pair with corresponding structure category $\mathcal{A}$.

For any $i\in \Emb(A,B)$, precomposition with $i$ gives a $G$-map $ p_i:X_B\longrightarrow X_A$. Recall that this construction induces an anti-equivalence between $\mathcal{A}^{mono}_f $ and the category of transitive $\V$-smooth $G$-sets, via $A\mapsto X_A$, $i\mapsto p_i$. 

Let $Q_{B/A,i}(T)$ be the polynomial describing the measure of the fiber of $p_i$ (that is, for any $t\in \kk$, we have $Q_{B/A,i}(t)=\mu_t(p_i)$). Then by the Tower Formula \cref{eq:tower-formula}, we have:
\begin{equation}\label{eq:polynomial-ratio}
 P_B(T)=P_A(T)Q_{B/A,i}(T).
\end{equation}

\begin{lemma}\label{lem:simple-root-regularity}
Fix $t_0\in\kk^\times$ and let $A\in \mathcal{A}_f$. Fix the stabilizer $G_A$ of a point in $X_A$. Suppose:
\begin{enumerate}[label=\textup{(\roman*)}]
\item $P_A(t_0)=0$;
\item for every $B\in \mathcal{A}_f$, every nonzero root of the polynomial $P_B(T)$ is simple;

\end{enumerate}
Then the restricted measure $\mu_{t_0}\rvert_{G_A}$ is regular.
\end{lemma}

\begin{proof}
Since the $\V$-stabilizer class is large, it is enough to check that the following statement holds: for any $B\in \mathcal{A}_f$ with a $G$-map $X_B\to X_A$, it holds that $\mu_{t_0}(G_A/G_B)\neq 0$, where $G_B$ denotes the stabilizer of a point in $X_B$ such that $G_B\leq {G_A}$.

Such a $G$-map $X_B\to X_A$ is necessarily of the form $p_i$ for some $i\in \Emb(A, B)$. Then $ \mu_{t_0}(G_A/G_B) = Q_{B/A,i}(t_0)$. 

By \eqref{eq:polynomial-ratio} and assumption (i), the value $t_0$ is a root of both $P_A$ and $P_B$. By assumption (ii), the order of the root $t_0$ of $P_B$ is $1$ so $ Q_{B/A,i}(t_0) \neq 0$, implying that $ \mu_{t_0}(G_A/G_B) \neq 0$ as required.

\end{proof}

\begin{proposition}
    \label{thm:finite-model}
    Fix a classical oligomorphic pair $(G, \V)$ and a measure $\mu_{t_0}$ for some singular parameter $t_0\neq 0$ (for the orthogonal group, we let $\mu_{t_0}=\mu_{t_0}^*$ for some parity $*$, and in what follows, write $P_A$ instead of $P_A^*$).
    
Assume that the following conditions hold:
\begin{itemize}
    \item There exists $M\in \mathcal{A}_f$ such that $P_A(t_0) = \sharp\Emb(A, M)$ for any $A\in \mathcal{A}_f$.
\item Every nonzero root of every polynomial $P_A$ ($A\in \mathcal{A}_f$) is simple.
\end{itemize}
Then the poset $\mathcal{P}_{\mu_{t_0}}^G$ has exactly two elements
\[
 \mathsf N=\{[X_A] \ \vert \  \Emb(A, M)=\emptyset\},\qquad
 [\bullet_G]=\{[X_A]  \ \vert \ \Emb(A, M)\neq \emptyset\}.
\]
Hence the negligible ideal is the unique nonzero proper tensor ideal.
\end{proposition}

\begin{proof} We will reduce this to \cref{cor:two-class}, or rather to its relative version (see \cref{ssec:rel_setting_ideals}). It is enough to check for every $A\in \mathcal{A}_f$, one of the following conditions holds: 
\begin{itemize}
    \item $\mu_{t_0}(X_A)\neq 0$ (then $[X_A] = [\bullet_G]$).
    \item the stabilizer $G_A$ of any given point in $X_A$ satisfies: the measure $\mu_{t_0}\rvert_{G_A}$ is regular (then $[X_A] =  \mathsf N$). 
\end{itemize}
Indeed, let $A\in \mathcal{A}_f$. 

If $\sharp \Emb(A, M) \neq 0 $ then $\mu_{t_0}(X_A) = P_A(t_0) = \sharp \Emb(A, M)\neq 0$, so we are in the first case.

If $\sharp \Emb(A, M) =0$ then $P_A(t_0) = \sharp \Emb(A, M)=0$. By \cref{lem:simple-root-regularity}, in this case  the measure $\mu_{t_0}\rvert_{G_A}$ is regular, so we are in the second case.

\end{proof}

\begin{proof}[Proof of \cref{thm:classify-classical-nonzero}] 
  Fix a classical oligomorphic pair $(G, \V)$ and a measure $\mu_{t_0}$ (respectively, $\mu_{t_0}=\mu_{t_0}^*$ in the orthogonal case) for some singular parameter $t_0\neq 0$. For each $A\in \mathcal{A}_f$, let $P_A$ be the corresponding polynomial (respectively, $P_A^*$ in the orthogonal case).  We need to verify that the conditions of \cref{thm:finite-model} hold.

  To see that all the non-zero roots of $P_A$ (respectively, $P_A^*$ the orthogonal case) are simple, we use the explicit formulas for the $P_A$ from \cite[Propositions~6.3, 8.3, 10.4]{HS-classical}, and the formulas \cite[Proposition~9.5]{HS-classical} for $P_A^*$ in the orthogonal case.

The choice of $M$, based on \cref{sec:counting}, is given in the following list, depending on the group $G$, the value $t_0$ of the parameter and the choice of measure (in the orthogonal case).  This list exhausts all the singular non-zero values of the parameter $t$, per \cite{HS-classical}.

\[
\begin{array}{c|c|c|c|c}
\text{family of groups}&t_0&\text{measure}& M &\Aut_{\mathcal{A}}(M)\\ \hline
\text{general linear (Levi)} &q^n,~ n\geq 0&\mu_{t_0}&V_{n,0,0}&GL_n(\F_q)\\
\text{symplectic} &q^{2n},~ n\geq 0&\mu_{t_0}&V_{n,0}&Sp_{2n}(\F_q)\\
\text{unitary}&(-q)^n,~ n\geq 0&\mu_{t_0}&V_{n,0}&U_n(\F_q)\\
\text{orthogonal}& \pm q^r,~ r\geq 1&\mu^{even}_{t_0}& V_{2r,\delta,0}& O_{2r,\delta }(\F_q)\\
& 1&\mu^{even}_{1}& \{0_{\F_q}\}& \{1\}\\
& \pm q^r, ~r\geq 0&\mu^{odd}_{t_0}& V_{2r+1,\delta,0}& O_{2r+1,\delta }(\F_q)
\end{array}
\]

  In the orthogonal case, $\chi: \F_q^{\times} \to \{\pm 1\}$ denotes again the quadratic character, and {$\delta\in\F_q^\times/(\F_q^\times)^2$} is such that $\chi({\delta})$ is the sign of $t_0$.

\end{proof}

\begin{remark}
    It is not difficult to present the negligible ideal $\mathcal{N}$ in this case as a kernel of a SM functor from $\Perm(G, \mu_{t_0})$ to a semisimple SM category, namely the category $ \Rep_{\kk}(\Aut_{\mathcal{A}}(M))$ for the above choice of $M$. The details of this construction will be discussed in a subsequent paper. 
\end{remark}
%%%%%%%%%%%%

%%%%%%%%%%%%%%%%%%%

\section{The singular \texorpdfstring{$t=0$}{t=0} case}
Let $(G, \V)$ be a classical oligomorphic pair (not parabolic general linear), and let $\mathcal{A}_f$ be the corresponding category of finite structures. 

For $t= 0$ the measures are no longer quasi-regular \cite[Propositions~6.18, 8.13, 9.14, 10.11]{HS-classical}, and so Corollary \ref{cor:two-class} cannot be applied anymore. In this setting, we compute the poset $\mathcal{P}_{\mu_0}^G$ directly. Keeping in mind that $\mu^{even}_0=\mu^{odd}_0$ (see \cite[§9.6]{HS-classical}), we will not distinguish between the even and the odd orthogonal branches in this case. 

By \cref{lem:valuation}, \cref{lem:zero-fibre} and the equality $\mu_0({X_A})=P_A(0)$ we have the following key fact:
\begin{corollary}\label{cor:when_mu_0_is_zero}
    Let $A, B\in \mathcal{A}_f$ and $i\in \Emb(A, B)$. Then \begin{align*}
        \mu_0(X_A)\neq 0 \;&\Longleftrightarrow \; \rk(A)=0,\\
        \mu_0(p_i)\neq 0 \;&\Longleftrightarrow \; \rk(A)=\rk(B).
    \end{align*}
\end{corollary}
\subsection{Computing the poset \texorpdfstring{$\mathcal{P}_{\mu_0}^G$}{P at t=0}}

\begin{definition}\label{def:amalgam} (See \cite[Section 2.7]{HS-classical}) 
An \emph{amalgam} of $A,B\in \mathcal{A}_f$ is a diagram of embeddings in $
\mathcal{A}_f$
\[
 \begin{tikzpicture}[baseline=-2pt,node distance=18mm]
  \node (A) {$A$};
  \node (C) [right=of A] {$C$};
  \node (B) [right=of C] {$B.$};
  \draw[-{Latex}] (A)-- node[above] {$\alpha$} (C);
  \draw[-{Latex}] (B)-- node[above] {$\beta$} (C);
 \end{tikzpicture}
\] 
The amalgam is called \emph{generated} by $A$ and $B$ if $C$ is the smallest subobject of {itself} containing $\alpha(A)\cup\beta(B)$.
\end{definition}
Given $A,B\in \mathcal{A}_f$, there is an obvious notion of an isomorphism of amalgams.

\begin{lemma}\label{lem:amalgam-orbits}
Let $A, B \in \mathcal{A}_f$. We have a surjective map
\begin{align*}
    \{\text{ generated amalgam } A\xrightarrow{\alpha} C\xleftarrow{\beta}B \text{ in } \mathcal{A}_f\}/iso\;&\longtwoheadrightarrow\; \{Z \subseteq X_A \times X_B~:~ Z\in \Set{G}^{tran}\} \\
    (A\xrightarrow{\alpha} C\xleftarrow{\beta}B) \;\; &\longmapsto \;\; p(X_C)
\end{align*}
where $p:X_C\to X_A\times X_B$ is given by the projections $p_{\alpha}:X_C\to X_A$, $p_{\beta}:X_C\to X_B$.

\end{lemma}

\begin{proof}
To prove that this is a surjection, let us construct the right inverse map.

Let $(i,j)\in X_A\times X_B$, where $i\in X_A=\Emb(A, \V)$ and $j\in X_B=\Emb(B, \V)$.  Let $C\subset\mathbf V$, $C\in \mathcal{A}_f$ be the finite substructure generated by $i(A)\cup j(B)$ as a subset of $\V$. We denote by $\alpha\in \Emb(A, C)$, $ \beta\in \Emb(B, C)$ the corresponding embeddings. Clearly, $(A\xrightarrow{\alpha} C\xleftarrow{\beta}B)$ is a generated amalgam.

This defines a map $G.(i,j) \;\longmapsto\; (A\xrightarrow{\alpha} C\xleftarrow{\beta}B)$. Let us show that this is indeed the right inverse to the map in the statement of the lemma. 

Fix $(i,j)\in X_A\times X_B$ and define an amalgam $(A\xrightarrow{\alpha} C\xleftarrow{\beta}B)$ as before, as well as a map 
\[
 p:X_C\longrightarrow X_A\times X_B,
 \qquad h\longmapsto(h\alpha,h\beta).
\]
This map is clearly $G$-equivariant. The construction of $C$ as a substructure of $\V$ gives us $h_0\in X_C=\Emb(C, \V)$ such that $(i,j)=(h_0\alpha,h_0\beta)$. Thus $(i, j)\in p(X_C)$, so $p(X_C)=G.(i,j) $.

Since $C$ is generated by $\alpha(A)$ and $\beta(B)$, two embeddings $C\to\mathbf V$ whose restrictions to both subspaces agree are equal. Thus $p:X_C\to O$ is a $G$-equivariant bijection, and the two projections to $X_A$ have the same fiber measure.

\end{proof}

Recall that for every $ A\in \mathcal{A}_f$ we defined the zero-form radical subspace $A_{\rmrad}$ (see \cref{ssec:setup_classcial_grps}, \cref{sec:classical-cases}).

\begin{proposition}\label{prop:core-preorder-revised}
Let $A, B\in \mathcal{A}_f$. Fix decompositions 
\[
 A=A_{\nd}\oplus A_{\rmrad}, \qquad  B=B_{\nd}\oplus B_{\rmrad}
\]
into a nondegenerate subspace and the radical subspace, as in \cref{ssec:setup_classcial_grps}.
Then
\[
 X_A\preceq_{\mu_0}X_B
 \quad\Longleftrightarrow\quad
 \Emb(B_{\nd}, A_{\nd})\neq \emptyset.
\]
\end{proposition}

\begin{proof}
Fix $A, B\in \mathcal{A}_f$.
Recall that by \cref{lem:succeq_mu}, $X_A\preceq_{\mu_0}X_B$ if and only if there exists a $G$-orbit $ O\subseteq X_A\times X_B
$
such that the projection $O\to X_A$ has non-zero measure of fibre.

Assume first that $X_A\preceq_{\mu_0}X_B$. By
the above argument, there is a $G$-orbit
$O\subseteq X_A\times X_B$ whose projection to $X_A$ has non-zero measure of fibre. 
By \cref{lem:amalgam-orbits}, there exists a generated amalgam $(A\xrightarrow{\alpha} C\xleftarrow{\beta}B)$ such that $O=p(X_C)$ where $p:X_C\to X_A\times X_B$ is given by the projections $p_{\alpha}:X_C\to X_A$, $p_{\beta}:X_C\to X_B$. Thus $\mu_0(p_{\alpha})\neq 0$. By \cref{cor:when_mu_0_is_zero}, we get
\[
 \rk(C)=\rk(A).
\]

The images of $A_{\nd}, B_{\nd}$ in $C$ are substructures with non-degenerate form, so $$\alpha(A_{\nd})\cap C_{\rad}=\beta(B_{\nd})\cap C_{\rad}=0.$$ Since we established that $ \rk(C)=\rk(A)$, we conclude that $\alpha(A_{\nd})$ is a maximal subspace of $C$ with non-degenerate form, so we have an orthogonal decomposition $$C=\alpha(A_{\nd})\oplus C_{\rmrad}.$$ and the projection $C \to \alpha(A_{\nd})$ gives an embedding $\beta(B_{\nd}) \hookrightarrow \alpha(A_{\nd})$, since $\beta(B_{\nd})\cap C_{\rad}=0$. This gives an embedding $B_{\nd}\to A_{\nd}$.

Conversely, let $j\in \Emb(B_{\nd}, A_{\nd})$. Consider the finite structure $$C=A_{\nd}\oplus A_{\rad}\oplus B_{\rad}$$
in $\mathcal{A}_f$. By construction, $C_{\rad} = A_{\rad}\oplus B_{\rad}$ and $A_{\nd}$ is a maximal nondegenerate subspace of $C$, so
\[
 \rk(C)=\rk(A).
\]

Let $\alpha:A\to C$, $\beta: B\to C$ be the corresponding embeddings in $\mathcal{A}_f$.
Then \cref{cor:when_mu_0_is_zero} shows that the fibre of $p_{\alpha}:X_C\to X_A$ has
nonzero measure and the triple $(X_C, p_{\alpha}, p_{\beta})$ witnesses the inequality
$X_A \preceq_{\mu_0} X_B$.
\end{proof}

\begin{corollary}\label{cor:poset_P_mu_0}
    Let $A, B\in \mathcal{A}_f$.
Then
\[
 X_A\preceq_{\mu_0}X_B
 \quad\Longrightarrow\quad
 \rk(A)\geq \rk(B).
\] Moreover, for a classical pair $(G, \V)$ which is either general linear (Levi version), symplectic or unitary these conditions are equivalent: $ X_A\preceq_{\mu_0}X_B
 \quad\Longleftrightarrow\quad
 \rk(A)\geq \rk(B).$

In the case of the orthogonal group, we have: \[X_A\preceq_{\mu_0}X_B
 \quad\Longleftrightarrow\quad
     \rk(A)> \rk(B) \text{ or } \left(\rk(A)=\rk(B) \wedge\epsilon(A)=\epsilon(B)\right)\]
     where $\epsilon(A)$ stands for the modified discriminant of $A_{\nd}$, and similarly for $\epsilon(B)$.
\end{corollary}
\begin{proof}
    By \cref{prop:core-preorder-revised}, \[
 X_A\preceq_{\mu_0}X_B \quad\Longleftrightarrow\quad
 \Emb(B_{\nd}, A_{\nd})\neq \emptyset
 \quad\Longrightarrow\quad
 \rk(A)\geq \rk(B).
\]  The structure of $\mathcal{A}_f$ in the general linear (Levi version), symplectic and unitary settings immediately gives $ \Emb(B_{\nd}, A_{\nd})\neq \emptyset
 \quad\Longleftrightarrow\quad
 \rk(A)\geq \rk(B)$.
 
 Let us prove the $(\Rightarrow)$ implication in the orthogonal case. Let $A, B\in \mathcal{A}_f$ and assume $\Emb(B_{\nd}, A_{\nd})\neq \emptyset$. Then $ \rk(A)\geq \rk(B)$; if there is an equality $ \rk(A)= \rk(B)$ then the embedding $ B_{\nd}\hookrightarrow A_{\nd}$ is an isomorphism so $\epsilon(A)=\epsilon(B)$.

 To prove the $(\Leftarrow)$ implication in the orthogonal case, let $A, B\in \mathcal{A}_f$.
 
 First, assume that $\rk(A)=\rk(B) $ and $\epsilon(A)=\epsilon(B)$. Then $A_{\nd}\cong B_{\nd}$ in $\mathcal{A}_f$, so $ \Emb(B_{\nd}, A_{\nd})\neq \emptyset$.

 Now, assume that $\rk(A)>\rk(B) $, and let $a=\rk(A)-\rk(B)$. We may choose an $a$-dimensional space $V_{a,\overline{\epsilon(A)}, 0} \in \mathcal{A}_f$ with a non-degenerate bilinear form and a suitable modified discriminant $\overline{\epsilon(A)}$ such that $B_{\nd}\oplus  V_{a,\overline{\epsilon(A)}, 0}$ has discriminant $\epsilon(A)$. Thus $ B_{\nd}\oplus  V_{a,\overline{\epsilon(A)}, 0}\cong A_{\nd}$ in $\mathcal{A}_f$ and so $ \Emb(B_{\nd}, A_{\nd})\neq \emptyset$, as required.
\end{proof}
\begin{corollary}\label{cor:equiv_classes_mu_0}
   For a classical pair $(G, \V)$ which is either general linear (Levi version), symplectic or unitary, we have: $$ X_A~\sim_{\mu_0}~X_B
 \quad\Longleftrightarrow\quad
 \rk(A)= \rk(B).$$

In the case of the orthogonal group, we have: \[X_A~\sim_{\mu_0}~X_B
 \quad\Longleftrightarrow\quad \left(\rk(A)=\rk(B) \wedge\epsilon(A)=\epsilon(B)\right)\]
     where $\epsilon(A)$ stands for the modified discriminant of $A_{\nd}$, and similarly for $\epsilon(B)$ (if $\rk(A)=\rk(B)=0$, the condition $\epsilon(A)=\epsilon(B)$ should be considered void).
\end{corollary}
\subsection{Tensor ideals at \texorpdfstring{$t=0$}{t=0}}\label{sec:ideals_at_zero}

Let $(G, \V)$ be a classical oligomorphic pair of either general linear (Levi version), symplectic or unitary types. 

Recall that each class in $\mathcal{P}_{\mu_0}^G$ contains some $X_A$ due to \cref{ssec:saturation}.  
For any $[X]\in \mathcal{P}_{\mu_0}^G$ such that $[X]=[X_A]$, we will write $\rk([X])=\rk(A)$ and $\eps([X])=\eps(A)$ (if $\rk([X])>0$) in the orthogonal case. These are well-defined (independent of the choice of $A$) due to \cref{cor:equiv_classes_mu_0}.

\begin{theorem}\label{thm:zero-chain}
Let $(G, \V)$ be a classical oligomorphic pair of either general linear (Levi version), symplectic or unitary types. The poset $\mathcal{P}_{\mu_0}^G$ is a linear poset
\[
 \mathsf C_0>\mathsf C_1>\mathsf C_2>\cdots,
\]
where
\[
 \mathsf C_a=[X_A] \text{ with } \rk(A)=a
\] (this is well-defined due to \cref{cor:equiv_classes_mu_0}).
The elements of $\Down(\mathcal{P}_{{\mu_0}}^G)$ are
\[
 \varnothing,\qquad D_m=\{\mathsf C_a \ | \ a\geq m\}\quad(m\geq0).
\]
Consequently the tensor ideals are $0$ and $\mathcal I_m$ ($m\geq0$), where 
\[
 \mathcal I_m(\mathcal{C}(X),\mathcal{C}(Y))
 =\Span_\kk\{\one_O \ : \ O\in G\backslash (Y\times X) , ~  \rk([O])\geq m\}.
\]
\end{theorem}
\begin{remark}
    The ideal $\mathcal I_0$ is the entire category $\Perm(G, \mu_0)$, while $\mathcal I_1=\N$ is the negligible ideal.
\end{remark}

\begin{proof}
The description of the poset $\mathcal{P}_{{\mu_0}}^G$ follows directly from \cref{cor:poset_P_mu_0} and \cref{cor:equiv_classes_mu_0}. The explicit description of the tensor ideals corresponding to each downward-closed subset $D_m$ follows from \cref{thm:classification}.

\end{proof}

We now describe the ideals in $\Perm(G, \mu_0)$ in the orthogonal case. Let $(G, \V)$ be the orthogonal oligomorphic pair and $\mathcal{A}_f$ the corresponding category of finite structures. %For $t=0$, the two branches of the measures on $G$ coincide.

  As before, $\chi: \F_q^{\times} \to \{\pm 1\}$ denotes the quadratic character. For any $A\in \mathcal{A}_f$ with $\rk(A)>0$ and for any $[X]\in \mathcal{P}^G_{\mu}$ with $\rk(X)>0$ we write
\[
 \sigma(A)=\chi(\epsilon(A))~\in~\{+,-\}, \;  \sigma(X)=\chi(\epsilon(X))~\in~\{+,-\}.
\]

\begin{theorem}\label{thm:orthogonal-zero}
    Set $\mathsf C_0=[X_A]$ where $ \rk(A)=0.$ For $a\geq1$ and $\sigma\in\{+,-\}$ set $$\mathsf C_{a,\sigma} =[X_A] \text{ where } \rk(A)=a, \; \sigma(A)=\sigma
$$ (this is well-defined due to \cref{cor:equiv_classes_mu_0}). The poset $\mathcal{P}_{\mu_0}^G$ is $$\mathcal{P}_{\mu_0}^G=\{\mathsf C_0, \mathsf C_{a,\sigma}~|~a\geq 1, ~\sigma\in\{+,-\}\}$$
with partial order given by
\[
 \mathsf C_{a,\sigma}\preceq_{\mu_0} \mathsf C_{b,\tau}
 \quad\Longleftrightarrow\quad
 a>b\ \text{ or }\ (a=b\text{ and }\sigma=\tau) 
\] and $\mathsf C_0$ being the maximal element.  Thus the two classes $\mathsf C_{a,\pm}$ for every $a\geq 1$ are incomparable, and $C_{a+1, \sigma}\preceq_{\mu_0} \mathsf C_{a, \pm}$ for any $\sigma$.

The elements of $\Down(\mathcal{P}_{{\mu_0}}^G)$ are
$
 \varnothing,~ \mathcal{P}_{\mu_0}^G,
$
and, for each $m\geq1$ and each non-empty $ S\subseteq\{+,-\}$,
\[
 D_{m,S}
 =\{\mathsf C_{m,\sigma}:\sigma\in S\}
  \cup\{\mathsf C_{a,\sigma}:a>m,\ \sigma\in\{+,-\}\}.
\]
The corresponding ideals are
\[
 \mathcal I_{m,S}(\mathcal{C}(X),\mathcal{C}(Y))
 =\Span_\kk\left\{\one_O\,:\, O\in  G\backslash (Y\times X) \, \text{ and } \,
 \begin{cases}
      \rk([O])>m,\ &\text{or}\\
 \rk([O])=m&\text{ and }\sigma(O)\in S
 \end{cases}
\right\}.
\]

\end{theorem}
\begin{proof}

The description of the poset $\mathcal{P}_{{\mu_0}}^G$ follows directly from \cref{cor:poset_P_mu_0} and \cref{cor:equiv_classes_mu_0}. The explicit description of the tensor ideals corresponding to each downward-closed subset $D_{m, S}$ follows from \cref{thm:classification}.

\end{proof}
\begin{remark}
  With these notations we have 
\[
 \mathcal I_{1,\{+,-\}}=\N,
 \qquad
 \mathcal I_{m,\{+\}}+\mathcal I_{m,\{-\}}=\mathcal I_{m,\{+,-\}},
 \qquad
 \mathcal I_{m,\{+\}}\cap\mathcal I_{m,\{-\}}=\mathcal I_{m+1,\{+,-\}}.
\]
\end{remark}

In the orthogonal case, the Hasse diagram (with larger ideals placed higher) is an infinite string of diamonds:

\begin{center}
\scalebox{0.7}{\begin{tikzpicture}[
  node distance=11mm and 22mm,
  every node/.style={font=\small},
  line/.style={thick}
]
\node (top) {$\mathcal{I}_0 = \Perm(G=\mathrm{O}, \mu_0)$};
\node (b1) [below=of top] {$\mathcal I_{1,\{+,-\}}=\N$};
\node (a1p) [below left=of b1] {$\mathcal I_{1,\{+\}}$};
\node (a1m) [below right=of b1] {$\mathcal I_{1,\{-\}}$};
\node (b2) [below=18mm of b1] {$\mathcal I_{2,\{+,-\}}$};
\node (a2p) [below left=of b2] {$\mathcal I_{2,\{+\}}$};
\node (a2m) [below right=of b2] {$\mathcal I_{2,\{-\}}$};
\node (b3) [below=18mm of b2] {$\mathcal I_{3,\{+,-\}}$};
\node (dots) [below=of b3] {$\vdots$};
%\node (zero) [below=of dots] {$0$};

\draw[line] (top)--(b1);
\draw[line] (b1)--(a1p);
\draw[line] (b1)--(a1m);
\draw[line] (a1p)--(b2);
\draw[line] (a1m)--(b2);
\draw[line] (b2)--(a2p);
\draw[line] (b2)--(a2m);
\draw[line] (a2p)--(b3);
\draw[line] (a2m)--(b3);
\draw[line] (b3)--(dots);
%\draw[line] (dots)--(zero);
\end{tikzpicture}}
\end{center}

In this diagram, the two prongs at level $m$ remember the two discriminant classes at the smallest rank present in the lower set.

The elementary explanation in terms of linear algebra is that for fixed positive dimension there are two nondegenerate quadratic spaces, distinguished by discriminant, and the two types are incomparable.  In a larger dimension, either type embeds into the other type: choose the discriminant of a nondegenerate complement so that the orthogonal sum has the prescribed target discriminant.

%%%%%%

%%%%%%%%

\subsection{Tensor ideals of abelian type at \texorpdfstring{$t=0$}{t=0}}

Our aim in this section is to give explicit descriptions of tensor ideals for $t=0$ as kernels of SM functors into abelian rigid SM categories. Tensor ideals which can be presented as such kernels are said to be \emph{of abelian type}.

Let $q$ be a prime power, and let $(G, \V)$ be a classical oligomorphic pair (not parabolic).
We will denote $\F=\F_q$ when in the general linear (Levi), symplectic and orthogonal cases, while $\F=\F_{q^2}$ in the unitary case.

\subsubsection{Semisimplification}
We start with the negligible ideal $\mathcal{N}$ in $\Perm(G, \mu_0)$. We will construct an explicit equivalence between $\Perm(G, \mu_0)^{\Kar}/\mathcal{N}$ and the semisimple SM category presented below.

For a finite field $\F$ of size $q'$, let
\[
 \mathscr G_{\F}=\Perm(GL^{par}_{\infty}(\F), \mu_0)^{\Kar}
\]
be the Karoubi envelope of the permutation category for the parabolic classical oligomorphic pair 
$ (GL_{\infty}(\F), \F^{\infty})$ and measure $\mu_0$. This category is denoted by $\mathcal T(\mathrm{Vect}_{\F_{q'}},\delta_0)$ in \cite{Knop} and by $\underline{\Rep}\bigl(GL_0(\F_{q'})\bigr)$ in \cite{EAHGL}. Since the parameter of the measure is $0\notin\{{(q')}^n~:~n\geq0\}$, the category $\mathscr G_{\F}$ is abelian semisimple in characteristic zero \cite{Knop, HS-classical}.

For our chosen classical oligomorphic pair $(G, \V)$, define
\[
 \mathscr S_G=
 \begin{cases}
 \mathscr G_{\F}\boxt\mathscr G_{\F},&\text{ if } G=GL,\\
 \mathscr G_\F,&\text{ if }G=Sp, O \text{ or } U.
 \end{cases}
\]

\begin{theorem}\label{thm:rankzero}
Let $q$ be a prime power, and let $(G, \V)$ be a classical oligomorphic pair (not parabolic){,} with $ \mathscr S_G$ defined as above. There is a SM equivalence
\[
 \overline{Q}_G: \mathscr S_G~\xrightarrow{\sim}~\mathcal \Perm(G, \mu_0)^{\Kar}/\N.
\]
\end{theorem}
The key ingredient is the universal property of $\mathscr S_G$, proved in \cite{EAHGL} and described briefly in \cref{sssec:F_q_lin}, allowing us to define the functor $\overline{Q}_G$. The proof that $\overline{Q}_G$ is an equivalence will be given in \cref{sssec:proof_of_semisimplification}.

\subsubsection{\texorpdfstring{$\F$}{F}-linear Frobenius spaces and the universal property of \texorpdfstring{$\mathscr S_G$}{the category S}}\label{sssec:F_q_lin}

The universal property of $ \mathscr G_{\F}$, proved in \cite{EAHGL}, states that this is the universal Karoubi additive $\kk$-linear SM category generated by a $\F$-linear Frobenius object of dimension $0$ (\cite[Definition~5.1.1]{EAHGL}). Such an object should be regarded as the algebra of functions on a finite scheme over $\F$, the scheme having the structure of {a} linear $\F$-space, and the algebra of functions should have categorical dimension $0$. We denote by $\underline{\V}$ the generating object of $ \mathscr G_{\F}$.

Recall that for a finitary smooth $G$-set $X$, the corresponding object
$\cC(X)$ carries the structure of a special commutative Frobenius algebra.

Let $V$ be an $\F$-vector space regarded with an $\F$-linear finitary smooth $G$-action. The maps
\[
 +:V\times V\longrightarrow V,
 \qquad
 0:\{*\}\longrightarrow V,
 \qquad
 m_a:V\longrightarrow V,\quad w\longmapsto aw
\]
are then $G$-equivariant.

Recall that for $X, Y\in \Set{G}$, a $G$-map $f:X\to Y$ defines a $G$-morphism $f_*:\cC(X)\to \cC(Y)$. Set
\[
 \dotplus=(+)_*:
 \cC(V)\otimes\cC(V)\longrightarrow\cC(V),
 \qquad
 z=(0)_*:\one\longrightarrow\cC(V),
\]
\[
(m_a)_*:\cC(V)\longrightarrow\cC(V)
 \qquad(a\in\F).
\]

\begin{lemma}\label{lem:linear-frob}
With the maps $\dotplus,z,(m_a)_*$ from above, $\cC(V)$ is an $\F$-linear Frobenius space in the sense of \cite[Definition~5.1.1]{EAHGL}.
\end{lemma}

\begin{proof}
    The Schwartz space $\cC(V)$ is always a special commutative Frobenius algebra by \cite{HS-oligomorphic}. All the defining relations described in \cite[Definition~5.1.1]{EAHGL} follow by applying either the $(-)_*$ construction to the corresponding identities of maps of the vector space $V$ or follow from graph-map functoriality together with the relevant diagonal and base-change identities.
\end{proof}

Recall our choice of the $\F$-linear $G$-space $\V$ as defined in \cref{sec:classical-cases}.

\begin{lemma}\label{lem:existence_of_func_Q_G}
    Let $q$ be a prime power, and let $(G, \V)$ be a classical oligomorphic pair (not parabolic) with $ \mathscr S_G$ defined as above. There exists a SM functor
\[
 Q_G: \mathscr S_G~\longrightarrow~\mathcal \Perm(G, \mu_0)^{\Kar}
\] 

\end{lemma}

\begin{proof}
By \cref{lem:linear-frob} and by the universal property of $\mathscr G_{\F}$ (see \cite[Theorem~8.1.1]{EAHGL}), a SM functor from $\mathscr G_{\F}$ to $\Perm(G, \mu_0)^{\Kar}$ is determined by a choice of an $\F$-linear $G$-set so that $\dim \cC(V)=\mu_0(V)$ would be $0$. This will be the image of the generating object $\underline{\V}$ of $ \mathscr G_{\F}$.

So in the symplectic, orthogonal and unitary case, we take $V=\V=\F^{\infty}$ and the Schwartz space $W=\cC(\mathbf V)$ is the required $\F$-linear Frobenius space of categorical dimension $0$. 

In the general linear (Levi) case, where $\mathbf V=\mathbf V_{\bar0}\oplus\mathbf V_{\bar1}$, we consider the two Schwartz spaces $\cC(\mathbf V_{\bar0})$, $\cC(\mathbf V_{\bar1})$. Both are $\F$-linear Frobenius spaces of categorical dimension $0$, so we obtain a functor 
\[
 Q_G: \mathscr S_G = \mathscr G_{\F}\boxtimes \mathscr G_{\F}~\longrightarrow~\mathcal \Perm(G, \mu_0)^{\Kar}
\] 
sending $\underline{\V}\boxtimes \triv $ to $\cC(\mathbf V_{\bar0})$ and $ \triv \boxtimes \underline{\V}$ to $\cC(\mathbf V_{\bar1})$.
\end{proof}
\begin{definition}\label{def:overline_Q_G}
    We define $\overline{Q}_G$ as the composition of $Q_G$ with the quotient (semisimplification) functor 
$$\Perm(G, \mu_0)^{\Kar} ~\longrightarrow~\Perm(G, \mu_0)^{\Kar}/\N.$$
\end{definition}

%%%%%%%%%%%%%%%%%%%%%%

%%%%%%%%%%%%%%%%%%%%%%%%%

%%%%%%%%%%%%%%%%%%%%%%%

%%%%%%%%%%%%%%%%%%%%%%%%

\subsubsection{Proof of the Semisimplification \cref{thm:rankzero}}\label{sssec:proof_of_semisimplification}

\begin{proof}[Proof of Theorem \ref{thm:rankzero}]

Consider the functor $\overline{Q}_G: \mathscr{S}_G\to \Perm(G, \mu_0)^{\Kar}/\N$ defined in \cref{def:overline_Q_G}. This is a $\kk$-linear SM functor so its kernel is a tensor ideal in $\mathscr{S}_G$. Yet $\mathscr{G}_{\F}$ and $\mathscr{S}_G$ are semisimple rigid SM categories with $\End(\triv)=\kk$, so they have no non-trivial ideals. Thus $\overline{Q}_G$ is faithful.

To prove that $\overline{Q}_G$ is full, it is enough to {prove} the fullness of $ \overline{Q}_G$ restricted to $\mathscr{S}_G^0$: this is the full subcategory of $\mathscr{S}_G$ whose objects are tensor powers of the generator (for $G=\mathrm{GL}$, we take the mixed tensor powers of the two generators). Passing to the additive Karoubi envelope $ \mathscr{S}_G=(\mathscr{S}_G^0)^{\Kar}$ will preserve the fullness of the functor.

To prove fullness on the above subcategory, we need to establish that the corresponding $\Hom$-spaces have the same dimension. 

In the category $ \mathscr{G}_{\F}$, the $\Hom$-space $\Hom_{\mathscr{G}_{\F}}(\underline{\V}^{\otimes m},
\underline{\V}^{\otimes n})$ has a basis parametrized by $\mathrm{R}_{m+n}=\{\F\text{-linear subspaces } R\subseteq \F^{m+n}\}$. Respectively, in $ \mathscr{G}_{\F}\boxtimes\mathscr{G}_{\F}$, the $\Hom$-space $$\Hom_{\mathscr{G}_{\F}\boxtimes\mathscr{G}_{\F}}(\underline{\V}^{\otimes m_0}\boxtimes\underline{\V}^{\otimes m_1},
\underline{\V}^{\otimes n_0} \boxtimes \underline{\V}^{\otimes n_1}) = \Hom_{\mathscr{G}_{\F}}(\underline{\V}^{\otimes m_0}, \underline{\V}^{\otimes n_0})\otimes_{\kk}
 \Hom_{\mathscr{G}_{\F}}( \underline{\V}^{\otimes m_1}, \underline{\V}^{\otimes n_1}) $$ has a basis parametrized by $$\mathrm{R}_{m_0+n_0}\times \mathrm{R}_{m_1+n_1}=\{(R_0, R_1)~:~ R_i\subseteq \F^{m_i+n_i} \text{ an } \F\text{-linear subspace }\}.$$

The corresponding $\Hom$-spaces in $ \Perm(G, \mu_0)^{\Kar}/\N$ are:
\begin{itemize}[itemindent=0pt,leftmargin=*]
    \item In the symplectic, orthogonal or unitary case:
   \begin{align*}
       \Hom_{ \Perm(G, \mu_0)/\N}(\overline{Q}_G(\underline{\V})^{\otimes m},\overline{Q}_G(\underline{\V})^{\otimes n}) &= \Hom_{ \Perm(G, \mu_0)/\N}(\cC(\V)^{\otimes m},\cC(\V)^{\otimes n})
   \end{align*}
   which has, by \cref{lem:negl-basis}, a basis given by $\left\{O\in  G\backslash\V^{m+n} ~:~~\mu_0(O)\neq 0\right\}.$
   In \cref{lem:rank-zero-orbits-formed} below, we will prove that the latter set is in bijection with $\mathrm{R}_{m+n}$.
       \item In the general linear case (Levi version):
   \begin{align*}
       &\Hom_{ \Perm(G, \mu_0)/\N}(\overline{Q}_G(\underline{\V}^{\otimes m_0}\boxtimes\underline{\V}^{\otimes m_1}),\overline{Q}_G(\underline{\V}^{\otimes n_0} \boxtimes \underline{\V}^{\otimes n_1})) \\
       &= \Hom_{ \Perm(G, \mu_0)/\N}(\cC(\V_{\bar 0})^{\otimes m_0}\otimes \cC(\V_{\bar 1})^{\otimes m_1},\cC(\V_{\bar 0})^{\otimes n_0}\otimes \cC(\V_{\bar 1})^{\otimes n_1})\\
       &\cong  {\Hom_{ \Perm(G, \mu_0)/\N}(\triv,\cC(\V_{\bar 0})^{\otimes (m_0+n_0)}\otimes\cC(\V_{\bar 1})^{\otimes (m_1+ n_1)})}
   \end{align*}
   which has, by \cref{lem:negl-basis}, a basis given by $\left\{O\in G\backslash( \V_{\bar 0}^{m_0+n_0}\times \V_{\bar 1}^{m_1+n_1} )~:~\mu_0(O)\neq 0\right\}.$
   In \cref{lem:rank-zero-orbits-levi} below, we will prove that the latter set is in bijection with $\mathrm{R}_{m_0+n_0}\times \mathrm{R}_{m_1+n_1}$.
\end{itemize}
This proves the fullness of $\overline{Q}_G$. 

We now show that $\overline{Q}_G$ is essentially surjective. 
Recall that by \cref{ssec:saturation}, the stabilizer class $\mathcal{E}_{\V}$ is saturating. By \cref{cor:tensor-generator}, this implies that $\Perm(G, \mu_0)^{\Kar}$ is generated (as a Karoubi additive monoidal category) by $\cC(\V)$. Since $\cC(\V)$ lies in the essential image of $\overline{Q}_G$ in each case, we conclude that $\overline{Q}_G$ is essentially surjective.

\end{proof}

\begin{remark}

In \cite{Knop, EAHGL}, an explicit basis $\{f_R~:~ R\in \mathrm{R}_{m+n}\}$ of the $\Hom$-space $\Hom_{{\mathscr{G}_{\F}}}(\underline{\V}^{\otimes m},
\underline{\V}^{\otimes n})$ was constructed. 

It is not difficult to write out the action of $\overline{Q}_G$ on the corresponding basis of the appropriate $\Hom$-space in $\mathscr{S}_G$. For $G=\mathrm{Sp},\mathrm O,\mathrm U$ and $R\subseteq \F^{m+n}$, one has in $\Perm(G, \mu_0)^{\Kar}/\N$
\begin{equation*}\label{eq:zeta-formula}
 \overline{Q}_G(f_R)=\sum_{S~:\,R\subseteq S\subseteq \F^{m+n}}\one_{O_S},
\end{equation*}
where $ \one_{O_S} \in \cC(\V^{m+n})$ is the indicator function on the orbit $O_S$ constructed in \cref{lem:rank-zero-orbits-formed}. 

For $G=\mathrm{GL}$ (Levi version) and $R_i\subseteq \F^{m_i+n_i}$ ($i=0,1$), one has in $\Perm(G, \mu_0)^{\Kar}/\N$
\begin{equation*}\label{eq:zeta-formula-product}
 \overline{Q}_G(f_{R_0}\boxt f_{R_1})
=\sum_{(S_0,~ S_1) ~:~ R_i\subseteq S_i\subseteq \F^{m_i+n_i}}
 \one_{O_{S_0,S_1}}
\end{equation*}
where $ \one_{O_{S_0,S_1}} \in \cC(\V_{\bar 0}^{ m_0+n_0}\times \V_{\bar 1}^{m_1+n_1})$ is the indicator function on the orbit $O_{S_0,S_1}$ constructed in \cref{lem:rank-zero-orbits-levi}.

\end{remark}

Before we proceed to the auxiliary lemmata, we remind the reader of the crucial criterion of \cref{cor:when_mu_0_is_zero}: for $A\in \mathcal{A}_f$,
\begin{equation}
 \mu_0(X_A)\neq 0
 \quad\Longleftrightarrow\quad
 \rk(A)=0 \quad\Longleftrightarrow\quad A=A_{\rad}.
\end{equation}

\begin{definition}
    For a tuple $\mathbf{v} = (v_1,\ldots,v_d) \in \V^d$, define its \emph{space of linear relations} by \[ \mathrm{Rel}(\mathbf{v}) = \{ (a_1,\ldots,a_d) \in \F^d \ | \ \sum_i a_i v_i = 0 \}. \] Equivalently, let $\alpha_{\mathbf{v}} :\F^d\to \mathbf{V}, \ \alpha_{\mathbf{v}} (a_1,\ldots,a_d)  = \sum a_i v_i.$ Then $\mathrm{Rel}(\mathbf{v}) = \ker(\alpha_{\mathbf{v}})$.
\end{definition}

Let $d\geq 0$. Recall the notation $\mathrm{R}_d= \{\F\text{-linear subspaces } R\subseteq \F^{d}\}$.
%%%%%%

%%%%%%%%%%%%

\begin{lemma}[The $\mathrm{Sp, O, U}$ case] \label{lem:rank-zero-orbits-formed}
Assume that $(G, \V)$ is either the symplectic, the orthogonal or the unitary oligomorphic classical pair.  
Let $d\geq 0$. For every linear subspace $R\subseteq \F^d$, define
\begin{equation}\label{eq:OR-explicit}
 O_R=
 \left\{
 \mathbf v=(v_1,\ldots,v_d)\in\mathbf V^d \ | \
 \Rel(\mathbf v)=R,
 \ \langle v_i,v_j\rangle=0\text{ for all }i,j
 \right\}.
\end{equation}
Then $O_R \in \Set{G}^{tran}$, ${\mu_0}(O_R)\neq 0$ and the correspondence
$$ \mathrm{R}_{d}\;\longrightarrow\;\left\{O\in G\backslash \V^{d} ~:~\mu_0(O)\neq 0\right\}, \qquad R~\longmapsto ~O_R$$
is a bijection.

\end{lemma}

\begin{proof} 
For any $\mathbf v=(v_1,\ldots,v_d)\in\mathbf V^d$, let $A_{\mathbf v} = \Span_{\F}\{v_1,\ldots,v_d\} \subseteq \mathbf{V}$ with the restricted alternating, symmetric or Hermitian form respectively.

Let $R\in \mathrm{R}_d$. Let us first show that the set $O_R$ defined in \eqref{eq:OR-explicit} is not empty. Equip the quotient space $\F^d/R$ with the zero form. We choose a totally isotropic embedding $ j:\F^d/R~\hookrightarrow~\mathbf V.$ 
 Let $(e_i)_{i=1}^d$ denote the standard basis of $\F^d$. Putting $v_i=j(e_i+R)$ gives a tuple $\mathbf v=(v_1,\ldots,v_d)\in\mathbf V^d$ whose relation space is $R$ and which has $\langle v_i,v_j\rangle=0$ for all $i,j$.

Now let us show that $O_R$ is transitive (it is clearly $G$-stable). Suppose $\mathbf v,\mathbf w$ both satisfy the condition \eqref{eq:OR-explicit} for the same $R\in \mathrm{R}_d$. This determines a linear isomorphism $\phi:A_{\mathbf v}\to A_{\mathbf w}$ in $\mathcal{A}_f$ sending $v_i$ to $w_i$. Any isomorphism in $\mathcal{A}_f$ extends to an automorphism of $\mathbf{V}$. Hence there exists $g\in G$ with $g\mathbf v=\mathbf w$, so $O_R \in \Set{G}^{tran}$.

Since the tuple $\mathbf v$ spans $A_{\mathbf v}$, the map
\[
 X_{A_{\mathbf v}}=\Emb(A_{\mathbf v},\mathbf V)
 \longrightarrow G\mathbf v,
 \qquad
 i\longmapsto(i(v_1),\ldots,i(v_d))
\]
is a $G$-equivariant bijection. Thus $\mu_0( X_{A_{\mathbf v}})=\mu_0(G\mathbf v)$. \cref{cor:when_mu_0_is_zero} mentioned above says that $\mu_0( X_{A_{\mathbf v}}){\neq} 0$ precisely when the restricted form on $A_{\mathbf v}$ has rank zero - that is, when it vanishes identically. Since the $v_i$ span $A_{\mathbf v}$, this condition translates to $\langle v_i,v_j\rangle=0$ for all $i,j$.

This proves that the correspondence $R~\mapsto~ O_R$ is well-defined. It is clearly injective, so it remains to prove that it is surjective.

Given $\mathbf v=(v_1,\ldots,v_d)\in  \V^{d}$ such that $\mu_0(G\mathbf v)\neq 0$, let $R=\Rel(\mathbf v)$. To prove that the correspondence is surjective, we need to show that $G\mathbf v=O_R$; that is, $\langle v_i,v_j\rangle=0$ for all $i,j$. The above argument shows that $\mu_0(G\mathbf v)\neq 0$ if and only if $\langle v_i,v_j\rangle=0$ for all $i,j$, so we are done.

\end{proof}
    In the general linear (Levi) case the arguments are the same, just adapted to $\mathbb{Z}/2$-graded spaces.

\begin{lemma}[The $\mathrm{GL}$ case]\label{lem:rank-zero-orbits-levi}

Assume that $(G, \V)$ is the general linear (Levi version) oligomorphic classical pair.  
Let $d_0, d_1\geq 0$. For every pair of linear subspaces $(R_0, R_1)$ with $R_i\subseteq \F^{d_i}$, define
\begin{align}\label{eq:OR-levi-explicit}
 O_{R_0,R_1}=\bigl\{(\mathbf x,\mathbf y)&\in
 \mathbf V_{\bar 0}^{d_0}\times\mathbf V_{\bar 1}^{d_1} \ | \ \\
 &\Rel(\mathbf x)=R_0,\quad \Rel(\mathbf y)=R_1,
 \quad \langle x_i,y_j\rangle=0\text{ for all }i,j\bigr\}.\notag
\end{align}
Then $ O_{R_0,R_1} \in \Set{G}^{tran}$, ${\mu_0}( O_{R_0,R_1})\neq 0$ and the correspondence
$$ {\mathrm{R}_{d_0}\times \mathrm{R}_{d_1}}\;\longrightarrow\;\left\{O\in G\backslash \V_{\bar 0}^{d_0}\times\mathbf V_{\bar 1}^{d_1} ~:~\mu_0(O)\neq 0\right\}, \qquad (R_0, R_1)~\longmapsto ~ O_{R_0,R_1}$$
is a bijection.

%Then $O_{R_0,R_1}$ is a single $G$-orbit, every pair $(R_0,R_1)$ occurs, and these are exactly the non-zero-measure orbits at $t=0$.
\end{lemma}

\begin{proof}
For any $\mathbf x=(x_1,\ldots,x_{d_0})\in\V_{\bar 0}^{d_0}$, $\mathbf y=(y_1,\ldots,y_{d_1})\in\V_{\bar 1}^{d_1}$, let $$A_{\mathbf x,\mathbf y}=A_{0}\oplus A_{1}\subseteq \V,
 \qquad
 A_{ 0}=\Span\{x_i\},\quad A_{ 1}=\Span\{y_j\},$$ with the restricted pairing.

Let $(R_0, R_1)\in \mathrm{R}_{d_0}\times \mathrm{R}_{d_1}$. We first show that the set ${O_{R_0, R_1}}$ defined in {\eqref{eq:OR-levi-explicit}} is not empty. Equip the quotient space $\F^{d_0}/R_0\oplus \F^{d_1}/R_1$ with the zero pairing. We choose a $\Z/2$-graded embedding $ j:\F^{d_0}/R_0\oplus \F^{d_1}/R_1~\hookrightarrow~\mathbf V$ respecting the pairing. 
 Let $(e^{(r)}_i)_{i=1}^{{d_r}}$ denote the standard basis of $\F^{d_r}$, $r=0,1$. The images of the standard bases in $\F^{d_0}, \F^{d_1}$ give tuples $\mathbf x=(x_1,\ldots,x_{d_0})\in\mathbf V_0^{d_0}$, $\mathbf y=(y_1,\ldots,y_{d_1})\in\mathbf V_1^{d_1}$ whose relation spaces are $R_0, R_1$ respectively, and which have $\langle x_i,y_j\rangle=0$ for all $i,j$.

Now let us show that ${O_{R_0, R_1}}$ is transitive (it is clearly $G$-stable). Suppose $(\mathbf x,\mathbf y),(\mathbf x',\mathbf y')$ both satisfy the condition \eqref{eq:OR-levi-explicit} for the same pair $(R_0, R_1)\in \mathrm{R}_{d_0}\times \mathrm{R}_{d_1}$. This determines a linear isomorphism $\phi:A_{\mathbf x,\mathbf y}\to A_{\mathbf x',\mathbf y'}$ in $\mathcal{A}_f$ sending $x_i$ to $x'_i$, ${y_j}$ to ${y'_j}$. Any isomorphism in $\mathcal{A}_f$ extends to an automorphism of $\mathbf{V}$, so there exists $g\in G$ with $g(\mathbf x,\mathbf y)=(\mathbf x',\mathbf y')$. Thus ${O_{R_0, R_1}} \in \Set{G}^{tran}$.

Since the tuples $ \mathbf x,\mathbf y$ span $A_0, A_1$ respectively, the map
\[
 X_{A_{\mathbf x,\mathbf y}}=\Emb(A_{\mathbf x,\mathbf y},\mathbf V)
 \longrightarrow G.(\mathbf x,\mathbf y),
 \qquad
 i\longmapsto{(i(x_1),\ldots,i(x_{d_0}),i(y_1),\ldots,i(y_{d_1}))}
\]
is a $G$-equivariant bijection. Thus $\mu_0( X_{A_{\mathbf x,\mathbf y}})=\mu_0(G.(\mathbf x,\mathbf y))$. \cref{cor:when_mu_0_is_zero} mentioned above says that $\mu_0( X_{A_{\mathbf x,\mathbf y}}){\neq} 0$ precisely when the restricted form on $A_{\mathbf x,\mathbf y}$ has rank zero - that is, when it vanishes identically. Since $ \mathbf x,\mathbf y$ span $A_0, A_1$ respectively, this condition translates to $\langle x_i,y_j\rangle=0$ for all $i,j$.

This proves that the correspondence $(R_0, R_1)~\mapsto~ O_{(R_0, R_1)}$ is well-defined. It is clearly injective, and the equivalence 
$$\mu_0({G.(\mathbf x,\mathbf y)})\neq 0 \qquad \Longleftrightarrow \qquad \forall ~i, j, \;\langle {x_i,y_j}\rangle=0 $$
proved above again shows that the correspondence $(R_0, R_1)~\mapsto~ O_{(R_0, R_1)}$ is surjective.

\end{proof}

%%%%%%%%%%%%%%%%%%%%%%%%

\subsubsection{Prime ideals and ideals of abelian type}\label{sec:abelian_ideals}

\begin{definition} 
\mbox{}

\begin{itemize}[itemindent=0pt,leftmargin=*]
    \item A tensor ideal in $\mathcal{T}$ is \emph{of abelian type} (see \cite[Definition 2.1.1]{CEO-abelian}) if it is the kernel of a tensor functor to an abelian finite-length rigid SM category with $\End(\one) = \kk$.
\item A proper thick tensor ideal $\I$ is \emph{prime} if $X\otimes Y \in \I$ implies that either $X \in\I$ or $Y \in \I$ (see \cite{Balmer}, \cite[Section 1.2.1]{CEO-abelian}).   
\end{itemize} 
\end{definition}

It was proven in \cite[Proposition 2.1.3]{CEO-abelian} that an ideal of abelian type is always \emph{faithfully prime}, and its associated thick tensor ideal is prime: this follows from the fact that in an abelian rigid SM category with $\End(\one) = \kk$, tensoring with a non-zero object (either on the left or on the right) is a faithful endofunctor.

\medskip

Let $A=A_{\nd} \in \mathcal{A}_f$, and let $U_A\leq G$ be its pointwise stabilizer. Since $A=A_{\nd}$, $U_A\cong G$ as oligomorphic groups, so the measures on $U_A$ are the ones described in \cref{sec:counting}. Restriction from $G$ to $U_A$ yields a SM functor
\[
 \Res_A:\Perm(G, \mu_0)^{\Kar}\longrightarrow \Perm(U_A,\mu_0|_{U_A})^{\Kar}.
\]
By the restriction formulas in \cite{HS-classical}, the measure $\mu_0|_{U_A}$ is the measure on $U_A$ corresponding to the parameter $t=0$ per \cite[§§6.7, 8.5, 9.6, 10.4, respectively]{HS-classical}. We define the SM functor $\Phi_A$ as a composition of $\Res_A$ with the quotient functor of \cref{thm:rankzero}: 
\[
 \Phi_A:\Perm(G, \mu_0)^{\Kar} ~\xrightarrow{\;\Res_A\;}~
 \Perm(U_A,\mu_0|_{U_A})^{\Kar}~\rightarrow~
 \Perm(U_A,\mu_0|_{U_A})^{\Kar}/\mathcal{N}~\xrightarrow{\;\overline{Q}_{U_A}^{~-1}\;} ~\mathscr S_{U_A}\cong \mathscr S_{G}.
\]

\begin{theorem}\label{thm:kernels}
\mbox{}

\begin{enumerate}[itemindent=0pt,leftmargin=*,label=\textup{(\alph*)}]
\item For $G=\mathrm{GL,~Sp,~U}$, let $A=A_{\nd}\in \mathcal{A}_f$ with $\rk(A)=r\geq0$.
Then
\[
 \ker(\Phi_{A})=\mathcal I_{r+1}.
\]
\item For $G=\mathrm{O}$, let $A=A_{\nd}\in \mathcal{A}_f$ with $\rk(A)=r\geq1$ and modified discriminant $\epsilon(A)$. If $\overline{\epsilon(A)}$ denotes
the other class in $\F_q^\times/(\F_q^\times)^2$ and $\chi: \F_q^{\times} \to \{\pm 1\}$ the quadratic character, then
\[
 \ker(\Phi_{A})=\mathcal I_{r,\{\chi(\overline{\epsilon(A)})\}}.
\]
For the non-degenerate space $A=A_{\nd}\in \mathcal{A}_f$ with $\rk(A)=0$, the corresponding functor $\Phi_{A}$ is the semisimplification functor $\Perm(G, \mu_0)^{\Kar} \to \mathscr S_{G}$, and $ \ker(\Phi_{A})=\mathcal{N}=\mathcal{I}_{1, \{+, -\}}$.
\end{enumerate}
\end{theorem}

The main input is the following lemma.

\begin{lemma}\label{lem:restriction-detects}
Let $A=A_{\nd}\in \mathcal{A}_f$. For every $B\in \mathcal{A}_f$,
\begin{equation}\label{eq:restriction-detects}
 \Phi_A\bigl(\mathcal{C}(X_B)\bigr)\neq0
 \quad\Longleftrightarrow\quad
 X_A\preceq_{\mu_0}X_B  \quad\Longleftrightarrow\quad \Emb(B_{\nd}, A)\neq \emptyset.
\end{equation}

\end{lemma}

\begin{proof}
The restriction of $X_B$ to $U_A$ is a finite disjoint union of $U_A$-orbits,
$
 X_B=Y_1\sqcup\cdots\sqcup Y_s.$ 
Hence
\[
 \Res_A\mathcal{C}(X_B)\cong\bigoplus_{i=1}^s\mathcal{C}(Y_i).
\]
\cref{ex:negligible} implies that $$\forall ~Y\in \Set{U_A}^{tran}, \qquad \id_{\mathcal{C}(Y)}\in \mathcal{N} \qquad \Longleftrightarrow \qquad\mu_0(\Delta Y)=0$$ where $\Delta Y\subset Y\times Y$ is the diagonal and $\mathcal{N}$ is the ideal of negligible morphisms in $\Perm(U_A, \mu_0\rvert_{U_A})$. If $\pi: \Perm(U_A,\mu_0|_{U_A})^{\Kar}~\rightarrow~
 \Perm(U_A,\mu_0|_{U_A})^{\Kar}/\mathcal{N}$ denotes the semisimplification functor, applying  $\overline{Q}_{U_A}^{~-1}$, we get:
\[
 \overline{Q}_{U_A}^{~-1}(\pi(\mathcal{C}(Y_i)))=0
 \quad\Longleftrightarrow\quad
 (\mu_0|_{{U_A}})(Y_i)=0.
\]
It follows that
\begin{equation}\label{eq:some-orbit}
 \Phi_A\bigl(\mathcal{C}(X_B)\bigr)\neq0
 \quad\Longleftrightarrow\quad
 \text{some $U_A$-orbit in $X_B$ has nonzero $\mu_0$-measure.}
\end{equation}

Let $i_A\in \Emb(A, \mathbf{V}) = X_A$ be the fixed embedding such that $U_A = \Stab_G(i_A)$. If $i_B\in X_B$, the $G$-orbit
\[
 \mathcal O=G\cdot(i_A,i_B)\subseteq X_A\times X_B
\]
has fibre over $i_A$ equal to
\[
 \{i_B'\in X_B:(i_A,i_B')\in\mathcal O\}=U_A\cdot i_B.
\]
Thus the $U_A$-orbits in $X_B$ are the fibres over $i_A$ of the diagonal $G$-orbits in $X_A\times X_B$, and their restricted measures are the corresponding fibre measures.  By \cref{lem:succeq_mu}, the right hand side of \eqref{eq:some-orbit} is equivalent to $X_A\preceq_{\mu_0}X_B$. This proves \eqref{eq:restriction-detects}. Since $A_{\nd}=A$, \cref{prop:core-preorder-revised} gives: $$X_A\preceq_{\mu_0}X_B \qquad \Longleftrightarrow \qquad \Emb(B_{\nd}, A)\neq \emptyset.$$
\end{proof}

\begin{proof}[Proof of \cref{thm:kernels}] Let $\mathcal K_A=\ker(\Phi_A)$.  The functor $\Phi_A$ is SM, so $\mathcal K_A$ is a tensor ideal. Since all the tensor ideals in $\Perm(G, \mu_0)$ and its Karoubi envelope are thick (\cref{cor:tensor_ideals_are_thick}), the ideal $\mathcal K_A$ is thick, so it is determined by the objects $\mathcal{C}(X_B)$ for which $\Phi_A(\mathcal{C}(X_B))=0$. By Lemma~\ref{lem:restriction-detects},
\begin{equation}\label{eq:killed-core}
 \Phi_A(\mathcal{C}(X_B))=0
 \quad\Longleftrightarrow\quad
 \Emb(B_{\nd}, A)=\emptyset.
\end{equation}

Consider first the case when $G=\mathrm{GL,~Sp,~U}$ and $A=A_{\nd}\in \mathcal{A}_f$ with $\rk(A)=r$.  In each of these three families, there is a unique isometry class of a non-degenerate structure $A$ with given rank $r$, so all such $A$ are isomorphic. In this case \eqref{eq:killed-core} simply means
\[
 \Phi_A(\mathcal{C}(X_B))=0
 \quad\Longleftrightarrow\quad
 \rk(B)\geq r+1.
\]

\cref{thm:zero-chain} states that the downward subset $D_{r+1}=\{[X_B]~:~\rk(B)\geq r+1\}$ of $\mathcal{P}^G_{\mu_0}$ corresponds to the tensor ideal $\I_{r+1}$. Therefore $\mathcal K_{A}=\I_{r+1}$ for $\rk(A)=r$ (in particular, if $\rk(A)=0$ then $\ker(\Phi_{A})=\I_1=\N$).

Now let $G=\mathrm{O}$ and $A=A_{\nd}\in \mathcal{A}_f$ with $\rk(A)=r{\geq 1}$ and modified discriminant $\epsilon(A)$.  By \cref{prop:core-preorder-revised}, \cref{cor:poset_P_mu_0}, for any $B\in \mathcal{A}_f$ with modified discriminant $\epsilon(B)$, we have: 
$$\Emb({B_{\nd}}, A)\neq \emptyset \qquad \Longleftrightarrow \qquad \rk(B)<\rk(A) \text{ or } (\rk(B)=\rk(A)\,\wedge \, \epsilon(B)=\epsilon(A)).$$
So \eqref{eq:killed-core} says that $\Phi_{A}$ kills exactly the classes
\[
 \{\mathsf C_{r,\chi(\overline{\epsilon(A)})}\}
 \cup
 \{\mathsf C_{s,\delta}:s>r,\ \delta\in\{+,-\}\}
\]
($\mathsf C_{s,\delta}$ as defined in \cref{thm:orthogonal-zero}). This is the downset $D_{r,\{\chi(\overline{\epsilon(A)})\}}$ of Theorem \ref{thm:orthogonal-zero}.  By definition, its tensor ideal is $\I_{{r,\{\chi(\overline{\epsilon(A)})\}}}$, and so  $\ker(\Phi_{A})=\I_{{r,\{\chi(\overline{\epsilon(A)})\}}}$. If $A=A_{\nd}\in \mathcal{A}_f$ with $\rk(A)=0$, then $A=0$ so $U_A=G$ and $\Res_A$ is the identity functor. Therefore $\ker (\Phi_A)=\N$.
\end{proof}

\begin{corollary}\label{cor:abelian-nonzero}

\mbox{}

\begin{itemize}[itemindent=0pt,leftmargin=*]
    \item In the general linear (Levi), symplectic and unitary cases, every tensor ideal $\mathcal I_m$ of $\Perm(G, \mu_0)$ for $m\geq1$ is of abelian type.  
    \item In the orthogonal case, $\mathcal I_{1,\{+,-\}}=\N$ and every $\mathcal I_{m,\{\sigma\}}$ ($\sigma\in \{\pm\}$) are of abelian type. In particular, all their associated thick tensor ideals are prime. 

The thick tensor ideals $\mathcal{I}_{m,\{+,-\}}$ for $m\geq 2$ are not prime and therefore the associated tensor ideals are not of abelian type. 
\end{itemize}

\end{corollary}

\begin{proof} The statement about ideals of abelian type follows from \cref{thm:kernels}. By \cite{CEO-abelian} if $\mathcal I$ is of abelian type then $\mathcal{I}$ is prime. So it remains to show that the ideals $\mathcal I_{m,\{+,-\}}$ in the orthogonal case are not prime for $m \geq 2$.

We now consider the ideal $\mathcal{I}_{m,\{+,-\}}$ in the orthogonal case. Recall that $$\I_{m-1,\{+\}} \cap \I_{m-1,\{-\}} = \I_{m,\{+,-\}}.$$ Since these are tensor ideals, the tensor product of two objects in $\I_{m-1,\{+\}}, \I_{m-1,\{-\}}$ lies in $\I_{m,\{+,-\}}$. So it is enough to find, for each $\sigma\in \{\pm\}$, an object in each $\I_{m-1,\{\sigma\}}$ which does not lie in $\I_{m,\{+,-\}}$. By the definition of the ideals $\I_{m-1,\{\sigma\}}$, $\I_{m,\{+,-\}}$, an example of such an object would be $\cC(X_A)$ where $A=A_{\nd}\in \mathcal{A}_f$ has rank $\rk(A)=m-1$ and a modified discriminant corresponding to the sign $\sigma$ under the character $\chi:\F_q^\times \to \{\pm 1\}$.

\end{proof}

\begin{remark}
    The zero ideal is a prime ideal in $\Perm(G, \mu)^{\Kar}$ for any oligomorphic group $G$ and any measure $\mu$. 
    
    Indeed, let $X,Y\in \Set{G}$ such that $ \cC(X),\cC(Y) \neq 0$, and let $e_X\in \End(\cC(X)), ~e_Y\in \End(\cC(Y))$ be non-zero idempotents. Then $X,Y\neq \varnothing$ and $e_X, e_Y$ are given by Schwartz functions $\lambda_X\in \cC(X\times X)^G,~ \lambda_Y\in \cC(Y\times Y)^G$. Now, $X\times Y\neq \varnothing$ and $e_X\otimes e_Y$ corresponds to the Schwartz function $\lambda_X\times \lambda_Y \in \cC({(X\times Y)\times(X\times Y)})$. Since $\lambda_X\times \lambda_Y$ is not identically zero, we conclude that $e_X\otimes e_Y\neq 0$. Thus $(\cC(X), e_X)\otimes (\cC(Y), e_Y)\neq 0$. This proves that the zero ideal is a prime ideal.
\end{remark}

\subsection{Generalized negligible ideals}\label{sec:higher-negligible}

We give another interpretation of the tensor ideals at the parameter $t=0$. These arise as generalized negligible ideals in the sense of \cite{HW}.

\subsubsection{Generalized negligible morphisms}

In \cite{HW} the authors construct generalizations of the negligible morphisms arising from a deformation of a given SM category over some local ring. Let $(R,\mfrak)$ be a commutative local ring.  Take a rigid spherical monoidal $R$-linear category $\mathcal{T}_R$ such
that $\End_{\mathcal{T}_R}(\one)=R$ and all Hom modules are free over $R$. Its mod-$\mfrak$
evaluation $\mathcal{T}$ is the category over $\kk = R/\mfrak$ with the same objects and
\[
 \Hom_{{\mathcal{T}}}(X,Y)=
 \Hom_{\mathcal{T}_R}(X,Y)/\mfrak\Hom_{\mathcal{T}_R}(X,Y).
\]

For an ideal $J\subseteq R$, define the collection of \emph{$J$-negligible morphisms} to be
\[
 \N^R_J(X,Y)=
 \{f:X\to Y \ | \ \Tr(gf)\in J\text{ for every }g:Y\to X\}.
\]
The $J$-negligible morphisms form a tensor ideal and the $J$-negligible
objects (objects whose identity morphism is $J$-negligible) form a thick tensor ideal \cite[Lemma~2.2]{HW}. The image of $\N^R_J$ in the residue category $\mathcal{T}$ will
be denoted by $\N_J$.

For $J=\mfrak^k$, write $\N^R_k$ for the corresponding tensor ideal and the associated thick tensor ideal. Their images in the residue category $\mathcal{T}$ will
be denoted by $\N_k$, and they define tensor ideals and thick tensor ideals respectively. For $k=1$ we recover the usual tensor ideal of negligible morphisms in $\mathcal{T}$.

Intuitively, these tensor ideals measure the vanishing order of the trace or dimension when going from $R$ to $R/\mfrak$. For example, if $\mathrm{trace}(f) \in \mfrak^k \backslash \mfrak^{k+1}$, then this means that the trace vanishes to the $k$-th order. For more examples see \cite{HW,HW2}.

\subsubsection{Oligomorphic permutation categories over local rings}

Let $R$ be an arbitrary commutative ring. Let $G$ {be an} oligomorphic group and fix an $R$-valued
measure $\widetilde\mu$ on $G$. The permutation category $\Perm_R(G,\widetilde\mu)$ is an additive rigid SM $R$-linear
category {whose} objects are $\cC_R(X)$, with $X$ a finitary smooth $G$-set.
For finitary $X,Y$, we have an isomorphism of free $R$-modules
\[
 \Hom\bigl(\cC_R(X),\cC_R(Y)\bigr)
 =\bigoplus_{O\in G\backslash(Y\times X)}R\,\one_O
\]
where $\one_O$ is the indicator function on the $G$-orbit $O\subset Y\times X$.
Furthermore, for a ring homomorphism $R\to R'$ and the specialized measure $\widetilde\mu'$, extension of scalars gives a SM equivalence
\[
 \Perm_R(G,\widetilde\mu)\otimes_R R'
 \simeq
 \Perm_{R'}(G,\widetilde\mu').
\]

The next corollary follows from the lifting conditions in \cite{HW2}.

\begin{corollary}\label{cor:kar-reduction}
Let $(R,\mfrak,\kk)$ be complete local, and let $\widetilde\mu$ be an
$R$-valued measure reducing to $\mu$. Then the {mod-$\mathfrak{m}$} evaluation of $ \Perm_R(G,\widetilde\mu)^{\Kar}$ is $\Perm_{\kk}(G,\mu)^{\Kar}$.
\end{corollary}

\begin{theorem}\label{thm:orbit J}
Let $(R,\mathfrak m)$ be a commutative local ring with residue field $\kk$,
let $\widetilde\mu$ be an $R$-valued measure, and let $J\subseteq R$ be an
ideal. For finitary smooth $G$-sets $X,Y$, one has
\begin{align}
 \mathcal N^R_J({\cC_R(X),\cC_R(Y)})
 &=\bigoplus_{O\in G\backslash(Y\times X)}
       (J:\widetilde\mu(O))\,\one_O,
       \label{equation:orbit}\\
 \mathcal N_J(\mathcal{C}(X),\mathcal{C}(Y))
 &=\Span_\kk\{\one_O \ | \ \widetilde\mu(O)\in J\}.
       \label{equation:image-orbits}
\end{align}
Here $(J:a)=\{c\in R:ca\in J\}$. 

If $R$ is a discrete valuation ring
with uniformizer $\pi$, write $v_\pi(0)=+\infty$. Then, for $r\geq1$,
\begin{equation}\label{equation:valuation}
 \mathcal N_{{r}}(\mathcal{C}(X),\mathcal{C}(Y))
 =\Span_\kk\{\one_O \ | \, O\in G\backslash (Y\times X), \ v_\pi(\widetilde\mu(O))\geq r\}.
\end{equation}
\end{theorem}

\begin{proof}
By Corollary \ref{cor:trace-diagonal} for orbits $O,O'$ in $Y \times X$ \[ \operatorname{Tr}\bigl(\one_{O'^T}\circ\one_O\bigr)
 =\begin{cases}
      \widetilde\mu(O),&O'=O,\\
      0,&O'\ne O.
   \end{cases}\]
    This calculation is also valid in $R$. Therefore for
$f=\sum_O a_O\one_O$ and $g=\sum_O b_O\one_{O^T}$,
\[
 \operatorname{Tr}(gf)=\sum_O a_Ob_O\widetilde\mu(O).
\]
Testing separately with each $\one_{O^T}$ proves \eqref{equation:orbit}. 
If $\widetilde\mu(O)\in J$, the coefficient ideal
$(J:\widetilde\mu(O))$ is all of $R$, so its image under $R\to \kk$ is $\kk$.
In the other case $(J:\widetilde\mu(O))$ is contained in $\mathfrak m$ and its image in $\kk$ is therefore zero. Reducing
\eqref{equation:orbit} mod $ \mathfrak{m}$ proves \eqref{equation:image-orbits}.

For a discrete valuation ring, an element is in $(\pi^r)$ if and only if its
valuation is at least $r$. This proves \eqref{equation:valuation}. More explicitly,
if $s_O=v_\pi(\widetilde\mu(O))<\infty$, the corresponding summand in the lift 
over $R$ is
\[
 \pi^{\max\{r-s_O,0\}}R\,\one_O.
\]
If on the other hand $\widetilde\mu(O)=0$, the corresponding summand over $R$ is $R\,\one_O$.

\end{proof}

For our given $\Perm$-categories over $\kk$ attached to the classical oligomorphic pairs $(G, \V)$, we now choose $R$ as follows: we set 
\[
 R=\kk[[h]],\qquad\mathfrak m=(h),\qquad R/\mathfrak m=\kk.
\] Equivalently, $R$ is the completion of $\kk[T]_{(T)}$, with $T$ identified
with $h$. We also need to have a lift of our fixed measure $\mu_0$. For this it is convenient to use the universal measure ring $\Theta(G)$, whose homomorphisms
to a commutative ring $R'$ classify $R'$-valued measures on $G$
\cite[Section~3]{HS-oligomorphic}. For $G=\mathrm{GL},~\mathrm{Sp},~\mathrm U$, the polynomial measures are induced by homomorphisms
\[
 \phi_G:\Theta(G)\otimes_{\mathbb Z}\mathbb Q
             \longrightarrow\mathbb Q[T]
\]
with $\phi_G([X_A])=P_A(T)$ \cite[Corollaries~6.17, 8.12 and~10.10]{HS-classical}.
Define $\widetilde\mu$ by the composite
\begin{equation}\label{sec68:measure}
 \Theta(G)\longrightarrow\Theta(G)\otimes_{\mathbb Z}\mathbb Q
 \xrightarrow{\phi_G}\mathbb Q[T]
 \longrightarrow {\kk}[[h]],\qquad T\longmapsto h.
\end{equation}
With this definition we have $\widetilde\mu(X_A)=P_A(h)$, and its reduction is $\mu_0$. 

For $G= \mathrm{O}$ \cite[Theorem~9.12 and Corollary~9.13]{HS-classical} identifies $\Theta(G)\otimes_{\mathbb Z}\mathbb Q\cong \mathbb Q[u,v]/(uv)$. We can then pass to $R$ either via  $(u,v)\mapsto(h,0)$ or via $(u,v)\mapsto(0,h)$. Both choices reduce to the same measure at $(u,v)=(0,0)$ (the measures of the $X_A$ are then given either by $P_A^{\mathrm{even}}(h)$ or 
$P_A^{\mathrm{odd}}(h)$).

\begin{theorem}\label{thm:classical-Nk}
In the general linear (Levi), symplectic and unitary cases, for every $k\geq1$ we have: $\mathcal{N}_k=\I_k$.
In the orthogonal case, for every $k\geq1$ we have: $\mathcal{N}_k=\I_{k,\{+,-\}}$.
\end{theorem}

\begin{proof} Recall that any $Z \in \Set{G}^{tran}$ has a finite cover $X_A \to Z$; let $d$ be the size of its fibre. Then $\tilde{\mu}(X_A)=d\cdot\tilde{\mu}(Z)$ and $d\in \kk^{\times}$. Denoting by $v_h$ the $h$-adic valuation, we get: $$ v_h(\tilde{\mu}(Z)) =v_h(\tilde{\mu}(X_A))= \rk(A).$$ Therefore, taking into account the description of ideals $\mathcal N_k$ in Theorem \ref{thm:orbit J}, we get  
\[
 \mathcal N_k(\mathcal{C}(X),\mathcal{C}(Y))
 =\Span_{\kk}\{\one_O \ | \ \rk([O])\geq k\}.
\]
For $G=\mathrm{GL},\mathrm{Sp},\mathrm U$, this is the description of $\mathcal I_{{k}}$. For $\mathrm O$, this is the description of 
$\mathcal I_{{k},\{+,-\}}$.
\end{proof}

\begin{remark} The definition of the $\mathcal{N}_{{J}}$ is not canonical: different deformations over different local rings may produce different tensor ideals $\mathcal{N}_{{J}}$. This should be seen as an advantage as it makes the concept more flexible. For example, one can show that replacing $R$ by $R' = \kk[[u,v]]/(uv)$, and taking suitable ideals $J'$ allows to interpret also the ideals $\I_{r,\{+\}}$ and $\I_{r,\{-\}}$ as generalized negligible ideals. 
\end{remark}
%%%%%%%%%%%%%%%%%%%%%%%%%%%

%%%%%%%%%%%%%%%%%%%%%%%%%%%%

\printbibliography
\end{document}